\documentclass[psamsfonts,reqno]{amsart}
\usepackage{amscd,amssymb}
\newcommand{\jj}{\vee}% join
\newcommand{\mm}{\wedge}% meet
\newcommand{\JJ}{\bigvee}% big join
\newcommand{\MM}{\bigwedge}% big meet
\newcommand{\uu}{\cup}% union
\newcommand{\ii}{\cap}% intersection
\newcommand{\UU}{\bigcup}% big union
\newcommand{\ci}{\subseteq}% contained in with equality
\newcommand{\sci}{\subset}% strictly contained in with equality
\newcommand{\ce}{\supseteq}% containing with equality
\newcommand{\es}{\varnothing}% the empty set
\newcommand{\set}[1]{\{#1\}}% set
\newcommand{\setm}[2]{\{\,#1\mid#2\,\}}% set with a middle
\def\vv<#1>{\langle#1\rangle}% vector

\newcommand{\ga}{\alpha}
\newcommand{\gb}{\beta}

\newcommand{\gf}{\varphi}
\renewcommand{\gg}{\gamma}% old use >>
\newcommand{\gh}{\eta}

\newcommand{\gk}{\kappa}
\newcommand{\gl}{\lambda}
\newcommand{\gm}{\mu}
\newcommand{\gn}{\nu}
\newcommand{\go}{\omega}
\newcommand{\gp}{\pi}
\newcommand{\gq}{\theta}
\newcommand{\gr}{\varrho}

\newcommand{\gx}{\xi}
\newcommand{\gy}{\psi}
\newcommand{\gz}{\zeta}

\newcommand{\gF}{\Phi}

\newcommand{\gO}{\Omega}

\newcommand{\gQ}{\Theta}

\newcommand{\tbf}{\textbf}% text bold
\newcommand{\tup}{\textup}% text upright

\newcommand{\E}[1]{\mathcal{#1}}% same as \C
\newcommand{\ol}[1]{\overline{#1}}

\providecommand{\bysame}{\makebox[3em]{\hrulefill}\thinspace}
\newcommand{\q}{\quad}% spacing
\newcommand{\qq}{\qquad}% more spacing
\newcommand{\iso}{\cong}% isomorphic
\def\con#1=#2(#3){#1\equiv#2\pod{#3}}
\theoremstyle{plain}

\newtheorem{lemma}{Lemma}[section]
\newtheorem{theorem}{Theorem}
\newtheorem{proposition}[lemma]{Proposition}
\newtheorem{corollary}[lemma]{Corollary}
\newtheorem{claim}{Claim}

\newtheorem*{stat}{\name}
\newcommand{\name}{testing}

\theoremstyle{definition}
\newtheorem{definition}[lemma]{Definition}
\newtheorem{example}[lemma]{Example}
\newtheorem{problem}{Problem}

\theoremstyle{remark}
\newtheorem{remark}[lemma]{Remark}
\newtheorem{notation}[lemma]{Notation}
\newtheorem*{note}{Note}

\newenvironment{all}[1]{\renewcommand{\name}{#1}\begin{stat}}
                        {\end{stat}}

\newcommand{\qedc}{{\qed}~{\rm Claim~{\theclaim}.}}

\newenvironment{cproof}
{\begin{proof}[Proof of Claim.]}
{\qedc\renewcommand{\qed}{}\end{proof}}

\numberwithin{equation}{section}

\newcommand{\oll}[1]{\,\overline{\!#1}}
\newcommand{\rd}{\mathrm{d}}
\newcommand{\zero}{\mathbf{0}}
\newcommand{\one}{\mathbf{1}}
\newcommand{\two}{\mathbf{2}}
\newcommand{\gos}{\go\setminus\set{0}}
\newcommand{\peq}{\preceq}
\newcommand{\dnw}{\mathbin{\downarrow}}
\newcommand{\upw}{\mathbin{\uparrow}}

\DeclareMathOperator{\Id}{Id}
\DeclareMathOperator{\Fil}{Fil}

\newcommand{\MPL}[1]{$#1$-meas\-ured par\-tial lat\-tice}
\newcommand{\ML}[1]{$#1$-meas\-ured lat\-tice}
\newcommand{\CMPL}[1]{$#1$-co\-meas\-ured par\-tial lat\-tice}
\newcommand{\ghost}{{}_{--}}
\newcommand{\ijj}{\jj^{\Id}}
\newcommand{\fjj}{\jj^{\Fil}}

\newcommand{\IDM}{$(\mathrm{Id}\ii)$}
\newcommand{\FIM}{$(\mathrm{Fil}\ii)$}
\newcommand{\IDJ}{$(\mathrm{Id}\jj)$}
\newcommand{\FIJ}{$(\mathrm{Fil}\jj)$}

\newcommand{\Set}[1]{\left\{\,#1\,\right\}}

\newcommand{\famm}[2]{\left\langle\,#1\mid#2\,\right\rangle}
\newcommand{\fJJ}[1]{\sideset{}{_{#1}^{\mathrm{fin}}}{\JJ}}
\newcommand{\fJJe}{\sideset{}{^{\mathrm{fin}}}{\JJ}}

\newcommand{\bv}[1]{\left\|#1\right\|}
\newcommand{\vbv}[1]{{[\![}#1{]\!]}}

\def\cls(#1,#2){{{#1}{}_{/{#2}}}}

\DeclareMathOperator{\dom}{dom}
\DeclareMathOperator{\Con}{Con}
\DeclareMathOperator{\Conc}{Con_c}
\DeclareMathOperator{\W}{W}
\DeclareMathOperator{\hgt}{ht}

\newcommand{\jz}{$\set{\jj,0}$}
\newcommand{\jzh}{\jz-ho\-mo\-mor\-phism}
\newcommand{\res}{\mathbin{\restriction}}
\newcommand{\fin}[1]{[#1]^{<\go}}
\newcommand{\fine}[1]{[#1]_*^{<\go}}

\newcommand{\ID}{\E{I}}
\newcommand{\FIL}{\E{F}}
\newcommand{\tP}{\widetilde{P}}

\newcommand{\tR}{\widetilde{R}}
\newcommand{\aID}{\ID_{\mathrm{aff}}}
\newcommand{\aFIL}{\FIL_{\mathrm{aff}}}

\DeclareMathOperator{\Fg}{F}

\newcommand{\FL}{\Fg_{\mathbf{L}}}

\newcommand{\xa}{\boldsymbol{a}}
\newcommand{\xb}{\boldsymbol{b}}
\newcommand{\xc}{\boldsymbol{c}}

\newcommand{\xX}{\boldsymbol{X}}

\newcommand{\xZ}{\boldsymbol{Z}}

\newcommand{\dx}{\dot{x}}
\newcommand{\dy}{\dot{y}}
\newcommand{\dz}{\dot{z}}

\DeclareMathOperator{\Dim}{Dim}

\author[F. Wehrung]{Friedrich Wehrung}
\address{CNRS, FRE 2271\\
D\'epartement de Math\'ematiques, BP 5186\\
Universit\'e de Caen, Campus 2\\
14032 Caen cedex\\
France}
\email{wehrung@math.unicaen.fr}
\urladdr{http://www.math.unicaen.fr/\~{}wehrung}

\date{\today}

\begin{document}

\title{Forcing extensions of partial lattices}

\begin{abstract}
We prove the following result:

\emph{Let $K$ be a lattice, let $D$ be a distributive lattice with zero, and
let\linebreak
$\gf\colon\Conc K\to D$ be a \jzh, where
$\Conc K$ denotes the \jz-semilattice of all finitely generated congruences of
$K$. Then there are a lattice~$L$, a
lattice homomorphism $f\colon K\to L$, and an isomorphism
$\ga\colon\Conc L\to D$ such that $\ga\circ\Conc f=\gf$.
}

Furthermore, $L$ and $f$ satisfy many additional properties, for example:

\begin{enumerate}
\item $L$ is relatively complemented.

\item $L$ has definable principal congruences.

\item If the range of~$\gf$ is cofinal in $D$, then the convex sublattice of
$L$ generated by $f[K]$ equals $L$.
\end{enumerate}

We mention the following corollaries, that extend many results obtained
in the last decades in that area:

\begin{itemize}
\em
\item[---] Every lattice $K$ such that $\Conc K$ is a lattice admits a
congruence-pre\-serv\-ing extension into a relatively complemented lattice.

\item[---] Every \jz-direct limit of a countable sequence of distributive
lattices with zero is isomorphic to the semilattice of compact congruences of
a relatively complemented lattice with zero.
\end{itemize}
\end{abstract}

\maketitle
\tableofcontents

\section*{Introduction}

\subsection*{Background}
The Congruence Lattice Problem (CLP), formulated by R.P. Dilworth in the
forties, asks whether every distributive \jz-semilattice is isomorphic to
the semilattice $\Conc L$ of all compact congruences of a lattice $L$.
Despite considerable work in this area, this problem is still open, see
\cite{GrScC} for a survey.

In \cite{Schm68}, E.T. Schmidt presents an important sufficient condition,
for a distributive \jz-semilattice $S$, to be isomorphic to the congruence
lattice of a lattice. This condition reads ``$S$ is the image of a
generalized Boolean lattice under a `distributive' \jzh''. As
an important consequence, Schmidt proves the following result:

\begin{theorem}[Schmidt, see \cite{Schm81}]\label{T:Schmidt}
Let $S$ be a distributive lattice with zero. Then there exists a lattice $L$
such that $\Conc L\iso S$.
\end{theorem}

This result is improved in \cite{Pudl85}, where P. Pudl\'ak proves that one
can take $L$ a direct limit of finite atomistic lattices. Although we will
not use this fact, we observe that A.P. Huhn proved in
\cite{Huhn89a,Huhn89b} that Schmidt's condition is also satisfied by every
distributive \jz-semilattice $S$ such that $|S|\leq\aleph_1$.

The basic statement of CLP can be modified by keeping among the assumptions
the distributive \jz-semilattice $S$, but by adding to them a \emph{diagram}
$\E{D}$ of lattices and a morphism (in the categorical sense) from the image
of $\E{D}$ under the $\Conc$ functor to $S$. (See the end of
Section~\ref{S:PpPl} for a precise definition of this functor.)
The new problem asks whether
one can lift the corresponding \emph{diagram} $\Conc\E{D}\to S$ by a diagram
$\E{D}\to L$, for some lattice $L$
(that may be restricted to a given class of lattices)
and lattice homomorphisms. We cite a few examples:

\begin{theorem}[Gr\"atzer and Schmidt, see \cite{GrSc99}]\label{T:GrScFin}
Let $K$ be a lattice. If the lattice $\Con K$ of all congruences of~$K$ is
finite, then $K$ embeds congruence-preservingly into a sectionally
complemented lattice.
\end{theorem}

(A lattice $L$ with zero is \emph{sectionally complemented}, if for all
$a\leq b$ in $L$, there exists $x\in L$ such that $a\mm x=0$ and $a\jj x=b$.)

Theorem~\ref{T:GrScFin} does not extend to the case where $\Con K$ is
infinite: by M.~Plo\v{s}\v{c}ica, J.~T\r{u}ma, and F.~Wehrung \cite{PTWe98},
the free lattice $\FL(\go_2)$ on $\aleph_2$ generators does not have a
congruence-preserving, sectionally complemented extension. In fact, it is
proved in J.~T\r{u}ma and F.~Wehrung \cite{TuWe1} that $\FL(\go_2)$ does not
embed congruence-preservingly into a lattice with \emph{permutable}
congruences.

\begin{theorem}[Gr\"atzer, Lakser, and Wehrung, see \cite{GrLW};
see also T\r{u}ma~\cite{Tuma98}]
\label{T:GLWe}
Let $S$ be a finite distributive \jz-semilattice, let $\E{D}$ be a diagram of
lattices and lattice homomorphisms consisting of lattices $K_0$, $K_1$, and
$K_2$, and lattice homomorphisms $f_l\colon K_0\to K_l$, for
$l\in\set{1,2}$. Then any morphism from
$\Conc\E{D}$ to $S$ can be lifted, with respect to the $\Conc$ functor, by a
commutative square of lattices and lattice homomorphisms that extends $\E{D}$.
\end{theorem}

The three-dimensional version of Theorem~\ref{T:GLWe}, obtained by replacing
the truncated square diagram $\E{D}$ by a truncated \emph{cube} diagram,
does \emph{not} hold,
see J.~T\r{u}ma and F.~Wehrung \cite{TuWe1}. On the other hand, the
one-dimensional version of Theorem~\ref{T:GLWe} holds, see
Theorem~2 in G.~Gr\"atzer, H.~Lakser, and E.T.~Schmidt \cite{GLS1}, or
Theorem~4 in G.~Gr\"atzer, H.~Lakser, and E.T.~Schmidt \cite{GLS2}.

As a consequence of Theorem~2, we mention the result that \emph{every
distributive \jz-semilattice of cardinality at most $\aleph_1$ is isomorphic
to the semilattice of compact congruences of a relatively
complemented lattice with zero}, see \cite{GrLW}. Hence, lifting results of
finite character make it possible to prove representation results of infinite
character. The proofs of Theorems
\ref{T:GrScFin} and \ref{T:GLWe} do not extend to infinite $S$---in fact, we
do have a counterexample for the analogue of Theorem~\ref{T:GLWe} for
countable~$S$.

In this paper, we prove positive lifting results for \emph{infinite} $S$,
similar to Theorems \ref{T:GrScFin} and \ref{T:GLWe}. The only additional
assumption is that $S$ is a \emph{lattice}, just as in \cite{Schm81}.

Our first, most general theorem is the following.

\begin{all}{Theorem A}
Let $D$ be a distributive lattice with zero,
let $P$ be a partial lattice, let $\gf\colon\Conc P\to D$ be a
\jzh. If $\gf$ is `balanced', then it extends
to a \jzh\ $\gy\colon\Conc L\to D$, for a certain lattice $L$
generated, as a lattice, by $P$.
\end{all}

We refer to Section~\ref{S:ThmA} for a precise statement of Theorem~A. At
this point, we observe two facts:
\begin{itemize}
\item[---] There are, scattered in the literature, quite a number of
nonequivalent definitions of a partial lattice. For example, our definition
(see Definition~\ref{D:PartLatt}) is
tailored to provide, for a partial lattice $P$, an \emph{embedding} from
$\Con P$ into $\Con\FL(P)$, where $\FL(P)$ denotes the free lattice on~$P$.
It is \emph{not} equivalent to the definition presented in \cite{GLT2}.

\item[---] The condition that $\gf$ be `balanced' (see
Definition~\ref{D:SFSampled}) is quite complicated, which explains to a
large extent the size of this paper.
\end{itemize}

Intuitively, the condition that $\gf$ be balanced means that the computation
of finitely generated \emph{ideals} and \emph{filters}, as well as finite
intersections and joins of these, in every quotient of
$P$ by a prime ideal $G$ of~$D$, can be captured by \emph{finite}
amounts of information, and this uniformly on $G$. This condition is so
difficult to formulate that it appears at first sight as quite unpractical.

However, it is satisfied in two important cases, namely: either $P$ is a
\emph{lattice} (and then the statement of Theorem~A trivializes, as it
should), or $P$ is \emph{finite} with nonempty domains for the meet and the
join, see Proposition~\ref{P:LatFinDS}.
Although this observation is quite easy, the next one is far less trivial. It
shows that a large amount of \emph{amalgams} of balanced partial lattices and
homomorphisms are balanced, see Proposition~\ref{P:SquareAmalg}.

\begin{all}{Theorem B}
Let $D$ be a distributive lattice with zero.
Let $K$ be a finite lattice, let $P$ and $Q$ be partial lattices each of
them is either a finite partial lattice or a lattice, let $f\colon K\to P$ and
$g\colon K\to Q$ be homomorphisms of partial lattices, let
$\gm\colon\Conc P\to D$ and
$\gn\colon\Conc Q\to D$ such that $\gm\circ\Conc f=\gn\circ\Conc g$. Then
there exist a lattice~$L$, homomorphisms of partial lattices
$\oll{f}\colon P\to L$ and $\ol{g}\colon Q\to L$, and 
a \jzh\ $\gf\colon\Conc L\to D$ such that $\oll{f}\circ f=\ol{g}\circ g$,
$\gm=\gf\circ\Conc\oll{f}$, and $\gn=\gf\circ\Conc\ol{g}$. Furthermore, the
construction can be done in such a way that the following additional
properties hold:

\begin{enumerate}
\item $L$ is generated, as a lattice, by~$\oll{f}[P]\uu\ol{g}[Q]$.

\item The map $\gf$ isolates $0$.
\end{enumerate}
\end{all}

(We say that a map $\gf$ \emph{isolates $0$}, if $\gf(\gq)=0$
if{f} $\gq=0$, for all $\gq$ in the domain of~$\gf$.)

Unlike what happens with Theorem~A, stating Theorem~B does not require any
complicated machinery---it is an immediately usable tool.

We can now state our one-dimensional lifting result:

\begin{all}{Theorem C}
Let $K$ be a lattice, let $D$ be a distributive lattice with zero, and let
$\gf\colon\Conc K\to D$ be a \jzh. There are a relatively complemented
lattice $L$ of cardinality $|K|+|D|+\aleph_0$, a lattice homomorphism
$f\colon K\to L$, and an isomorphism $\ga\colon\Conc L\to D$ such that the
following assertions hold:
\begin{enumerate}
\item $\gf=\ga\circ\Conc f$.

\item The range of~$f$ is coinitial (resp., cofinal) in $L$.

\item If the range of~$\gf$ is cofinal in $D$, then the range of~$f$ is
internal in $L$.
\end{enumerate}
\end{all}

We observe that for a distributive \emph{semi}lattice $D$ with zero,
Theorem~C characterizes $D$ being a lattice, see \cite{TuWe2}.

Here, we say that a subset $X$ of a lattice $L$ is \emph{coinitial}
(\emph{cofinal}, \emph{internal}, resp.) if the upper subset (lower subset,
convex subset, resp.) generated by $X$ equals $L$.

The information that $L$ be relatively complemented in the statement of
Theorem~C reflects only part of the truth. It turns out that $L$ satisfies
certain strong closure conditions---we say that $\vv<L,\ga>$ is
\emph{internally saturated}, see Definition~\ref{D:EqnSat}. This
statement implies the following properties of~$L$, see
Proposition~\ref{P:SummSat} for details:

\begin{enumerate}
\item $L$ is relatively complemented.

\item $L$ has \emph{definable principal congruences}. More precisely, there
exists a positive existential formula $\gF(x,y,u,v)$ of the language of
lattice theory such that for every internally saturated
$\vv<L,\ga>$ and all $a$, $b$, $c$, $d\in L$,
 \[
 \gQ_L(a,b)\ci\gQ_L(c,d)\q\text{if{f}}\q L\text{ satisfies }\gF(a,b,c,d).
 \]
A similar result is easily seen to hold for statements of the form
$\gQ_L(a,b)\ci\JJ_{i<n}\gQ_L(c_i,d_i)$.
\end{enumerate}

Then Theorems B and C together imply easily the following two-dimensional
lifting result, that widely extends the main result of G. Gr\"atzer, H.
Lakser, and F. Wehrung \cite{GrLW}:

\begin{all}{Theorem D}
Let $K$, $P$, $Q$, $f$, $g$, $\gm$, $\gn$ satisfy the assumptions of
Theorem~\textup{B}. Then there are a relatively complemented lattice~$L$
of cardinality $|P|+|Q|+|D|+\aleph_0$, homomorphisms of partial lattices
$\oll{f}\colon P\to L$ and $\ol{g}\colon Q\to L$, and 
an \emph{isomorphism} $\gf\colon\Conc L\to D$ such that $\oll{f}\circ
f=\ol{g}\circ g$,
$\gm=\gf\circ\Conc f$, and $\gn=\gf\circ\Conc g$. Furthermore, the
construction can be done in such a way that the following additional
properties hold:

\begin{enumerate}
\item The subset $\oll{f}[P]\uu\ol{g}[Q]$ generates $L$ as an ideal (resp.,
filter).

\item If the subsemilattice of $D$ generated by
$\gm[\Conc P]\cup\gn[\Conc Q]$ is cofinal in~$D$, then
$\oll{f}[P]\uu\ol{g}[Q]$ generates $L$ as a convex sublattice.
\end{enumerate}
\end{all}

Furthermore, once Theorem~C is proved, easy corollaries follow. For example,

\begin{all}{Corollary~\ref{C:EmbRelCpl}}
Every lattice $K$ such that $\Conc K$ is a lattice has an
internal, congruence-preserving embedding into a relatively complemented
lattice.
\end{all}

\begin{all}{Corollary~\ref{C:omegaDistr}}
Every \jz-semilattice that is a direct limit of a countable
sequence of distributive lattices with zero is isomorphic to the semilattice
of compact congruences of a relatively complemented lattice with zero
\end{all}

\subsection*{Methods}
Our methods of proof, especially for Theorems A and B, are radically
different from the usual `finite' methods, for example, those used in the
proofs of Theorems~\ref{T:GrScFin} and \ref{T:GLWe}.
In some sense, we take the most naive possible approach of the problem. We
are given a partial lattice $P$, a distributive lattice $D$ with zero, a
homomorphism $\gf\colon\Conc P\to D$, and we wish to ``extend $P$ to a
relatively complemented lattice $L$, and make $\gf$ an isomorphism'', as in
the statement of Theorem~A. So we ``add new joins and meets'' in order to
make $P$ a \emph{total} lattice (we use Theorem~A), we ``add relative
complements'' in order to make $P$ relatively complemented (see
Lemma~\ref{L:RelCpl}), we ``force projectivity of intervals'' in order to
make $\gf$ an embedding (see Lemmas \ref{L:Persp}--\ref{L:gfEmb}), and we
``add new intervals'' in order to make $\gf$ surjective (see
Lemma~\ref{L:gfVemb}). Of course, the main problem is then to confine the
range of~$\gf$ within~$D$.

In this sense, this approach is
related to G.~Gr\"atzer and E.T.~Schmidt's \cite{GrSc} proof of the
representation problem of congruence lattices of algebras, see also P.
Pudl\'ak \cite{Pudl76} and Section~2.3 in E.T. Schmidt \cite{Schm82}: given an
algebraic (not necessarily distributive) lattice $A$, a partial algebra $U$
is constructed such that $\Con U\iso A$, then $U$ is extended to a total
algebra with the same congruence lattice.

However, there is an important difference between this approach and ours,
namely: in Gr\"atzer and Schmidt's proof, infinitely many new operations need
to be incorporated to the signature of the algebra. This restriction is
absolutely unavoidable, as proves W.~Lampe's result
(a stronger version was proved independently by R.~Freese and W.~Taylor) that
certain algebraic lattices require many operations to be represented, see
R.~Freese, W.~Lampe, and W.~Taylor \cite{FLT}, or Section~2.4 in E.T.~Schmidt
\cite{Schm82}. In the present paper, we are restricted to the language of
lattice theory $\vv<\jj,\mm>$. This may partly explain our restriction to
algebraic lattices which are ideal lattices of distributive \emph{lattices}.
That the latter restriction is necessary is established in the forthcoming
paper J.~T\r{u}ma and F.~Wehrung
\cite{TuWe2}.\smallskip

To get around this difficulty, we borrow the notations and methods of the
theory of forcing and Boolean-valued models. Although it has been recognized
that the latter are, in universal algebra, a more convenient framework than
the usual sheaf representation results, see, for example, Chapter~IV
in S.~Burris and H.P.~Sankappanavar~\cite{BuSa}, one can probably not say
that they are, at the present time, tools of common use in lattice theory.
For this reason, our presentation will assume no familiarity with
Boolean-valued models. We refer the reader, for example, to T.~Jech
\cite{Jech89} for a presentation of this topic.

The basic idea of the present paper is, actually, quite simple. For
a partial lattice $P$, we consider the standard construction of the
free lattice $\FL(P)$ on~$P$. More specifically, $\FL(P)$ is constructed as
the set of \emph{words} on $P$, using the binary operations $\jj$ and $\mm$.
The \emph{ordering} on $\FL(P)$ is defined inductively, see
Definition~\ref{D:peq}. This can be done by assigning to every statement
of the form $\dx\leq\dy$, where $\dx$ and $\dy$ are words on~$P$, a
`truth value' $\bv{\dx\leq\dy}$, equal either to $0$ (false) or $1$ (true).
So $\bv{\dx\leq\dy}$ equals $1$ if $\dx\leq\dy$, $0$ otherwise. So, for
example, rule (ii) of Definition~\ref{D:peq} may be stated as
 \begin{equation}\label{Eq:TruthVal01}
 \bv{\dx_0\jj\dx_1\leq\dy_0\mm\dy_1}=\MM_{i,j<2}\bv{\dx_i\leq\dy_j}.
 \end{equation}
If the truth values of statements are no longer confined to
$\set{0,1}$ but, rather, to elements of a given distributive lattice $D$
(which has to be thought as the \emph{dual} lattice of the lattice $D$ of
the statements of Theorems A and B), \eqref{Eq:TruthVal01} becomes part
of the inductive definition of a map that with
every pair $\vv<\dx,\dy>$ of words on~$P$ associates the `truth value'
$\bv{\dx\leq\dy}\in D$, that we shall still call ``Boolean value'' (after
all, $D$ embeds into a Boolean algebra).

In this way, it seems at first sight a trivial task to extend
Definition~\ref{D:peq} to a $D$-valued context. However, the major obstacle
remains of the computation of~$\bv{\dx\leq\dy}$ at the
\emph{bottom} level, that is, for $\dx$ and $\dy$ finite meets or joins of
elements of~$P$. This situation is not unlike what happens in set theory,
where the main problem in defining Boolean values in the Scott-Solovay
Boolean universe is to define them on the \emph{atomic} formulas, see
\cite{Jech78}. In fact, the method used in \cite{Wehr93} reflects more
closely what is done in the present paper, namely, the domain of the Boolean
value function is extended from a set of `urelements' to the universe of set
theory that they generate.

The condition that $\vv<P,\gf>$ be balanced is designed to ensure that these
Boolean values belong to $D$, while they would typically, in the general
case, belong to the completion of the universal Boolean algebra of~$D$.

The reader may feel at this point a slight uneasiness, because the
distributive lattice $D$ in which the Boolean values live is related to
the \emph{dual} of~$\Conc P$, rather than to $\Conc P$ itself. It seems,
indeed, pointless to dualize $D$, prove a large amount of results on the
dual, and then dualize again to recover $D$. Why bother doing this? The
alternative would be to stick with the original $D$, and so, to interpret
the Boolean values by $\bv{\dx\leq\dx}=0$ (instead of~$1$, `true'), and
$\bv{\dx\leq\dz}\leq\bv{\dx\leq\dy}\jj\bv{\dy\leq\dz}$ (instead of the
dual, see Definition~\ref{D:DvalPoset}(ii)). Furthermore, one would have
to interpret the propositional connective `and' by the join $\jj$, and
`or' by the meet $\mm$, and so on. This is definitely unattractive to
the reader familiar with Boolean models. Of course, the last decisive
argument for one way or the other is merely related to a matter of taste.
\medskip

We now give a short summary of the paper, part by part.

\textbf{Part~\ref{Pt:PartLatt}} introduces partial lattices and their
congruences, and also the free lattice on a partial lattice. A noticeable
difference between our definition of a congruence and the usual definition
of a congruence is that our congruences are \emph{not symmetric} in
general. The reason for this is very simple, namely, if $f\colon P\to Q$
is a homomorphism of partial lattices, then its \emph{kernel}, instead of
being defined as usual as the set of all pairs $\vv<x,y>$ such that
$f(x)=f(y)$, is defined here as the set of all $\vv<x,y>$ such that
$f(x)\leq f(y)$. Of course, for (total) lattices, the two resulting
definitions of a congruence are essentially equivalent---in particular, they
give isomorphic congruence lattices.

In Section~\ref{S:AmalgPart}, we interpret the classical operation of
`pasting' two partial lattices above a lattice as a \emph{pushout} in the
category of partial lattices and homomorphisms of partial lattices.
Although the description of the pushout, Proposition~\ref{P:PushStV}, is
fairly straightforward, it paves the way for its $D$-valued analogue,
Proposition~\ref{P:DPushStV}.

\textbf{Part~\ref{Pt:DVal}} begins with the simple definitions,
in Section~\ref{S:DvalPL}, of a $D$-valued poset or of a $D$-valued partial
lattice. The purpose of Section~\ref{S:JMSample} is to introduce the
important definition of a \emph{sample}, that makes it possible,
\emph{via} additional assumptions, to extend to the $D$-valued world the
classical notions of \emph{ideal} and \emph{filter} of a partial lattice
$P$, see Section~\ref{S:IdFil}. The corresponding $D$-valued notions,
instead of corresponding to \emph{subsets} of~$P$, correspond to
\emph{functions} from $P$ to $D$.

However, the objects we wish to solve problems about are not $D$-valued
partial lattices, but plain partial lattices. Thus we present in
\textbf{Part~\ref{Pt:DComeas}} a class of structures that live
simultaneously in both worlds, the \emph{\CMPL{D}s}, to which we
translate the results of Part~\ref{Pt:DVal}. In order to extend to a
\emph{lattice} the Boolean values defined on the original partial lattice,
we introduce the definition of a \emph{balanced} \CMPL{D}, see
Definition~\ref{D:SFSampled}. Then we prove, in Sections \ref{S:IdFiM}
and \ref{S:IdFiJ}, that all our finiteness conditions (they add up to the
condition of being balanced) are preserved under amalgamation above a finite
lattice.

Now that all this hard technical work is completed, we start applying it in
\textbf{Part~\ref{Pt:CongDistr}}. Most arguments used in this part are
based on simple amalgamation constructions of partial lattices above
finite lattices, that all yield, by our previous work, balanced \CMPL{D}s.

\section*{Notation and terminology}

For a set $X$, we denote by $\fin X$ (resp., $\fine X$) the set of all
finite (resp., nonempty finite) subsets of~$X$.

We put $\two=\set{0,1}$, endowed with its canonical structure of lattice.
For a nonnegative integer $n$, we identify $n$ with $\set{0,1,\ldots,n-1}$.

Let $P$ be a preordered set. For subsets $X$, $Y$ of~$P$, let $X\leq Y$ be
the statement $\forall x\in X$, $\forall y\in Y$, $x\leq y$. We shall write
$a\leq X$ (resp., $X\leq a$) instead of~$\set{a}\leq X$ (resp.,
$X\leq\set{a}$).
A subset $X$ of~$P$ is a \emph{lower subset} (resp., \emph{upper subset})
of~$P$ if for all $x\leq y$ in $P$, $y\in X$ (resp. $x\in X$) implies
that $x\in X$ (resp., $y\in X$). We say that $X$ is a \emph{convex
subset} of~$P$, if $a\leq x\leq b$ and $\set{a,b}\ci X$ implies that
$x\in X$, for all $a$, $b$, $x\in P$.

If $X\ci P$, we denote by $\dnw X$ (resp., $\upw X$) the
lower subset (resp., upper subset) of~$P$ generated by $X$. For $a\in P$,
we put $\dnw a=\dnw\set{a}$ and $\upw a=\upw\set{a}$.

For $a\in P$ and $X\ci P$, let $a=\sup X$ be the statement
 \[
 X\leq a\quad\text{and}\quad\forall x,\ X\leq x\Rightarrow a\leq x.
 \]
The statement $a=\inf X$ is defined dually. Note that if $a=\sup X$, then
$a'=\sup X$ for all $a'$ equivalent to $a$ with respect to the preordering
$\leq$ (that is, $a\leq a'\leq a$).

For a preordering $\ga$ of a set $P$ and for $x$, $y\in P$, the statement
$\vv<x,y>\in\ga$ will often be abbreviated $x\leq_\ga y$.

For a lattice $L$, $L^\rd$ denotes the \emph{dual lattice} of~$L$.

\part{Partial lattices}\label{Pt:PartLatt}

\section{Partial prelattices and partial lattices}\label{S:PpPl}

\begin{definition}\label{D:PartLatt}
A \emph{partial prelattice} is a structure $\vv<P,\leq,\JJ,\MM>$, where
$P$ is a nonempty set,
$\leq$ is a preordering on~$P$, and $\JJ$, $\MM$ are partial functions from
$\fine P$ to $P$ satisfying the following properties:
\begin{enumerate}
\item $a=\JJ X$ implies that $a=\sup X$,
for all $a\in P$ and all $X\in\fine P$.

\item $a=\MM X$ implies that $a=\inf X$,
for all $a\in P$ and all $X\in\fine P$.
\end{enumerate}

We say that $P$ is a \emph{partial lattice}, if $\leq$ is antisymmetric.

A \emph{congruence} of~$P$ is a preordering $\preceq$ of~$P$ containing
$\leq$ such that $\vv<P,\preceq,\JJ,\MM>$ is a partial prelattice.

If $P$ and $Q$ are partial prelattices, a \emph{homomorphism of partial
prelattices} from $P$ to $Q$ is an order-preserving map
$f\colon P\to Q$ such that $a=\JJ X$
(resp., $a=\MM X$) implies that $f(a)=\JJ f[X]$ (resp., $f(a)=\MM f[X]$),
for all $a\in P$ and all $X\in\fine P$. We say that a homomorphism $f$ is an
\emph{embedding}, if $f(a)\leq f(b)$ implies that $a\leq b$, for all
$a$, $b\in P$.
\end{definition}

We shall naturally identify \emph{lattices} with partial lattices $P$ such
that $\JJ$ and $\MM$ are defined everywhere on $\fine P$.

\begin{remark}
For an embedding $f\colon P\to Q$ of partial lattices, we do \emph{not}
require that $\JJ f[X]$ be defined implies that $\JJ X$ is defined (and
dually), for $X\in\fine P$.
\end{remark}

\begin{proposition}\label{P:ConP}
Let $P$ be a partial prelattice. Then the set $\Con P$ of all congruences of
$P$ is a closure system in the powerset lattice of~$P\times P$, closed under
directed unions. In particular, it is an algebraic lattice.
\end{proposition}

We denote by $\Conc P$ the \jz-semilattice of all \emph{compact}
congruences of~$P$, by $\zero_P$ the least congruence of~$P$
(that is, $\zero_P$ is the preordering of~$P$), and by $\one_P$ the largest
(coarse) congruence of~$P$. The map $P\mapsto\Conc P$ can be extended in a
natural way in a \emph{functor}, as follows. For a homomorphism
$f\colon P\to Q$ of partial lattices, we define a \jzh\
$\Conc f\colon\Conc P\to\Conc Q$ as the map that with every congruence
$\alpha$ of $P$ associates the congruence of $Q$ generated by all pairs
$\vv<f(x),f(y)>$, for $\vv<x,y>\in\alpha$.

If $P$ is a lattice, then $\Con P$ is distributive, but
this may not hold for a general partial lattice $P$.

For $a$, $b\in P$, we denote by $\gQ_P^+(a,b)$ the least congruence
$\gq$ of~$P$ such that $a\leq_\gq b$, and we put
$\gQ_P(a,b)=\gQ_P^+(a,b)\jj\gQ_P^+(b,a)$, the least congruence $\gq$ of~$P$
such that $a\equiv_\gq b$. Of course, the congruences of the form
$\gQ_P^+(a,b)$ are generators of the join-semilattice $\Conc P$.

\section{The free lattice on a partial lattice}

We present in this section an explicit construction, due to R.A. Dean
\cite{Dean64}, of the free lattice on a partial lattice, see also
\cite[Page~249]{FJN}. For the needs of this paper, the definitions are
slightly modified (in particular, the relation $\peq$ defined below, see
Definition~\ref{D:peq}), but it is easy to verify that they are, in fact,
equivalent to the original ones.

Throughout this section, we shall fix a partial lattice $P$.

\subsection{Ideals, filters}\label{S:IdFil}

\begin{definition}\label{D:IdFil}
An \emph{ideal} of~$P$ is a lower subset $I$ of~$P$ such that $X\ci I$ and
$a=\JJ X$ imply that $a\in I$, for all $X\in\fine P$ and all $a\in P$.
Dually, a \emph{filter} of~$P$ is an upper subset $F$ of~$P$ such that $X\ci
F$ and $a=\MM X$ imply that $a\in I$, for all $X\in\fine P$ and all $a\in P$.
\end{definition}

We observe that both $\es$ and $P$ are simultaneously an ideal and a filter
of~$P$. For $a\in P$, $\dnw a$ is an ideal of~$P$ (\emph{principal ideal}),
while $\upw a$ is a filter of~$P$ (\emph{principal filter}). In case $P$ is
a \emph{lattice} (that is, $\JJ$ and $\MM$ are everywhere defined), the
ideals of the form $\dnw a$ are the only nonempty finitely generated ideals
of~$P$.

\begin{lemma}\label{L:IdFilAL}
The set $\ID(P)$ (resp., $\FIL(P)$) of all ideals (resp., filters) of~$P$ is
a closure system in the powerset lattice $\mathcal{P}(P)$ of~$P$, closed
under directed unions. Hence, both $\ID(P)$ and $\FIL(P)$ are
\emph{algebraic lattices}.
\end{lemma}

\subsection{Description of the free lattice on~$P$}\label{S:DescFL(P)}

\begin{notation}
For a set $\gO$, let $\W(\gO)$ denote the set of \emph{terms} on
$\gO$ and the two binary operations $\jj$ and $\mm$.
\end{notation}

So, the elements of~$\W(\gO)$ are formal ``polynomials'' on the elements of
$\gO$, such as $((a\jj b)\jj c)\mm(d\jj e)$,
where $a$, $b$, $c$, $d$, $e\in\gO$, \emph{etc.}.
The \emph{height} of an element $\dx$ of~$\W(\gO)$ is defined inductively by
$\hgt(a)=0$ for $a\in\gO$, and
$\hgt(\dx\mm\dy)=\hgt(\dx\jj\dy)=\hgt(\dx)+\hgt(\dy)+1$.

We shall now specialize to the case where $\gO$ is the underlying set of
the partial lattice $P$. (In Section~\ref{S:BVtoW(P)}, the notation $\W(P)$
will be used for structures $P$ that are not necessarily partial lattices.)

\begin{definition}\label{D:x-x+class}
For $\dx\in\W(P)$, we define, by induction on the
height $\hgt(\dx)$ of~$\dx$, an ideal $\dx^-$ of~$P$ and a filter $\dx^+$ of
$P$ as follows:

\begin{enumerate}
\item $\dx^-=\dnw a$ and $\dx^+=\upw a$, if $\dx=a\in P$.

\item If $\dx=\dx_0\jj\dx_1$, we put $\dx^-=\dx_0^-\jj\dx_1^-$ (the join
being computed in $\ID(P)$), and $\dx^+=\dx_0^+\ii\dx_1^+$.

\item If $\dx=\dx_0\mm\dx_1$, we put $\dx^-=\dx_0^-\ii\dx_1^-$, and
$\dx^+=\dx_0^+\jj\dx_1^+$  (the join being computed in $\FIL(P)$).
\end{enumerate}
\end{definition}

\begin{definition}\label{D:xlly}
For $\dx$, $\dy\in\W(P)$, we define $\dx\ll\dy$ to hold, if
$\dx^+\ii\dy^-\ne\es$.
\end{definition}

\begin{definition}\label{D:peq}
We define inductively a binary relation $\peq$ on $\W(P)$, as follows:
\begin{enumerate}
\item $\dx\peq\dy$ iff $\dx\ll\dy$, for all $\dx$, $\dy\in\W(P)$
such that $\dx\in P$ or $\dy\in P$.

\item $\dx_0\jj\dx_1\peq\dy_0\mm\dy_1$ if{f} $\dx_i\peq\dy_j$,
for all $i$, $j<2$.

\item $\dx_0\jj\dx_1\peq\dy_0\jj\dy_1$ if{f} $\dx_i\peq\dy_0\jj\dy_1$, for
all $i<2$.

\item $\dx_0\mm\dx_1\peq\dy_0\mm\dy_1$ if{f} $\dx_0\mm\dx_1\peq\dy_j$, for
all $j<2$.

\item $\dx_0\mm\dx_1\peq\dy_0\jj\dy_1$ if{f} either
$\dx_0\mm\dx_1\ll\dy_0\jj\dy_1$ or $\dx_i\peq\dy_j$, for some
$i$, $j<2$.
\end{enumerate}
\end{definition}

The relevant observations can be summarized in the following form:

\begin{lemma}\label{L:x-x+llpeq}
Let $a$, $b\in P$ and let $\dx$, $\dy$, $\dz\in\W(P)$.
Then the following assertions hold:
\begin{enumerate}
\item $a\in\dx^-$ and $b\in\dx^+$ imply that $a\leq b$;

\item $a\in\dx^-$ and $\dx\peq\dy$ imply that $a\in\dy^-$;

\item $a\in\dy^+$ and $\dx\peq\dy$ imply that $a\in\dx^+$;

\item $\dx\ll\dy$ implies that $\dx\peq\dy$;

\item $\dx\peq\dx$;

\item $\dx\peq\dy$ and $\dy\peq\dz$ imply that $\dx\peq\dz$.

\end{enumerate}
\end{lemma}

Let $\equiv$ denote the equivalence relation associated with the preordering
$\peq$. We define $\FL(P)=\vv<\W(P),\peq>/{\equiv}$. Let
$j_P\colon P\to\FL(P)$, the \emph{natural map}, be defined by
$j_P(a)=\cls(a,\equiv)$, for all $a\in P$.

\begin{proposition}\label{P:DescFL}
The poset $\FL(P)$ is a lattice and $j_P$ is an embedding of partial
lattices. Furthermore, $j_P$ is universal among all the homomorphisms of
partial lattices from $P$ to a lattice.
\end{proposition}

So we identify $\FL(P)$ (together with the natural map $j_P$) with the
\emph{free lattice on the partial lattice $P$}, that is, the lattice defined
by generators $\check{a}$ ($a\in P$) and relations
$\check{a}=\check{x}_0\jj\cdots\jj\check{x}_{n-1}$ (resp.,
$\check{a}=\check{x}_0\mm\cdots\mm\check{x}_{n-1}$) if
$a=\JJ\set{x_0,\ldots,x_{n-1}}$ (resp., $a=\MM\set{x_0,\ldots,x_{n-1}}$) in
$P$.

\subsection{Generation of ideals and filters}\label{S:GenIdFil}

In any lattice, the finitely generated ideals are exactly the principal
ideals, and similarly for filters. In general partial lattices, the situation
is much more complicated. The somewhat more precise description of ideals
and filters that we shall give in this section will be used later
in Section~\ref{S:IdFilSamples}.

\begin{definition}\label{D:IdFiln01}
Let $X$ and $U$ be subsets of~$P$. For $n<\go$, we define, by induction on
$n$, a subset $\Id_n(X,U)$ of~$P$, as follows:
\begin{enumerate}
\item $\Id_0(X,U)=\dnw X$.

\item $\Id_{n+1}(X,U)$ is the union of~$\Id_n(X,U)$ and the
lower subset of~$P$ generated by all elements
of the form $\JJ Z$, where $\es\sci Z\ci U\ii\Id_n(X,U)$ and $\JJ Z$ is
defined ($\sci$ denotes \emph{proper} inclusion).
\end{enumerate}
Dually, we define, by induction on $n$, a subset $\Fil_n(X,U)$ of~$P$, as
follows:
\begin{enumerate}
\item $\Fil_0(X,U)=\upw X$.

\item $\Fil_{n+1}(X,U)$ is the union of~$\Fil_n(X,U)$ and the
upper subset of~$P$ generated by all elements
of the form $\MM Z$, where $\es\sci Z\ci U\ii\Fil_n(X,U)$ and $\MM Z$ is
defined.
\end{enumerate}
\end{definition}

We observe, in particular, that $\UU_{n<\go}\Id_n(X,P)$
is the ideal $\Id(X)$ of~$P$ generated by $X$. The subsets $\Id_n(X,U)$, for
finite $U$, can be viewed as ``finitely generated approximations'' of
$\Id(X)$. Similar considerations hold for $\Fil_n(X,U)$ and
$\Fil(X)=\UU_{n<\go}\Fil_n(X,P)$.

\section{Amalgamation of partial lattices above a lattice}
\label{S:AmalgPart}

Most of the results of this section are folklore, we recall them here for
convenience.

\begin{definition}\label{D:V-Form}
A \emph{V-formation} of partial lattices is a structure
$\vv<K,P,Q,f,g>$ subject to the following conditions:
\begin{itemize}
\item[(V1)] $K$, $P$, $Q$ are partial lattices.

\item[(V2)] $f\colon K\hookrightarrow P$ and $g\colon K\hookrightarrow Q$ are
embeddings of partial lattices.
\end{itemize}

A V-formation $\vv<K,P,Q,f,g>$ is \emph{standard}, if the following
conditions hold:
\begin{itemize}
\item[(SV1)] $K$ is a lattice.

\item[(SV2)] $K=P\ii Q$ (set-theoretically), and $f$ and $g$ are,
respectively, the inclusion map from $K$ into $P$ and the inclusion map from
$K$ into $Q$.
\end{itemize}

\end{definition}

Of course, we identify a V-formation $\vv<K,P,Q,f,g>$ with the
\emph{diagram} of partial lattices that consists of two arrows from $K$, one
of them $f\colon K\to P$, the other $g\colon K\to Q$.

Furthermore, the homomorphisms in standard V-formations are understood
(they are the inclusion maps), so, in that case, we shall write $\vv<K,P,Q>$
instead of~$\vv<K,P,Q,f,g>$.

The following lemma is a set-theoretical triviality:

\begin{lemma}\label{L:V-FormSt}
Every V-formation $\vv<K,P,Q,f,g>$ of partial lattices, with $K$ a lattice,
is isomorphic to a standard V-formation.
\end{lemma}

\begin{definition}\label{D:AmalgVform}
Let $\E{D}=\vv<K,P,Q,f,g>$ be a V-formation of partial lattices. An
\emph{amalgam} of~$\E{D}$ is a triple $\vv<R,f',g'>$, where $R$ is a partial
lattice and $f'\colon P\hookrightarrow R$, $g'\colon Q\hookrightarrow R$ are
embeddings of partial lattices such that $f'\circ f=g'\circ g$.

As usual in category-theoretical terminology, we say that a \emph{pushout}
of~$\E{D}$ is any initial object in the category of amalgams of~$\E{D}$ with
their \emph{homomorphisms} (not only embeddings). Of course, if the pushout
of~$\E{D}$ exists, then it is unique up to isomorphism.
\end{definition}

We shall be concerned about not only the existence but also the
\emph{description} of pushouts in a very precise context:

\begin{proposition}\label{P:PushStV}
Let $\E{D}=\vv<K,P,Q,f,g>$ be a V-formation of partial lattices, with $K$ a
lattice. Then $\E{D}$ has a pushout. Furthermore, assume that $\E{D}$ is a
standard V-formation. Then the pushout $\vv<R,f',g'>$ of~$\E{D}$ can be
described by the following data:
\begin{itemize}
\item[(a)] $R=P\uu Q$, endowed with the partial ordering $\leq$ consisting
of all pairs $\vv<x,y>$ of elements of~$R$ satisfying the following
conditions:

\begin{itemize}
\item[(a1)] $x$, $y\in P$ and $x\leq_Py$.

\item[(a2)] $x$, $y\in Q$ and $x\leq_Qy$.

\item[(a3)] $x\in P$, $y\in Q$, and there exists $z\in K$ such that
$x\leq_Pz$ and $z\leq_Qy$.

\item[(a4)] $x\in Q$, $y\in P$, and there exists $z\in K$ such that
$x\leq_Qz$ and $z\leq_Py$.
\end{itemize}

\item[(b)] For $a\in R$ and $X\in\fine R$, $a=\JJ X$ holds in $R$ if{f}
either $X\uu\set{a}\ci P$ and $a=\JJ X$ in $P$ or
$X\uu\set{a}\ci Q$ and $a=\JJ X$ in $Q$.

\item[(b*)] For $a\in R$ and $X\in\fine R$, $a=\MM X$ holds in $R$ if{f}
either $X\uu\set{a}\ci P$ and $a=\MM X$ in $P$ or
$X\uu\set{a}\ci Q$ and $a=\MM X$ in $Q$.

\item[(c)] $f'$ (resp., $g'$) is the inclusion map from $P$ into $R$ (resp.,
from $Q$ into $R$).
\end{itemize}

\end{proposition}

\begin{note}
It is easy to prove that \emph{any} diagram of partial lattices admits a
colimit. In particular, pushouts always exist. However, we are, in
Proposition~\ref{P:PushStV}, more interested in the \emph{description} of
the pushout.

\end{note}

\begin{proof}
The fact that the binary relation $\leq$ defined above on $R$ is a partial
ordering is folklore (and easy to verify).

Now we prove that $R$ is a partial lattice. We first observe that since $K$
is a partial sublattice of both $P$ and $Q$, the partial operations $\JJ$
and $\MM$ on $R$ described in (b) and (b*) above are, indeed, partial
functions.

Let $\vv<a,X>\in R\times\fine R$ such that $a=\JJ X$ in $R$, we prove that
$a=\sup X$ in $R$. By the definition of~$\JJ$ in $R$, $a=\JJ X$ holds either
in $P$ or in $Q$, so, without loss of generality, $X\cup\set{a}\ci P$ and
$a=\JJ X$ in $P$. Since $P$ is a partial lattice, it follows that
 \begin{equation}\label{Eq:a=supXinP}
 a=\sup X\q\text{in }P.
 \end{equation}
{F}rom $X\leq_Pa$ follows that $X\leq a$. Now let $b\in R$
such that $X\leq b$, we prove that $a\leq b$. If $b\in P$, then $X\leq_Pb$,
thus, by \eqref{Eq:a=supXinP}, $a\leq_Pb$, so $a\leq b$.

Now suppose that $b\in Q$. For all $x\in X$, $x\leq b$ with $x\in P$ and
$b\in Q$, thus there exists $x^*\in K$ such that
 \begin{align}
 x&\leq_Px^*\label{Eq:xleqx*}\\
 x^*&\leq_Qb.\label{Eq:x*leqb}
 \end{align}
Since $K$ is a lattice, $c=\JJ_{x\in X}x^*$ is defined in $K$.
Since $K$ is a partial sublattice of~$Q$, the equality $c=\JJ_{x\in X}x^*$
also holds in $Q$. Thus, by \eqref{Eq:x*leqb}, we obtain the inequality
 \begin{equation}\label{Eq:cleqQb}
 c\leq_Qb.
 \end{equation}
Furthermore, for $x\in X$, $x^*\leq_Kc$, thus $x^*\leq_Pc$, hence, by
\eqref{Eq:xleqx*}, $x\leq_Pc$. This holds for all $x\in X$, thus, by
\eqref{Eq:a=supXinP}, we obtain the inequality
 \begin{equation}\label{Eq:aleqPc}
 a\leq_Pc.
 \end{equation}
Hence, by \eqref{Eq:aleqPc} and \eqref{Eq:cleqQb}, $a\leq b$. Therefore,
$a=\sup X$ in $R$.

The proof for $\MM$ and $\inf$ is similar.

Finally, the proof that $\vv<R,f',g'>$ is a pushout of~$\E{D}$ is
straightforward.
\end{proof}

\begin{notation}\label{Not:PamalgQ}
In the context of Proposition~\ref{P:PushStV}, in the case of a standard
V-formation $\vv<K,P,Q>$, we shall write $R=P\amalg_KQ$.
\end{notation}

\part{$D$-valued posets and partial lattices}\label{Pt:DVal}

\section{$D$-valued posets}\label{S:DvalPL}

We shall fix in this section a distributive lattice $D$ with unit (largest
element) $1$.

The following definition is similar to the classical definition of a
Boolean-valued model, see, for example, \cite{Jech89}.

\begin{definition}\label{D:DvalPoset}
A \emph{$D$-valued poset} is a nonempty set $P$, together with a map
$P\nobreak\times\nobreak P\nobreak\to\nobreak D$,
$\vv<a,b>\mapsto\bv{a\leq b}$, that satisfies the following properties:
\begin{enumerate}
\item $\bv{a\leq a}=1$, for all $a\in P$.

\item $\bv{a\leq b}\mm\bv{b\leq c}\leq\bv{a\leq c}$, for all $a$, $b$,
$c\in P$.
\end{enumerate}

\end{definition}

If $P$ is a $D$-valued poset, then we define
$\bv{a=b}=\bv{a\leq b}\mm\bv{b\leq a}$, for all $a$, $b\in P$.
Furthermore, for $a\in P$ and nonempty, finite subsets $X$ and $Y$ of~$P$, we
put
 \begin{align*}
 \bv{a\in Y}&=\JJ_{y\in Y}\bv{a=y},\\
 \bv{X\ci Y}&=\MM_{x\in X}\bv{x\in Y}.
 \end{align*}
and we put $\bv{X=Y}=\bv{X\ci Y}\mm\bv{Y\ci X}$.

We observe that a $D$-valued poset is \emph{not} given with a partial
ordering on~$P$---there is no such thing as ``the binary relation $\leq$ on
$P$''. Instead, $\bv{a\leq b}$ denotes an
\emph{element of~$D$}, as opposed to a \emph{statement}.

\begin{example}\label{Ex:01ValPos}
Let $\vv<P,\leq>$ be a poset. Then $P$ can be canonically endowed with a
structure of $\two$-valued poset, by putting $\bv{a\leq b}=1$ if $a\leq b$,
$0$ otherwise.
\end{example}

Hence, the rule $\bv{a\leq b}\mm\bv{b\leq c}\leq\bv{a\leq c}$ may be
interpreted as the $D$-valued version of the transitivity of the partial
ordering.

We record below some basic facts about Boolean values. For the remainder of
this section, we fix a $D$-valued poset $P$. In fact, many of the results
below hold for $D$-valued models of equality (a set $P$ with a map
$\vv<x,y>\mapsto\bv{x=y}$), with the same proofs.

\begin{lemma}\label{L:XciYciZ}
The following assertions hold:
\begin{enumerate}
\item $\bv{x=y}\mm\bv{y\in Z}\leq\bv{x\in Z}$, for all $x$, $y\in P$
and all $Z\in\fine P$.

\item $\bv{x\in Y}\mm\bv{Y\ci Z}\leq\bv{x\in Z}$, for all $x\in P$ and all
$Y$, $Z\in\fine P$.

\item $\bv{X\ci Y}\mm\bv{Y\ci Z}\leq\bv{X\ci Z}$, for all $X$, $Y$,
$Z\in\fine P$.

\item $\bv{X=Y}\mm\bv{Y=Z}\leq\bv{X=Z}$, for all $X$, $Y$,
$Z\in\fine P$.
\end{enumerate}
\end{lemma}

\begin{proof}
(i) We compute:
 \[
 \bv{x=y}\mm\bv{y\in Z}=\JJ_{z\in Z}\bv{x=y}\mm\bv{y=z}\leq
 \JJ_{z\in Z}\bv{x=z}=\bv{x\in Z}.
 \]

(ii) We compute, by using (i):
 \[
 \bv{x\in Y}\mm\bv{Y\ci Z}\leq\JJ_{y\in Y}\bv{x=y}\mm\bv{Y\ci Z}
 \leq\JJ_{y\in Y}\bv{x=y}\mm\bv{y\in Z}\leq
 \bv{x\in Z}.
 \]

(iii) We compute, by using (ii):
 \[
 \bv{X\ci Y}\mm\bv{Y\ci Z}=\MM_{x\in X}\bv{x\in Y}\mm\bv{Y\ci Z}
 \leq\MM_{x\in X}\bv{x\in Z}=\bv{X\ci Z}.
 \]

(iv) is an obvious consequence of (iii).
\end{proof}

\begin{lemma}\label{L:XciYya}
Let $a\in P$, let $X$, $Y\in\fine P$, let $\gf(z,a)$ be one of the formulas
$z\leq a$ or $a\leq z$.
Then the following inequalities hold:
\begin{enumerate}
\item
$\bv{X\ci Y}\mm\MM_{y\in Y}\bv{\gf(y,a)}\leq
\bv{X\ci Y}\mm\MM_{x\in X}\bv{\gf(x,a)}$.

\item
$\bv{X=Y}\mm\MM_{y\in Y}\bv{\gf(y,a)}=
\bv{X=Y}\mm\MM_{x\in X}\bv{\gf(x,a)}$.

\end{enumerate}
\end{lemma}

\begin{proof}
(i) Put $\gg=\bv{X\ci Y}\mm\MM_{y\in Y}\bv{\gf(y,a)}$. For $x\in X$,
 \[
 \gg\leq\bv{X\ci Y}\leq\bv{x\in Y}=\JJ_{y\in Y}\bv{x=y},
 \]
so, to prove (i), it is sufficient to prove that
$\gg\mm\bv{x=y}\leq\bv{\gf(x,a)}$, for all $y\in Y$. But this follows from
the fact that $\gg\mm\bv{x=y}\leq\bv{x=y}\mm\bv{\gf(y,a)}$ and the
definition of a $D$-valued poset.

(ii) follows immediately from (i).
\end{proof}

\begin{lemma}\label{L:XciYX=Z}
Let $X$, $Y\in\fine P$. Then the following equality holds:
 \[
 \bv{X\ci Y}=\JJ_{\es\sci Z\ci Y}\bv{X=Z}.
 \]
\end{lemma}

\begin{proof}
For $\es\sci Z\ci Y$, the inequality $\bv{X=Z}\leq\bv{X\ci Y}$ is clear.
Conversely, we compute:
 \begin{align*}
 \bv{X\ci Y}&=\MM_{x\in X}\JJ_{y\in Y}\bv{x=y}\\
 &=\JJ_{\gn\colon X\to Y}\MM_{x\in X}\bv{x=\gn(x)},
 \end{align*}
so, to conclude the proof, it suffices to prove that for every map
$\gn\colon X\to Y$, there exists $Z$ such that $\es\sci Z\ci Y$ and
 \begin{equation}\label{Eq:xXxgnxXZ}
\MM_{x\in X}\bv{x=\gn(x)}\leq\bv{X=Z}.
 \end{equation}
We define $Z$ as the range of~$\gn$. So,
 \[
 \MM_{x\in X}\bv{x=\gn(x)}\leq\MM_{x\in X}\bv{x\in Z}=\bv{X\ci Z}.
 \]
Furthermore, if $z\in Z$, so, $z=\gn(x^*)$ for some $x^*\in X$, then
 \[
 \MM_{x\in X}\bv{x=\gn(x)}\leq\bv{x^*=\gn(x^*)}\leq\bv{z\in X},
 \]
thus $\MM_{x\in X}\bv{x=\gn(x)}\leq\bv{Z\ci X}$, so, finally,
\eqref{Eq:xXxgnxXZ} holds. This concludes the proof.
\end{proof}

Every $D$-valued poset $P$ can be ``localized'' at every
prime filter of~$D$, in a classical fashion that we shall recall here.
Let $G$ be any \emph{filter} of~$D$, that is, a nonempty upper subset
of~$D$ closed under finite meet. We define binary relations,
$\leq_G$ and $\equiv_G$, on~$P$, by the rule
 \begin{align*}
 a\leq_Gb&\Longleftrightarrow\bv{a\leq b}\in G,\\
 a\equiv_Gb&\Longleftrightarrow\bv{a=b}\in G,
 \end{align*}
for all $a$, $b\in P$. It is easy to verify that the relation $\leq_G$ is a
preordering on~$P$, and that $\equiv_G$ is the associated equivalence
relation. Hence, the quotient structure $P/G=\vv<P,\leq_G>/{\equiv_G}$ may
be endowed with a partial ordering, defined by the rule
 \[
 \cls(a,G)\leq\cls(b,G)\Longleftrightarrow a\leq_Gb,
 \]
for all $a$, $b\in P$, where we write, of course,
$\cls(a,G)=\cls(a,\equiv_G)$.

The abundance of prime filters may be recorded in the following classical
result, that we shall use most of the time without mentioning:

\begin{lemma}\label{L:ManyG's}
Let $a$, $b\in D$. Then $a\leq b$ if{f} $a\in G$ implies that $b\in G$ for
all prime filters $G$ of~$D$.
\end{lemma}

As a rule, handling $D$-valued posets is very similar to handling
Boolean-valued posets. We point out two important differences with the
classical context:

\begin{itemize}
\item[---] The ``value set'' $D$ is no longer a complete Boolean algebra as
it is usually the case in the theory of Boolean-valued models. It is only a
distributive lattice, not even necessarily complete.

\item[---] No analogue of ``fullness'', as it is ordinarily defined for
Boolean models, will be assumed or even considered throughout this paper.
\end{itemize}

\section{$D$-valued partial lattices}\label{S:DvalPos}

\begin{definition}\label{D:DvalPL}
A \emph{$D$-valued partial lattice} is a $D$-valued poset $P$, endowed with
two maps from $P\times\fine P\to D$, denoted respectively by
$\vv<a,X>\mapsto\bv{a=\JJ X}$ and $\vv<a,X>\mapsto\bv{a=\MM X}$, such that
for all $a$, $b\in P$ and all $X$, $Y\in\fine P$, the following equalities
hold:
\begin{itemize}
\item[(1)]
$\bv{a=\JJ X}\mm\bv{a\leq b}=\bv{a=\JJ X}\mm\MM_{x\in X}\bv{x\leq b}$;

\item[(1*)]
$\bv{a=\MM X}\mm\bv{b\leq a}=\bv{a=\MM X}\mm\MM_{x\in X}\bv{b\leq x}$;

\item[(2)] $\bv{a=\JJ X}\mm\bv{X=Y}\leq\bv{a=\JJ Y}$;

\item[(2*)] $\bv{a=\MM X}\mm\bv{X=Y}\leq\bv{a=\MM Y}$.

\item[(3)] $\bv{a=\JJ X}\mm\bv{a=b}\leq\bv{b=\JJ X}$;

\item[(3*)] $\bv{a=\MM X}\mm\bv{a=b}\leq\bv{b=\MM X}$.
\end{itemize}

\end{definition}

\begin{example}\label{Ex:01ValPL}
Every partial lattice $P$ can be viewed as a $\two$-valued poset, as in
Example~\ref{Ex:01ValPos}. This structure can be extended to a structure of
$\two$-valued partial lattice, by putting
 \[
 \bv{a=\JJ X}=1\text{ if }a=\JJ X,\ 0\text{ otherwise},
 \]
and similarly for $\MM$.
\end{example}

\emph{For the remainder of this section, we shall fix a $D$-valued partial
lattice $P$}.

\begin{lemma}\label{L:DsiUnique}
Let $a$, $b\in P$ and let $X\in\fine P$. Then the following assertions hold:
\begin{enumerate}
\item $\bv{a=\JJ X}\mm\bv{b=\JJ X}\leq\bv{a=b}$;

\item $\bv{a=\MM X}\mm\bv{b=\MM X}\leq\bv{a=b}$.
\end{enumerate}
\end{lemma}

\begin{proof}
We only prove (i). Put $\gg=\bv{a=\JJ X}\mm\bv{b=\JJ X}$. By
(1) of Definition~\ref{D:DvalPL},
 \[
 \bv{b=\JJ X}\mm\MM_{x\in X}\bv{x\leq b}=\bv{b=\JJ X}\mm\bv{b\leq b}
 =\bv{b=\JJ X},
 \]
thus $\gg\leq\bv{b=\JJ X}\leq\MM_{x\in X}\bv{x\leq b}$.
Furthermore,
 \begin{align*}
 \gg\mm\bv{a\leq b}&=\gg\mm\MM_{x\in X}\bv{x\leq b}&&
 \text{(by (1) of Definition~\ref{D:DvalPL})}\\
 &=\gg&&\text{(by the above paragraph),}
 \end{align*}
so $\gg\leq\bv{a\leq b}$. Symmetrically, $\gg\leq\bv{b\leq a}$, so the
conclusion follows.
\end{proof}

If $G$ is a filter of~$D$, we have seen that we can define a
quotient poset $P/G$. We shall now show how to extend the structure of~$P/G$
to a structure of partial lattice.

\begin{definition}\label{D:P/GPartLatt}
Let $\xX\in\fine{P/G}$ and let $\xa\in P/G$. We define $\xa=\JJ\xX$ (resp.,
$\xa=\MM\xX$) to hold, if there are $a\in P$ and $X\in\fine P$ such that
$\xa=\cls(a,G)$, $\xX=X/G$, and $\bv{a=\JJ X}\in G$ (resp.,
$\bv{a=\MM X}\in G$).
\end{definition}

As an immediate consequence of Definition~\ref{D:DvalPL}(2,2*,3,3*), we
obtain the following lemma:

\begin{lemma}\label{L:P/GfromP,G}
Let $a\in P$, let $X\in\fine P$. Then $\cls(a,G)=\JJ X/G$
(resp., $\cls(a,G)=\MM X/G$) if{f} $\bv{a=\JJ X}\in G$
(resp., $\bv{a=\MM X}\in G$).
\end{lemma}

\begin{proposition}\label{P:P/GPartLatt}
The poset $P/G$, endowed with $\JJ$ and $\MM$ of
Definition~\tup{\ref{D:P/GPartLatt}}, is a partial lattice.
\end{proposition}

\begin{proof}
We first have to prove that $\JJ$ and $\MM$ are \emph{functions}. We do it
for $\JJ$. So let $\xX\in\fine{P/G}$ and let $\xa$, $\xb\in P/G$ such that
$\xa=\JJ\xX$ and $\xb=\JJ\xX$. Let $a$, $b\in P$ and let $X\in\fine P$ such
that $\xa=\cls(a,G)$, $\xb=\cls(b,G)$, and $\xX=X/G$. By
Lemma~\ref{L:P/GfromP,G}, both $\bv{a=\JJ X}$ and $\bv{b=\JJ X}$ belong to
$G$, hence, by Lemma~\ref{L:DsiUnique}, $\bv{a=b}\in G$, so $\xa=\xb$. Hence
$\JJ$ is a function on $P/G$. The same argument applies to $\MM$.

To conclude the proof, it is sufficient to prove that for $\xa\in P/G$
and $\xX\in\fine{P/G}$, $\xa=\JJ\xX$ implies that $\xa=\sup\xX$ (for the
partial ordering of~$P/G$), and similarly for $\MM$. We present the proof
for $\JJ$. Let $a\in P$ and $X\in\fine P$ such that
$\xa=\cls(a,G)$ and $\xX=X/G$. By Lemma~\ref{L:P/GfromP,G},
$\bv{a=\JJ X}\in G$. For $x\in X$, it follows from
Definition~\ref{D:DvalPL}(1) that
 \[
 \bv{a=\JJ X}=\bv{a=\JJ X}\mm\bv{a\leq a}\leq\bv{a=\JJ X}\mm\bv{x\leq a}
 \leq\bv{x\leq a},
 \]
so $\bv{x\leq a}\in G$, that is, $\cls(x,G)\leq\cls(a,G)=\xa$. So,
$\xX\leq\xa$. Now let $\xb\in P/G$ such that
$\xX\leq\xb$. Pick $b\in\xb$. For $x\in X$,
$\cls(x,G)\leq\xb=\cls(b,G)$, so $\bv{x\leq b}\in G$; hence
$\MM_{x\in X}\bv{x\leq b}\in G$. By Definition~\ref{D:DvalPL}(1),
 \[
 \bv{a=\JJ X}\mm\bv{a\leq b}=
 \bv{a=\JJ X}\mm\MM_{x\in X}\bv{x\leq b}\in G,
 \]
hence $\bv{a\leq b}\in G$, that is, $\xa\leq\xb$. So we have proved that
$\xa=\sup\xX$.
\end{proof}

\section{Join-samples and meet-samples}\label{S:JMSample}

Let $D$ be a distributive lattice with unit, let $P$
be a $D$-valued partial lattice.

We introduce one of the most important definitions of the whole paper:

\begin{definition}\label{D:MJsample}
Let $X$ be a nonempty finite subset of~$P$. A \emph{join-sample} (resp.,
\emph{meet-sample}) of~$X$ is a nonempty finite subset $U$ of~$P$ such that
 \begin{align*}
 \bv{x=\JJ X}&\leq\JJ_{u\in U}\bv{u=\JJ X},\q\text{for all }x\in P\\
 \biggl(\text{resp., }
 \bv{x=\MM X}&\leq\JJ_{u\in U}\bv{u=\MM X},\q\text{for all }x\in P\biggr).
 \end{align*}
\end{definition}

\begin{definition}\label{D:FinSam}
A $D$-valued partial lattice $P$ is \emph{finitely join-sampled} (resp.,
\emph{finitely meet-sampled}), if every nonempty finite subset of~$P$ has a
join-sample (resp., a meet-sample). We say that $P$ is \emph{finitely
sampled}, if it is both finitely join-sampled and finitely meet-sampled.
\end{definition}

Of course, if $U$ is a join-sample of~$X$ and $V$ is a meet-sample of~$X$,
then $U\uu V$ (or anything larger) is both a join-sample and a meet-sample
of~$X$.

\begin{lemma}\label{L:aleqjjX}
Let $X$, $U$, $V\in\fine P$.
\begin{enumerate}
\item If $U$ and $V$ are join-samples of~$X$, then the equality
 \[
 \JJ_{u\in U}\bv{a\leq u}\mm\bv{u=\JJ X}=
 \JJ_{v\in V}\bv{a\leq v}\mm\bv{v=\JJ X}
 \]
holds, for all $a\in P$.
\item If $U$ and $V$ are meet-samples of~$X$, then the equality
 \[
 \JJ_{u\in U}\bv{u\leq a}\mm\bv{u=\MM X}=
 \JJ_{v\in V}\bv{v\leq a}\mm\bv{v=\MM X}
 \]
holds, for all $a\in P$.
\end{enumerate}
\end{lemma}

\begin{proof}
We provide a proof for (i); (ii) is dual.
For $u\in U$,
\begin{align*}
\bv{a\leq u}\mm\bv{u=\JJ X}&=
\JJ_{v\in V}\left(\bv{a\leq u}\mm\bv{u=\JJ X}\mm\bv{v=\JJ X}\right)\\
\intertext{(because $V$ is a join-sample of~$X$)}
&\leq\JJ_{v\in V}\left(\bv{a\leq u}\mm\bv{u=v}\mm\bv{v=\JJ X}\right)\\
\intertext{(by Lemma~\ref{L:DsiUnique})}
&\leq\JJ_{v\in V}\bv{a\leq v}\mm\bv{v=\JJ X}.
\end{align*}
hence $\JJ_{u\in U}\bv{a\leq u}\mm\bv{u=\JJ X}\leq
\JJ_{v\in V}\bv{a\leq v}\mm\bv{v=\JJ X}$. The proof of the converse
inequality is similar.
\end{proof}

Lemma~\ref{L:aleqjjX} makes it possible to define, for all $a\in P$ and all
$X\in\fine P$,
 \begin{align*}
 \bv{a\leq\JJ X}&=\JJ_{u\in U}\bv{a\leq u}\mm\bv{u=\JJ X},
 &&\text{for every join-sample }U\text{ of }X,\\
 \bv{\MM X\leq a}&=\JJ_{u\in U}\bv{u\leq a}\mm\bv{u=\MM X},
 &&\text{for every meet-sample }U\text{ of }X.
 \end{align*}
We recall that a filter $G$ of~$D$ is \emph{prime}, if $x\jj y\in G$ implies
that $x\in G$ or $y\in G$, for all $x$, $y\in D$.

\begin{lemma}[The Basic Truth Lemma]\label{L:TruthLemma}
Assume that $P$ is finitely sampled.
Let $\gf(z,Z)$ be one of the following formulas:
\begin{itemize}
\item $z=\JJ Z$;

\item $z\leq\JJ Z$;

\item $z=\MM Z$;

\item $\MM Z\leq z$.

\end{itemize}
Let $a\in P$, let $X\in\fine P$, let $G$ be a prime filter of~$D$.
Then the following equivalence holds:
 \[
 P/G\text{ satisfies }\gf(\cls(a,G),X/G)\q\text{if{f}}\q
 \bv{\gf(a,X)}\in G.
 \]
\end{lemma}

\begin{proof}
By duality, it is sufficient to prove the result in case $\gf(z,Z)$ is
either $z=\JJ Z$ or $z\leq\JJ Z$. The first case follows from
Lemma~\ref{L:P/GfromP,G}. So, suppose that $\gf(z,Z)$ is $z\leq\JJ Z$.

Let $U$ be a join-sample of~$X$. Suppose first that
$\cls(a,G)\leq\JJ X/G$ (in $P/G$). In particular, $\JJ(X/G)$ is
defined, so, by Definition~\ref{D:P/GPartLatt} and by
Lemma~\ref{L:P/GfromP,G}, there exists
$a'\in P$ such that
\begin{align}
\bv{a'=\JJ X}&\in G\label{Eq:a'JXG},\\
\bv{a\leq a'}&\in G\label{Eq:alea'G}.
\end{align}
Since $U$ is a join-sample of~$X$,
$\bv{a'=\JJ X}\leq\JJ_{u\in U}\bv{u=\JJ X}$, thus, since $G$ is prime and
$U$ is finite, there exists, by \eqref{Eq:a'JXG}, $u\in U$ such that
\begin{equation}\label{Eq:uJXG}
\bv{u=\JJ X}\in G.
\end{equation}
{F}rom \eqref{Eq:a'JXG} and Lemma~\ref{L:DsiUnique}, it follows that
$\bv{a'=u}\in G$, thus, by \eqref{Eq:alea'G}, $\bv{a\leq u}\in G$. Hence, by
\eqref{Eq:uJXG}, $\bv{a\leq\JJ X}\in G$.

Conversely, suppose that $\bv{a\leq\JJ X}\in G$. Since $U$ is finite
and $G$ is prime, there exists $u\in U$ such that $\bv{a\leq u}\in G$ and
$\bv{u=\JJ X}\in G$. Hence, by Lemma~\ref{L:P/GfromP,G},
$\cls(a,G)\leq\cls(u,G)$ and $\cls(u,G)=\JJ X/G$, so
$\cls(a,G)\leq\JJ X/G$.
\end{proof}

Further analogues of Lemma~\ref{L:TruthLemma} will be met in
\ref{L:TruthIdFiln}, \ref{L:TruthIdFil}, \ref{P:HomIdFilP/G},
\ref{C:xllyClassBool}, \ref{P:xleqyClassBool}, \ref{L:PDmeasP/G}.

\section{Ideal and filter samples}\label{S:IdFilSamples}

In this section, we fix a distributive lattice $D$ with unit and a
$D$-valued partial lattice $P$.

\begin{definition}\label{D:IDFIMsample}
Let $X\in\fine P$. An \emph{\IDM-sample} of~$X$ is an element $U$ of
$\fine P$ such that
 \begin{equation}\label{Eq:IdmSample}
 \MM_{x\in X}\bv{a\leq x}=
 \JJ_{u\in U}\left(\bv{a\leq u}\mm\MM_{x\in X}\bv{u\leq x}\right)
 \end{equation}
holds for all $a\in P$.

Dually, a \emph{\FIM-sample} of~$X$ is an element $U$ of~$\fine P$
such that
 \begin{equation}\label{Eq:FilmSample}
 \MM_{x\in X}\bv{x\leq a}=
 \JJ_{u\in U}\left(\bv{u\leq a}\mm\MM_{x\in X}\bv{x\leq u}\right)
 \end{equation}
holds for all $a\in P$.
\end{definition}

We observe that the $\geq$ half of both equalities
\eqref{Eq:IdmSample} and \eqref{Eq:FilmSample} always holds, thus it is
sufficient to verify the $\leq$ half. As a consequence of this, we observe
that \emph{every finite subset of $P$ that contains an \IDM-sample of $X$
is an \IDM-sample of $X$}.

\begin{definition}\label{D:IDFIM}
We say that \emph{$P$ has \IDM} (resp., \emph{\FIM}), if
every \emph{pair} of elements of~$P$ has an \IDM-sample (resp., a
\FIM-sample).
\end{definition}

We state without proof the following easy result, that will not be used
later:

\begin{proposition}\label{P:2tonIDFI}
If $P$ has \IDM\ (resp., \FIM), then every nonempty finite
subset of~$P$ has an \IDM-sample (resp., a \FIM-sample).
\end{proposition}

The last two finiteness properties about $P$ that we shall consider are
harder to define. To prepare for this task, we first define new $D$-valued
functions on~$P$.

\begin{definition}\label{D:BVainId}
Suppose that $P$ is finitely join-sampled.
For $a\in P$, for nonempty, finite subsets $X$ and $U$ of~$P$, and for
$n<\go$, we define an element $\bv{a\in\Id_n(X,U)}$ of~$D$, by induction on
$n$, as follows:
\begin{itemize}
\item[(i)] $\bv{a\in\Id_0(X,U)}=\bv{a\in\dnw X}=\JJ_{x\in X}\bv{a\leq x}$.

\item[(ii)] The induction step:
 \[
 \bv{a\in\Id_{n+1}(X,U)}=\bv{a\in\Id_n(X,U)}\jj
 \JJ_{\es\sci Z\ci U}\bv{a\leq\JJ Z}\mm\bv{Z\ci\Id_n(X,U)},
 \]
where we put $\bv{Z\ci\Id_n(X,U)}=\MM_{z\in Z}\bv{z\in\Id_n(X,U)}$.
\end{itemize}

Dually, suppose that $P$ is finitely meet-sampled.
For $a\in P$, for nonempty, finite subsets $X$ and $U$ of~$P$, and for
$n<\go$, we define an element $\bv{a\in\Id_n(X,U)}$ of~$D$, by induction on
$n$, as follows:
\begin{itemize}
\item[(i*)] $\bv{a\in\Fil_0(X,U)}=\bv{a\in\upw X}=\JJ_{x\in X}\bv{x\leq a}$.

\item[(ii*)] The induction step:
 \[
 \bv{a\in\Fil_{n+1}(X,U)}=\bv{a\in\Fil_n(X,U)}\jj
 \JJ_{\es\sci Z\ci U}\bv{\MM Z\leq a}\mm\bv{Z\ci\Fil_n(X,U)},
 \]
where we put $\bv{Z\ci\Fil_n(X,U)}=\MM_{z\in Z}\bv{z\in\Fil_n(X,U)}$.
\end{itemize}

\end{definition}

We observe that the condition that $P$ be finitely join- or meet-sampled is
necessary in order to \emph{define} the elements $\bv{a\in\Id_n(X,U)}$ and
$\bv{a\in\Fil_n(X,U)}$, since the elements $\bv{a\leq\JJ Z}$ and
$\bv{\MM Z\leq a}$ need to be defined. Our next result relates the
$D$-valued $\Id_n$ and $\Fil_n$ with their corresponding classical versions,
see Definition~\ref{D:IdFiln01}.

\begin{lemma}[Truth Lemma for $\Id_n(X,U)$ and $\Fil_n(X,U)$]
\label{L:TruthIdFiln}
Let $a\in P$, let $X$, $U\in\fine P$, and
let $G$ be a prime filter of~$D$. Then the following assertions hold:
\begin{enumerate}
\item Suppose that $P$ is finitely join-sampled.
Then $\cls(a,G)\in\Id_n(X/G,U/G)$ in $P/G$ if{f}
$\bv{a\in\Id_n(X,U)}\in G$, for any $n<\go$.

\item Suppose that $P$ is finitely meet-sampled.
Then $\cls(a,G)\in\Fil_n(X/G,U/G)$ in $P/G$ if{f}
$\bv{a\in\Fil_n(X,U)}\in G$, for any $n<\go$.
\end{enumerate}

\end{lemma}

\begin{proof}
We provide a proof for (i). We argue by induction on $n$.
The result for $n=0$ follows immediately from the finiteness of~$X$ and the
fact that $G$ is prime.

Now suppose the statement proved for $n$, we prove it for $n+1$.
Suppose first that $\bv{a\in\Id_{n+1}(X,U)}\in G$. Since $\mathcal{P}(U)$ is
finite and since $G$ is prime, either
$\bv{a\nobreak\in\nobreak\Id_n(X,U)}\in G$, or there exists a nonempty
finite subset $Z$ of~$U$ such that
$\bv{a\nobreak\leq\nobreak\JJ\nobreak Z}\in G$ and $\bv{Z\ci\Id_n(X,U)}\in G$.
In the first case, if follows from the induction hypothesis that
$\cls(a,G)\in\Id_n(X/G,U/G)\ci\Id_{n+1}(X/G,U/G)$,
so we are done. In the second case, $\cls(a,G)\leq\JJ(Z/G)$
by Lemma~\ref{L:TruthLemma}, $Z/G\ci\Id_n(X/G,U/G)$ by the induction
hypothesis, and $\es\sci Z/G\ci U/G$, hence,
$\cls(a,G)\in\Id_{n+1}(X/G,U/G)$.

Conversely, suppose that $\cls(a,G)\in\Id_{n+1}(X/G,U/G)$. If
$\cls(a,G)\in\Id_n(X/G,U/G)$, then, by the induction
hypothesis, $\bv{a\in\Id_n(X,U)}\in G$, hence $\bv{a\in\Id_{n+1}(X,U)}\in G$.
Otherwise, there exists a nonempty $\xZ\ci U/G$ such that
$\cls(a,G)\leq\JJ\xZ$ and $\xZ\ci\Id_n(X/G,U/G)$. Since
$\es\sci\xZ\ci U/G$, there exists a nonempty subset $Z$ of~$U$ such
that $\xZ=Z/G$. So $\cls(a,G)\leq\JJ Z/G$, thus, by
Lemma~\ref{L:TruthLemma}, $\bv{a\leq\JJ Z}\in G$. Since
$Z/G=\xZ\ci\Id_n(X/G,U/G)$, it follows from the induction
hypothesis that $\bv{Z\ci\Id_n(X,U)}\in G$. Since $\es\sci Z\ci U$, $Z$
witnesses the fact that $\bv{a\in\Id_{n+1}(X,U)}\in G$.
\end{proof}

\begin{definition}\label{D:IDFIJsample}
Let $X$ be a nonempty finite subset of~$P$. An \emph{\IDJ-sample} (resp.,
\emph{\FIJ-sample}) of~$X$ is a nonempty finite subset $U$ of~$P$ such that
there exists $n<\go$ such that
\begin{align*}
\bv{a\in\Id_n(X,U)}&=\bv{a\in\Id_{n+1}(X,Y)}\\
(\text{resp., }\bv{a\in\Fil_n(X,U)}&=\bv{a\in\Fil_{n+1}(X,Y)}),
\end{align*}
for all $a\in P$ and all $Y\in\fine P$ containing $U$.

We call any such $n$ an \emph{ideal index} (resp., \emph{filter index}) of
$\vv<X,U>$.
\end{definition}

If $U$ is an \IDJ-sample of~$X$, with ideal index $n$, then it is easy to
verify that $\bv{a\in\Id_k(X,Y)}=\bv{a\in\Id_n(X,U)}$, for all $a\in P$,
all $k\geq n$, and all finite $Y\ce U$. Hence this expression is independent
of the chosen sample $U$ and index $n$, we denote it by $\bv{a\in\Id(X)}$.

Dually, we define $\bv{a\in\Fil(X)}$ as the common value of
$a\in\Fil_n(X,U)$, for every \FIJ-sample $U$ of~$X$, with filter index $n$.

\begin{definition}\label{D:IDFIJ}
We say that $P$ \emph{has \IDJ} (resp., \emph{\FIJ}), if
$P$ is finitely join-sampled (resp., finitely meet-sampled) and every
nonempty finite subset of~$P$ has an \IDJ-sample (resp., a \FIJ-sample).
\end{definition}

As an easy consequence of the remarks following Definition~\ref{D:IdFiln01}
and of Lemma~\ref{L:TruthIdFiln}, we obtain the following:

\begin{lemma}[Truth Lemma for $\Id(X)$ and $\Fil(X)$]
\label{L:TruthIdFil}
Let $a\in P$, let $X\in\fine P$, and
let $G$ be a prime filter of~$D$. Then the following equivalences hold:
\begin{enumerate}
\item Suppose that $P$ has \IDJ. Then
 \[
\cls(a,G)\in\Id(X/G)\text{ in }P/G\text{ if{f} }
\bv{a\in\Id(X)}\in G.
 \]

\item Suppose that $P$ has \FIJ. Then
 \[
\cls(a,G)\in\Fil(X/G)\text{ in }P/G\text{ if{f} }
\bv{a\in\Fil(X)}\in G.
 \]
\end{enumerate}
\end{lemma}

The definition of \IDJ\ and of \FIJ\ for $D$-valued partial
lattice presented in Definition~\ref{D:IDFIJ} is quite unwieldy, because it
involves the Boolean values $\bv{a\in\Id_n(X,U)}$ or $\bv{a\in\Fil_n(X,U)}$
presented in Definition~\ref{D:IDFIJsample}. However,
Lemma~\ref{L:TruthIdFiln} makes it possible to find a useful equivalent
form:

\begin{lemma}\label{L:EquivIdFilSG}
The finitely join-sampled $D$-valued partial lattice $P$ has \IDJ\ if{f}
for all $X\in\fine P$, there are $U\in\fine P$ and $n<\go$ such that
 \[
 \Id_n(X/G,U/G)=\Id_{n+1}(X/G,Y/G),
 \]
for every $Y\in\fine P$ containing $U$ and every prime filter $G$ of~$D$.
The dual statement holds, about \FIJ\ and $\Fil_n$, for finitely
meet-sampled $P$.
\end{lemma}

\section{Affine lower and upper functions; ideal and filter functions}

Let $D$ be a distributive lattice with unit. The $D$-valued analogue of the
notions of a lower set and an upper set are provided by the following
definition.

\begin{definition}\label{D:LoUpFct}
Let $P$ be a $D$-valued poset.
A map $f\colon P\to D$ is a \emph{lower function}, if
$f(y)\mm\bv{x\leq y}\leq f(x)$, for all $x$, $y\in P$. Dually, $f$ is an
\emph{upper function}, if $f(x)\mm\bv{x\leq y}\leq f(y)$, for all $x$,
$y\in P$.
\end{definition}

For example, if $P$ is a poset, viewed, as in Example~\ref{Ex:01ValPos},
with its canonical structure of~$\two$-valued poset, then the lower (resp.,
upper) functions on~$P$ are exactly the characteristic functions of the
lower (resp., upper) subsets of~$P$.

It is obvious that for $a\in P$, the map $x\mapsto\bv{x\leq a}$ (resp.,
$x\mapsto\bv{a\leq x}$) is a lower function (resp., upper function) on
$P$---we shall call these functions \emph{principal lower functions} (resp.,
\emph{principal upper functions}). Furthermore, any constant function is both
a lower function and an upper function, and any finite meet or join of lower
functions (resp., upper functions) is a lower function (resp., an upper
function). This gives a class of ``simple'' lower functions and upper
functions, an analogue of finitely generated lower subsets or upper subsets of
a poset.

\begin{definition}\label{D:AffLoUp}
Let $P$ be a $D$-valued poset. An \emph{affine lower function} on~$P$ is a
map $f\colon P\to D$ defined by a rule of the form
 \[
 f(x)=\JJ_{i<n}\bv{x\leq u_i}\mm\ga_i,\qq\text{for all }x\in P,
 \]
where $n\in\gos$, $u_0$,\dots, $u_{n-1}\in P$, and
$\ga_0$,\dots, $\ga_{n-1}\in D$.
Dually, an \emph{affine upper function} on~$P$ is a
map $f\colon P\to D$ defined by a rule of the form
 \[
 f(x)=\JJ_{i<n}\bv{u_i\leq x}\mm\ga_i,\qq\text{for all }x\in P,
 \]
where $n\in\gos$, $u_0$,\dots, $u_{n-1}\in P$, and
$\ga_0$,\dots, $\ga_{n-1}\in D$.
\end{definition}

In particular, any affine lower function is a lower function, and any affine
upper function is an upper function.

\begin{definition}\label{D:IdFilFct}
Let $P$ be a $D$-valued partial lattice. An \emph{ideal function} on~$P$ is
a lower function $f\colon P\to D$ such that
 \[
 \bv{a=\JJ X}\mm\MM_{x\in X}f(x)\leq f(a),\qq\text{for all }a\in P
 \text{ and all }X\in\fine P.
 \]
Dually, a \emph{filter function} on~$P$ is an upper function
$f\colon P\to D$ such that
 \[
 \bv{a=\MM X}\mm\MM_{x\in X}f(x)\leq f(a),\qq\text{for all }a\in P
 \text{ and all }X\in\fine P.
 \]
An \emph{affine ideal function} on $P$ is a function that is simultaneously
an affine lower function and an ideal function on $P$.
Dually, an \emph{affine filter function} on $P$ is a function that is
simultaneously an affine upper function and a filter function on $P$.
\end{definition}

We observe that the set of all ideal functions (resp., filter functions) on
$P$ is closed under componentwise meet, but not under componentwise join as
a rule, just the same way as the \emph{union} of two ideals of a partial
lattice is not necessarily an ideal.

\begin{example}
Every partial lattice $P$ can be viewed as a $\two$-valued partial lattice,
see Example~\ref{Ex:01ValPL}. If $I$ is an ideal of~$P$, then the
characteristic function of~$I$ is an ideal function on~$P$, and, dually, a
similar statement holds for filters.
\end{example}

\begin{lemma}\label{L:IdFilFct}
Let $P$ be a $D$-valued partial lattice, let $f\colon P\to D$.
\begin{enumerate}
\item If $P$ is finitely join-sampled and $f$ is an ideal function, then
 \[
 \bv{a\leq\JJ X}\mm\MM_{x\in X}f(x)\leq f(a),\qq\text{for all }a\in P
 \text{ and all }X\in\fine P.
 \]
\item If $P$ is finitely meet-sampled and $f$ is a filter function, then
 \[
 \bv{\MM X\leq a}\mm\MM_{x\in X}f(x)\leq f(a),\qq\text{for all }a\in P
 \text{ and all }X\in\fine P.
 \]
\end{enumerate}
\end{lemma}

\begin{proof}
We provide a proof for (i). Let $U$ be a join-sample of
$X$. Then, for $a\in P$,
 \begin{align*}
 \bv{a\leq\JJ X}\mm\MM_{x\in X}f(x)&=
 \JJ_{u\in U}\left(\bv{a\leq u}\mm\bv{u=\JJ X}\mm\MM_{x\in X}f(x)\right)\\
 &\leq\JJ_{u\in U}\bv{a\leq u}\mm f(u)\\
\intertext{(because $f$ is an ideal function)}
 &\leq f(a), 
 \end{align*}
(because $f$ is a lower function).
\end{proof}

In a $D$-valued partial lattice $P$, it is easy to prove that any principal
lower function is an affine ideal function, and any principal upper function
is an affine filter function. Our next result provides extensions of this
simple fact.

\begin{proposition}\label{P:IdFilnAff}
Let $P$ be a $D$-valued partial lattice, let $X$ and $U$ be
nonempty finite subsets of~$P$. Then the following assertions hold:

\begin{enumerate}
\item Suppose that $P$ is finitely join-sampled. Then the map
$a\mapsto\bv{a\in\Id_n(X,U)}$ is an affine lower function on~$P$.
Furthermore, if $P$ has \IDJ, then the map $a\mapsto\bv{a\in\Id(X)}$
is an affine ideal function on~$P$.

\item Suppose that $P$ is finitely meet-sampled. Then the map
$a\mapsto\bv{a\in\Fil_n(X,U)}$ is an affine upper function on~$P$.
Furthermore, if $P$ has \FIJ, then the map $a\mapsto\bv{a\in\Fil(X)}$
is an affine filter function on~$P$.
\end{enumerate}
\end{proposition}

\begin{proof}
We provide a proof for (i). For $n<\go$, let
$f_n\colon a\mapsto\bv{a\in\Id_n(X,U)}$. We prove, by induction on $n$, that
$f_n$ is an affine lower function.

For $n=0$, $f_0(a)=\JJ_{x\in X}\bv{a\leq x}$ for all $a$, thus $f_0$ is an
affine lower function.

Before proceeding to the induction step, we prove a claim:

\setcounter{claim}{0}
\begin{claim}
The map $a\mapsto\bv{a\leq\JJ Y}$ is an affine lower function, for all
$Y\in\fine P$.
\end{claim}

\begin{cproof}
Let $V$ be a join-sample of~$Y$. Then
 \[
 \bv{a\leq\JJ Y}=\JJ_{v\in V}\bv{a\leq v}\mm\gb_v,\q\text{for all }a\in P
 \q\text{(use Lemma~\ref{L:aleqjjX}(i))},
 \]
where we put $\gb_v=\bv{v=\JJ Y}$, for all $v\in V$.
\end{cproof}

Now suppose that $f_n$ is an affine lower function on
$P$. Then
 \[
 f_{n+1}(a)=f_n(a)\jj\JJ_{\es\sci Z\ci U}\bv{a\leq\JJ Z}\mm\gg_Z,
 \]
for all $a\in P$, where we put $\gg_Z=\bv{Z\ci\Id_n(X,U)}$, for all nonempty
$Z\ci U$. Therefore, by Claim~1 and the induction hypothesis, $f_{n+1}$ is an
affine lower function. So all $f_n$ are affine lower functions.

Now let $f\colon a\mapsto\bv{a\in\Id(X)}$.
Suppose that $P$ has \IDJ. Let $U$ be an \IDJ-sample of~$X$, with index
$n$. So $\bv{a\in\Id(X)}=\bv{a\in\Id_n(X,U)}$, for all $a\in P$. Hence $f$
is an affine lower function. So, to conclude the proof, it is sufficient to
prove that $f$ is an ideal function. So let $Z\in\fine P$, let $a\in P$. We
compute:
 \begin{align*}
 \bv{a=\JJ Z}\mm\MM_{z\in Z}f(z)
 &=\bv{a=\JJ Z}\mm\MM_{z\in Z}\bv{z\in\Id_n(X,U)}\\
 &=\bv{a=\JJ Z}\mm\bv{Z\ci\Id_n(X,U)}\\
 &\leq\bv{a=\JJ Z}\mm\bv{Z\ci\Id_n(X,U\uu Z)}\\
 &\leq\bv{a\in\Id_{n+1}(X,U\uu Z)}\\
\intertext{(because $\bv{a=\JJ Z}\leq\bv{a\leq\JJ Z}$)}
 &=f(a).
 \end{align*}
So $f$ is an ideal function on~$P$.
\end{proof}

\begin{notation}
Let $P$ be a $D$-valued partial lattice. We denote by $\aID(P)$ (resp.,
$\aFIL(P)$) the set of all affine ideal functions (resp., affine filter
functions) on~$P$, partially ordered componentwise.
\end{notation}

For the remainder of this section, we assume that $P$ is a
$D$-valued partial lattice.

\begin{notation}
Let $f\colon P\to D$. If there exists a least ideal (resp., filter) function
$g$ such that $f\leq g$, then we denote this function by $f^{\Id}$ (resp.,
$f^{\Fil}$).
\end{notation}

We observe that $f\leq f^{\Id}$ and also $f\leq f^{\Fil}$.

\begin{lemma}\label{L:(dnwXmma)*}
Let $X\in\fine P$.
\begin{enumerate}
\item Suppose that $P$ has \IDJ. Let
$f\colon\nobreak a\mapsto\bv{a\nobreak\in\nobreak\dnw X}\mm\ga$.
Then $f^{\Id}\colon a\mapsto\bv{a\nobreak\in\nobreak\Id(X)}\mm\ga$.

\item Suppose that $P$ has \FIJ. Let
$f\colon\nobreak a\mapsto\bv{a\nobreak\in\nobreak\upw X}\mm\ga$.
Then $f^{\Fil}\colon a\mapsto\bv{a\nobreak\in\nobreak\Fil(X)}\mm\ga$.
\end{enumerate}
\end{lemma}

\begin{proof}
We provide a proof for (i). Let
$g\colon a\mapsto\bv{a\in\Id(X)}\mm\ga$. It is obvious that $f\leq g$. By
Proposition~\ref{P:IdFilnAff}, $g$ is an affine ideal function on~$P$.

It remains to prove that $g\leq h$, for every ideal function $h$ on~$P$ such
that $f\leq h$. Since $g$ has the form $a\mapsto\bv{a\in\Id_n(X,U)}\mm\ga$,
for some $U$ and some $n$, it suffices to prove that
 \begin{equation}\label{Eq:aIdnXuleh}
 \bv{a\in\Id_n(X,U)}\mm\ga\leq h(a),\q\text{for all }a\in P,
 \text{ all }U\in\fine P,\text{ and all }n<\go. 
 \end{equation}
For $n=0$, $\bv{a\in\Id_n(X,U)}\mm\ga=f(a)\leq h(a)$, so
\eqref{Eq:aIdnXuleh} holds. Assume that \eqref{Eq:aIdnXuleh} holds for $n$.
For nonempty $Z\ci U$, we compute:
 \begin{align*}
 \bv{a\leq\JJ Z}\mm\bv{Z\ci\Id_n(X,U)}\mm\ga
 &\leq\bv{a\leq\JJ Z}\mm\MM_{z\in Z}h(z)\\
\intertext{(by the induction hypothesis)}
 &\leq h(a)
 \end{align*}
by Lemma~\ref{L:IdFilFct}. Hence $\bv{a\in\Id_{n+1}(X,U)}\mm\ga\leq h(a)$,
for all $a\in P$. This concludes the proof of \eqref{Eq:aIdnXuleh}.
\end{proof}

As a consequence, $f^{\Id}$ can be computed explicitly, for any affine lower
function $f$ (and dually):

\begin{proposition}\label{P:Compf*}
Let $n\in\gos$, let $u_0$,\ldots,
$u_{n-1}\in P$, let $\ga_0$,\ldots, $\ga_{n-1}\in D$.
For all nonempty $I\ci n$, we put
 \begin{align*}
 u^{(I)}&=\setm{u_i}{i\in I},\\
 \ga_{(I)}&=\MM_{i\in I}\ga_i.
 \end{align*}
\begin{enumerate}
\item Suppose that $P$ has \IDJ.
Let $f\colon a\mapsto\JJ_{i<n}\bv{a\leq u_i}\mm\ga_i$. Then $f^{\Id}$ is
defined, and
 \[
 f^{\Id}(a)=\JJ_{\es\sci I\ci n}\bv{a\in\Id(u^{(I)})}\mm\ga_{(I)},
 \]
for all $a\in P$. In particular, $f^{\Id}$ is an affine ideal function.

\item Suppose that $P$ has \FIJ.
Let $f\colon a\mapsto\JJ_{i<n}\bv{u_i\leq a}\mm\ga_i$. Then $f^{\Fil}$ is
defined, and
 \[
 f^{\Fil}(a)=\JJ_{\es\sci I\ci n}\bv{a\in\Fil(u^{(I)})}\mm\ga_{(I)},
 \]
for all $a\in P$. In particular, $f^{\Fil}$ is an affine filter function.
\end{enumerate}
\end{proposition}

\begin{proof}
We provide a proof for (i). By Lemma~\ref{L:(dnwXmma)*}, for
$\es\sci I\ci n$, the map
$g_I\colon a\mapsto\bv{a\in\Id(u^{(I)})}\mm\ga_{(I)}$ is an
affine ideal function, so $g=\JJ_{\es\sci I\ci n}g_I$ is an affine lower
function on~$P$. Furthermore, $g_{\set{i}}(a)=\bv{a\leq u_i}\mm\ga_i$,
for all $i<n$ and all $a\in P$, thus $f\leq g$. Let $h$ be an ideal function
on~$P$ such that $f\leq h$. In order to verify that
$g\leq h$, it suffices to verify that $g_I\leq h$ for all nonempty $I\ci n$.
For $a\in P$,
 \[
 \bv{a\in\dnw u^{(I)}}\mm\ga_{(I)}
 =\JJ_{i\in I}\bv{a\leq u_i}\mm\ga_{(I)}\leq f(a)\leq h(a).
 \]
Therefore, by Lemma~\ref{L:(dnwXmma)*}, $g_I\leq h$. This holds for all $I$,
therefore, $g\leq h$.

To conclude the proof, it suffices to prove that $g$ is an ideal function on
$P$. So, let $a\in P$ and let $X\in\fine P$. We shall prove that
 \begin{equation}\label{Eq:gIdFct}
 \bv{a=\JJ X}\mm\MM_{x\in X}g(x)\leq g(a).
 \end{equation}
To prove \eqref{Eq:gIdFct}, it suffices to prove that
$\bv{a=\JJ X}\mm\MM_{x\in X}g(x)\in G$ implies that $g(a)\in G$, for any
prime filter $G$ of~$D$. By Lemmas \ref{L:TruthLemma}
and \ref{L:TruthIdFil} and by the definition of~$g$,
 \begin{equation}\label{Eq:a/GX/G}
 \cls(a,G)=\JJ X/G,
 \end{equation}
and, for all $x\in X$, there exists a nonempty $I_x\ci n$ such that
 \begin{align}
 \ga_{(I_x)}&\in G,\label{Eq:aIxG}\\
 \bv{x\in\Id(u^{(I_x)})}&\in G.\label{Eq:xuIxG}
 \end{align}
Now put $I=\UU_{x\in X}I_x$. Then, by \eqref{Eq:aIxG},
$\ga_{(I)}=\MM_{x\in X}\ga_{(I_x)}\in G$. Furthermore, $I_x\ci I$, hence
$\bv{x\in\Id(u^{(I_x)})}\leq\bv{x\in\Id(u^{(I)})}$, so, by \eqref{Eq:xuIxG},
$\bv{x\in\Id(u^{(I)})}\in G$, for all $x\in X$. So we have proved that
$\bv{x\in\Id(u^{(I)})}\mm\ga_{(I)}\in G$, for all $x\in X$. By
Lemma~\ref{L:TruthIdFil}, $\cls(x,G)\in\Id\left(u^{(I)}/G\right)$.
This holds for all $x\in X$, thus, by \eqref{Eq:a/GX/G},
$\cls(a,G)\in\Id\left(u^{(I)}/G\right)$. By Lemma~\ref{L:TruthIdFil},
$\bv{a\in\Id(u^{(I)})}\mm\ga_{(I)}\in G$, whence $g(a)\in G$. This completes
the proof of \eqref{Eq:gIdFct}.
\end{proof}

Note the following immediate corollary of Proposition~\ref{P:Compf*}:

\begin{corollary}\label{C:Compf*}
Let $P$ be a $D$-valued partial lattice, let
$f\colon P\to D$, let $\ga\in D$.
\begin{enumerate}
\item Suppose that $P$ has \IDJ.
If $f$ is an affine lower function, then $(f\mm\ga)^{\Id}=f^{\Id}\mm\ga$.

\item Suppose that $P$ has \FIJ.
If $f$ is an affine upper function, then $(f\mm\ga)^{\Fil}=f^{\Fil}\mm\ga$.
\end{enumerate}
\end{corollary}

\section{The lattices of affine ideal functions and affine filter functions}

In this section, we fix a distributive lattice $D$ with unit, and a
$D$-valued partial lattice $P$.

\begin{lemma}\label{L:fmmg,fgIdFil}\hfill
\begin{enumerate}
\item Suppose that $P$ has \IDM. Then the meet of any two affine lower
functions on~$P$ is an affine lower function on~$P$. In particular,
$\aID(P)$ is closed under meet.

\item Suppose that $P$ has \FIM. Then the meet of any two affine upper
functions on~$P$ is an affine upper function on~$P$. In particular,
$\aFIL(P)$ is closed under meet.
\end{enumerate}
\end{lemma}

\begin{proof}
We provide a proof for (i). Let $f$, $g\in\aID(P)$. Write
 \[
 f\colon x\mapsto\JJ_{i<m}\bv{x\leq u_i}\mm\ga_i,\qq
 g\colon x\mapsto\JJ_{j<n}\bv{x\leq v_j}\mm\gb_j.
 \]
By \IDM, there exists a common \IDM-sample $W$ for all pairs
$\set{u_i,v_j}$, for $\vv<i,j>\in m\times n$. This means that
 \[
 \bv{x\leq u_i}\mm\bv{x\leq v_j}=\JJ_{w\in W}\bv{x\leq w}\mm\gg_{i,j,w},
 \qq\text{for all }x\in P,
 \]
where we put $\gg_{i,j,w}=\bv{w\leq u_i}\mm\bv{w\leq v_j}$, for all
$\vv<i,j,w>\in m\times n\times W$. Hence, for $x\in P$,
 \begin{align*}
 f(x)\mm g(x)&=
 \JJ_{\substack{i<m\\ j<n}}\bv{x\leq u_i}\mm\bv{x\leq v_j}\mm\ga_i\mm\gb_j\\
 &=\JJ_{w\in W}\bv{x\leq w}\mm\gg_w,
 \end{align*}
where we put $\gg_w=\JJ_{\vv<i,j>\in m\times n}\gg_{i,j,w}\mm\ga_i\mm\gb_j$,
for all $w\in W$.
Therefore, $f\mm g$ is an affine lower function.

Since the meet of any two ideal functions on~$P$ is an ideal function on
$P$, it follows that $\aID(P)$ is closed under meet.
\end{proof}

\begin{corollary}\label{C:IFaffLatt}\hfill
\begin{enumerate}
\item If $P$ has both \IDM\ and \IDJ, then $\aID(P)$, with componentwise
ordering, is a lattice.

\item If $P$ has both \FIM\ and \FIJ, then $\aFIL(P)$, with componentwise
ordering, is a lattice.
\end{enumerate}
\end{corollary}

\begin{proof}
We provide a proof for (i).
By Lemma~\ref{L:fmmg,fgIdFil}, $\aID(P)$ is closed under meet. If $f$,
$g\in\aID(P)$, then $f\jj g$ (the componentwise join of~$f$ and $g$) is an
affine lower function on~$P$, thus, by Proposition~\ref{P:Compf*},
$(f\jj g)^{\Id}$ is defined, and it belongs to $\aID(P)$. So
$(f\jj g)^{\Id}$ is the join of~$\set{f,g}$ in $\aID(P)$.
\end{proof}

In order to differentiate between the componentwise join $f\jj g$
and the join of~$\set{f,g}$ in $\aID(P)$ (or $\aFIL(P)$), we introduce a
notation:

\begin{notation}
Under the assumptions of Corollary~\ref{C:IFaffLatt}, we denote by
$f\ijj g$ (resp., $f\fjj g$) the join of~$\set{f,g}$ in $\aID(P)$
(resp., in $\aFIL(P)$).
\end{notation}

Our next goal is to relate the meet and the join in $\aID(P)$ and $\aFIL(p)$
on the one hand, and the meet (intersection) and the join in $\ID(P/G)$ and
$\FIL(P/G)$ on the other hand, for a prime filter $G$ of~$D$.
For a lower function $f\colon P\to D$, the inverse image
$f^{-1}G$ of~$G$ has the property that if $y\in f^{-1}G$ and $x\leq_Gy$ (the
preordering $\leq_G$ has been introduced in Section~\ref{S:DvalPos}), then
$x\in f^{-1}G$. Hence, $\cls(x,G)\in f^{-1}G/G$ if{f} $f(x)\in G$. This also
holds for upper functions on~$P$. Our next result analyzes in more detail
the map $f\mapsto f^{-1}G/G$.

\begin{proposition}\label{P:HomIdFilP/G}
Let $G$ be a prime filter of~$D$.
\begin{enumerate}
\item Suppose that $P$ has \IDM\ and \IDJ. Then the rule $f\mapsto f^{-1}G/G$
determines a lattice homomorphism from $\vv<\aID(P),\mm,\ijj>$ to
$\vv<\ID(P/G),\ii,\jj>$.

\item Suppose that $P$ has \FIM\ and \FIJ. Then the rule $f\mapsto f^{-1}G/G$
determines a lattice homomorphism from $\vv<\aFIL(P),\mm,\fjj>$ to
$\vv<\FIL(P/G),\ii,\jj>$.
\end{enumerate}
\end{proposition}

\begin{proof}
We provide a proof for (i). We denote by $\gp_G$ the map
$f\mapsto f^{-1}G/G$. If $f\colon x\mapsto\bv{x\leq u}\mm\ga$, for fixed
$u\in P$ and $\ga\in D$, then $\gp_G(f)$ equals $\dnw\cls(u,G)$ if
$\ga\in G$, $\es$ otherwise, so $\gp_G(f)$ is in both cases an ideal of
$P/G$.

To prove that $\gp_G$ is a join-homomorphism with range contained in
$\Id(P/G)$, it suffices to prove that if
$f=\sideset{}{_{i<n}^{\Id}}{\JJ}f_i$, where $n\in\gos$
and $f_i\colon a\mapsto\bv{a\leq u_i}\mm\ga_i$ for all $i<n$ (where
$u_i\in P$ and $\ga_i\in D$), then $\gp_G(f)$ is the join of
$\setm{\gp_G(f_i)}{i<n}$ in $\Id(P/G)$, that is, we must prove that
 \begin{equation}\label{Eq:f-1Gfi-1G}
 f^{-1}G/G=\Id\left(\UU_{i<n}f_i^{-1}G/G\right).
 \end{equation}
So, let $a\in P$. Suppose first that $\cls(a,G)\in f^{-1}G/G$, that is,
$f(a)\in G$. By the formula given for $f$ in Proposition~\ref{P:Compf*}(i),
there exists a nonempty subset $I$ of~$n$ such that, using the same
notations as in Proposition~\ref{P:Compf*}(i), $\ga_{(I)}\in G$ and
$\bv{a\in\Id(u^{(I)})}\in G$. Therefore, $\ga_i\in G$ for all $i\in I$, and,
by Lemma~\ref{L:TruthIdFil}, $\cls(a,G)\in\Id(u^{(I)}/G)$. But, for $i\in I$,
$f_i(u_i)=\ga_i\in G$, thus $\cls(u_i,G)\in f_i^{-1}G/G$. Therefore,
$\cls(a,G)\in\Id\left(\UU_{i<n}f_i^{-1}G/G\right)$.

Conversely, suppose that $\cls(a,G)\in\Id\left(\UU_{i<n}f_i^{-1}G/G\right)$.
We observe that $\UU_{i<n}f_i^{-1}G/G$ is generated, as a lower subset of
$P/G$, by $u^{(I)}/G$, where $I=\setm{i<n}{\ga_i\in G}$. Thus,
$\cls(a,G)\in\Id(u^{(I)}/G)$, so, by Lemma~\ref{L:TruthIdFil},
$\bv{a\in\Id(u^{(I)})}\in G$. Since $\ga_{(I)}\in G$,
$\bv{a\in\Id(u^{(I)})}\mm\ga_{(I)}\in G$, whence $f(a)\in G$, that is,
$\cls(a,G)\in f^{-1}G/G$.

So we have proved that \eqref{Eq:f-1Gfi-1G} holds. As remarked above, this
shows that $\gp_G$ is a join-homomorphism with range a subset of~$\Id(P/G)$.

To conclude the proof, it is sufficient to prove that $\gp_G$ is a
meet-homomorphism. This is easy: for $a\in P$,
 \begin{align*}
 \cls(a,G)\in\gp_G(f\mm g)&\q\text{if{f}}\q f(a)\mm g(a)\in G\\
 &\q\text{if{f}}\q f(a)\in G\text{ and }g(a)\in G\\
 &\q\text{if{f}}\q\cls(a,G)\in\gp_G(f)\text{ and }\cls(a,G)\in\gp_G(g)\\
 &\q\text{if{f}}\q\cls(a,G)\in\gp_G(f)\ii\gp_G(g).
 \end{align*}
Therefore, $\gp_G(f\mm g)=\gp_G(f)\ii\gp_G(g)$.
\end{proof}

\section{The elements $\vbv{f\leq g}$}

We first introduce a convenient notation.

\begin{notation}
Let $\famm{\ga_i}{i\in I}$ be a family of elements of a lattice $D$, let
$\ga\in D$. Let $\ga=\fJJ{i\in I}\ga_i$ hold, if there exists a \emph{finite}
subset $J$ of~$I$ such that $\ga_i\leq\JJ_{j\in J}\ga_j$, for all $i\in I$,
and $\ga=\JJ_{j\in J}\ga_j$.
\end{notation}

Hence, $\ga=\fJJ{i\in I}\ga_i$ means that the supremum of the $\ga_i$ is,
really, the supremum of a \emph{finite} subfamily of~$\famm{\ga_i}{i\in I}$.

For the remainder of this section, let $D$ be a distributive lattice with
unit and let $P$ be a $D$-valued poset.

\begin{lemma}\label{L:Basic[fleqg]}
Let $m$, $n\in\gos$, let $u_0$,\dots, $u_{m-1}$, $v_0$, \dots,
$v_{n-1}\in P$, $\ga_0$,\dots, $\ga_{m-1}$, $\gb_0$,\dots,
$\gb_{n-1}\in D$. We define maps $f$ and $g$ from $P$ to $D$ by the rules
 \[
 f(x)=\JJ_{i<m}\bv{u_i\leq x}\mm\ga_i,\qq
 g(x)=\JJ_{j<n}\bv{x\leq v_j}\mm\gb_j,
 \]
for all $a\in P$. Put
 \[
 \gg=\JJ_{\vv<i,j>\in m\times n}\ga_i\mm\gb_j\mm\bv{u_i\leq v_j}.
 \]
Then $\gg=\fJJ{x\in P}f(x)\mm g(x)$.
\end{lemma}

\begin{proof}
For $x\in P$,
 \begin{align*}
 f(x)\mm g(x)&=
 \JJ_{\vv<i,j>\in m\times n}
 \bv{u_i\leq x}\mm\bv{x\leq v_j}\mm\ga_i\mm\gb_j\\
 &\leq\JJ_{\vv<i,j>\in m\times n}
 \bv{u_i\leq v_j}\mm\ga_i\mm\gb_j\\
 &=\gg.
 \end{align*}
Conversely, for $\vv<i,j>\in m\times n$,
 \[
 \ga_i\mm\gb_j\mm\bv{u_i\leq v_j}\leq\ga_i\mm g(v_j)\mm\bv{u_i\leq v_j}
 \leq f(v_j)\mm g(v_j).
 \]
The conclusion follows, with $\gg=\JJ_{j<n}f(v_j)\mm g(v_j)$.
\end{proof}

\begin{definition}\label{D:vbvfleqg}
For an affine upper function $f\colon P\to D$
and an affine lower function $g\colon P\to D$, we put
 \[
 \vbv{f\leq g}=\fJJ{x\in P}f(x)\mm g(x).
 \]
\end{definition}

By Lemma~\ref{L:Basic[fleqg]}, $\vbv{f\leq g}$ is always defined, and it is
an element of~$D$.

\begin{remark}\label{Rk:Expr[fleqg]}
With $f$ and $g$ defined as in the statement of Lemma~\ref{L:Basic[fleqg]},
we have obtained that
 \[
 \vbv{f\leq g}=\JJ_{j<n}f(v_j)\mm g(v_j).
 \]
We could have obtained, similarly, that
 \[
 \vbv{f\leq g}=\JJ_{i<m}f(u_i)\mm g(u_i).
 \]
These expressions will be used in Lemma~\ref{L:x+-(a)}.
\end{remark}

\section{Extension of the Boolean values to $\W(P)$}\label{S:BVtoW(P)}

Throughout this section, let $D$ be a distributive lattice with unit, let
$P$ be a $D$-valued partial lattice with \IDM, \FIM, \IDJ, and \FIJ. We
shall extend the notation $\bv{a\leq b}$, for $a$, $b\in P$, to all pairs of
elements of~$\W(P)$.

\begin{definition}\label{D:x+x-Dval}
For $\dx\in\W(P)$, we define, by induction on the height of~$\dx$, an affine
ideal function $\dx^-$ and an affine filter function $\dx^+$ on $P$ by the
following rules:
\begin{enumerate}
\item If $\dx=a\in P$, then $\dx^-\colon t\mapsto\bv{t\leq a}$
and $\dx^+\colon t\mapsto\bv{a\leq t}$.

\item $(\dx\mm\dy)^-=\dx^-\mm\dy^-$, and
$(\dx\jj\dy)^-=\dx^-\ijj\dy^-$, for all $\dx$, $\dy\in\W(P)$.

\item $(\dx\mm\dy)^+=\dx^+\fjj\dy^+$, and
$(\dx\jj\dy)^+=\dx^+\mm\dy^+$, for all $\dx$, $\dy\in\W(P)$.
\end{enumerate}

\end{definition}

We can now provide the $D$-valued analogue of the notation $\dx\ll\dy$
introduced in Definition~\ref{D:xlly}, by using the elements $\vbv{f\leq g}$,
see Definition~\ref{D:vbvfleqg}:

\begin{definition}
For $\dx$, $\dy\in\W(P)$, we put
 \[
 \bv{\dx\ll\dy}=\vbv{\dx^+\leq\dy^-}.
 \]
\end{definition}

As an easy consequence of Remark~\ref{Rk:Expr[fleqg]}, we
record the following:

\begin{lemma}\label{L:x+-(a)}
For $a\in P$ and $\dx\in\W(P)$, the following equalities hold:
 \[
 \bv{a\ll\dx}=\dx^-(a),\qq\bv{\dx\ll a}=\dx^+(a)
 \]
\end{lemma}

\begin{definition}\label{D:bvdxleqdy}
We define $\bv{\dx\leq\dy}$, for $\dx$, $\dy\in\W(P)$, by induction on
$\max\set{\hgt(\dx),\hgt(\dy)}$, as follows:
\begin{enumerate}
\item $\bv{\dx\leq\dy}=\bv{\dx\ll\dy}$, if $\dx\in P$ or $\dy\in P$,

\item $\bv{\dx_0\jj\dx_1\leq\dy_0\mm\dy_1}=\MM_{i,j<2}\bv{\dx_i\leq\dy_j}$,

\item
$\bv{\dx_0\jj\dx_1\leq\dy_0\jj\dy_1}=\MM_{i<2}\bv{\dx_i\leq\dy_0\jj\dy_1}$.

\item
$\bv{\dx_0\mm\dx_1\leq\dy_0\mm\dy_1}=\MM_{j<2}\bv{\dx_0\mm\dx_1\leq\dy_j}$.

\item
$\bv{\dx_0\mm\dx_1\leq\dy_0\jj\dy_1}=\bv{\dx_0\mm\dx_1\ll\dy_0\jj\dy_1}
\jj\JJ_{i,j<2}\bv{\dx_i\leq\dy_j}$.
\end{enumerate}
\end{definition}

\begin{proposition}\label{P:W(P)DvalLatt}
Let $\dx$, $\dx_0$, $\dx_1$, $\dy$, $\dy_0$, $\dy_1\in\W(P)$. Then the
following equalities hold:
 \begin{align}
 \bv{\dx_0\jj\dx_1\leq\dy}&=\bv{\dx_0\leq\dy}\mm\bv{\dx_1\leq\dy};
 \label{Eq:x0x1dy}\\
 \bv{\dx\leq\dy_0\mm\dy_1}&=\bv{\dx\leq\dy_0}\mm\bv{\dx\leq\dy_1}.
 \label{Eq:dxy0y1}
 \end{align}
\end{proposition}

\begin{proof}
By induction on $n<\go$, we prove that \eqref{Eq:x0x1dy} (resp.,
\eqref{Eq:dxy0y1}) holds for all $\dx_0$, $\dx_1$, $\dy$ such that
$\hgt(\dx_0)+\hgt(\dx_1)+\hgt(\dy)\leq n$ (resp., for all $\dx$, $\dy_0$,
$\dy_1$ such that $\hgt(\dx)+\hgt(\dy_0)+\hgt(\dy_1)\leq n$). Let us prove,
for example, that \eqref{Eq:x0x1dy} holds.

Suppose first that $\dy=a\in P$. We compute:
 \begin{align*}
 \bv{\dx_0\jj\dx_1\leq a}&=\bv{\dx_0\jj\dx_1\ll a}&&
 (\text{by the definition of }\bv{\dx\leq a})\\
 &=(\dx_0\jj\dx_1)^+(a)&&(\text{by Lemma~\ref{L:x+-(a)}})\\
 &=\dx_0^+(a)\mm\dx_1^+(a)&&\\
 &=\bv{\dx_0\ll a}\mm\bv{\dx_1\ll a}&&\\
 &=\bv{\dx_0\leq a}\mm\bv{\dx_1\leq a}.
 \end{align*}
Next, suppose that $\dy=\dy_0\mm\dy_1$. We compute:
 \begin{align*}
 \bv{\dx_0\jj\dx_1\leq\dy}&=\MM_{i,j<2}\bv{\dx_i\leq\dy_j}\\
\intertext{(by Definition~\ref{D:bvdxleqdy}(ii))}
 &=\bv{\dx_0\leq\dy_0\mm\dy_1}\mm\bv{\dx_1\leq\dy_0\mm\dy_1}
\intertext{(by the induction hypothesis about \eqref{Eq:dxy0y1})}
 &=\bv{\dx_0\leq\dy}\mm\bv{\dx_1\leq\dy}.
 \end{align*}
Finally, the case $\dy=\dy_0\jj\dy_1$ is trivial
(see Definition~\ref{D:bvdxleqdy}(iii)).
\end{proof}

For any prime filter $G$ on~$D$, we consider the canonical map
 \[
 \gr_G\colon\W(P)\twoheadrightarrow\W(P/G),\ \dx\mapsto\cls(\dx,G).
 \]
For $a\in P$ and for $\dx$, $\dy\in\W(P)$, $\gr_G$ satisfies
 \begin{align*}
 \gr_G(a)&=\cls(a,G)\in P/G,\\
 \gr_G(\dx\jj\dy)&=\gr_G(\dx)\jj\gr_G(\dy),\\
 \gr_G(\dx\mm\dy)&=\gr_G(\dx)\mm\gr_G(\dy).
 \end{align*}

We now relate the $D$-valued $\dx^-$ and $\dx^+$ with their corresponding
classical versions, see Definition~\ref{D:x-x+class}:

\begin{lemma}\label{L:x+-(a)inG}
Let $G$ be a prime filter of~$D$.
For any $a\in P$ and $\dx\in\W(P)$, the following equivalences hold:
 \begin{align}
 \dx^-(a)\in G&\q\text{if{f}}\q\cls(a,G)\in(\cls(\dx,G))^-,
 \label{Eq:aGdx-}\\
 \dx^+(a)\in G&\q\text{if{f}}\q\cls(a,G)\in(\cls(\dx,G))^+
 \label{Eq:aGdx+}
\end{align}
In other words, $(\dx^-)^{-1}G/G=(\cls(\dx,G))^-$ and
$(\dx^+)^{-1}G/G=(\cls(\dx,G))^+$.
\end{lemma}

\begin{proof}
We provide a proof for \eqref{Eq:aGdx-}. We argue by induction on the height
of~$\dx$. If $\dx=b\in P$, we must prove $\bv{a\leq b}\in G$ if{f}
$\cls(a,G)\in(\cls(b,G))^-$, which is the definition of the ordering in
$P/G$.

Suppose that $\dx=\dx_0\jj\dx_1$. Then $\dx^-=\dx_0^-\ijj\dx_1^-$. We
compute further:
 \begin{align*}
 (\dx^-)^{-1}G/G&=(\dx_0^-)^{-1}G/G\jj(\dx_1^-)^{-1}G/G&&
 \text{in }\ID(P/G)\\
\intertext{(by Proposition~\ref{P:HomIdFilP/G})}
 &=(\cls(\dx_0,G))^-\jj(\cls(\dx_1,G))^-\\
\intertext{(by the induction hypothesis)}
 &=(\cls(\dx_0,G)\jj\cls(\dx_1,G))^-\\
\intertext{(see Definition~\ref{D:x-x+class})}
 &=(\cls((\dx_0\jj\dx_1),G))^-\\
 &=(\cls(\dx,G))^-.
 \end{align*}
The proof for the case $\dx=\dx_0\mm\dx_1$ is similar.
\end{proof}

It is now easy to relate the symbols $\dx\ll\dy$
(see Definition~\ref{D:xlly}) and $\bv{\dx\ll\dy}$:

\begin{corollary}\label{C:xllyClassBool}
Let $G$ be a prime filter of~$D$,
let $\dx$, $\dy\in\W(P)$. Then $\cls(\dx,G)\ll\cls(\dy,G)$ 
(in $\W(P/G)$) if{f} $\bv{\dx\nobreak\ll\nobreak\dy}\nobreak\in\nobreak G$.
\end{corollary}

We are now ready to extend Corollary~\ref{C:xllyClassBool} to the case
$\bv{\dx\leq\dy}$:

\begin{proposition}\label{P:xleqyClassBool}
Let $G$ be a prime filter of~$D$,
let $\dx$, $\dy\in\W(P)$. Then $\cls(\dx,G)\peq\cls(\dy,G)$
(in $\W(P/G)$) if{f} $\bv{\dx\nobreak\leq\nobreak\dy}\nobreak\in\nobreak G$.
\end{proposition}

\begin{proof}
We argue by induction on the pair $\vv<\hgt(\dx),\hgt(\dy)>$. If $\dx\in P$
or $\dy\in P$, the conclusion follows from Corollary~\ref{C:xllyClassBool}.
The cases $\dx=\dx_0\jj\dx_1$ and $\dy=\dy_0\mm\dy_1$,
$\dx=\dx_0\jj\dx_1$ and $\dy=\dy_0\jj\dy_1$, and
$\dx=\dx_0\mm\dx_1$ and $\dy=\dy_0\mm\dy_1$ are obvious by the induction
hypothesis. The case $\dx=\dx_0\mm\dx_1$ and $\dy=\dy_0\jj\dy_1$ is easy by
Corollary~\ref{C:xllyClassBool} and the induction hypothesis.
\end{proof}

By using the fact that the relation $\peq$ (see Definition~\ref{D:peq}) is
reflexive and transitive in each $P/G$ (see Lemma~\ref{L:x-x+llpeq}) and by
Lemma~\ref{L:ManyG's}, we obtain the following consequence:

\begin{corollary}\label{C:BVleqTrans}
For $\dx$, $\dy$, $\dz\in\W(P)$, the following inequalities hold:
 \begin{align*}
 \bv{\dx\leq\dx}&=1;\\
 \bv{\dx\leq\dy}\mm\bv{\dy\leq\dz}&\leq\bv{\dx\leq\dz}.
 \end{align*}
\end{corollary}

\begin{remark}
Of course, the identity $\bv{\dx\leq\dx}=1$ is easy to prove directly.
However, proving the inequality
$\bv{\dx\leq\dy}\mm\bv{\dy\leq\dz}\leq\bv{\dx\leq\dz}$ directly is much less
intuitive (though, of course, possible) if one has to avoid the use of prime
filters of~$D$.
\end{remark}

In particular, $\W(P)$ is a $D$-valued poset.

\part{\CMPL{D}s}\label{Pt:DComeas}

\section{Finitely covered \CMPL{D}s}

Our next definition will be a combination between the definition of a
$D$-valued poset and a partial lattice.
In this section, we fix a distributive lattice $D$ with unit.

\begin{definition}\label{D:DvalMPL}
A \emph{\CMPL{D}} is a structure
\linebreak
$\vv<P,\bv{\ghost},\leq,\JJ,\MM>$ that satisfies the following axioms:
\begin{enumerate}
\item $\vv<P,\leq,\JJ,\MM>$ is a partial lattice.

\item $\vv<P,\bv{\ghost}>$ is a $D$-valued poset.

\item $x\leq y$ implies that $\bv{x\leq y}=1$.

\item $a=\JJ X$ implies that $\bv{a\leq b}=\MM_{x\in X}\bv{x\leq b}$,
for all $a$, $b\in P$ and all $X\in\dom\JJ$.

\item $a=\MM X$ implies that $\bv{b\leq a}=\MM_{x\in X}\bv{b\leq x}$,
for all $a$, $b\in P$ and all $X\in\dom\MM$.
\end{enumerate}

\end{definition}

\begin{definition}\label{D:D-cover}
Let $P$ be a $D$-valued poset, let $\E{D}$ be a subset of
$\fine P$. If $X\in\fine P$, a \emph{$\E{D}$-cover} of~$X$ is a
\emph{nonempty}, \emph{finite} subset $\E{D}'$ of~$\E{D}$ such that
 \[
 \bv{X=Y}\leq\JJ_{Z\in\E{D}'}\bv{X=Z},\q\text{for all }Y\in\E{D}.
 \]
We say that $\E{D}$ is \emph{finitely covering}, if every element
of~$\fine P$ has a $\E{D}$-cover.

We say that a \CMPL{D} $P$ is \emph{finitely covering}, if
both $\dom\JJ$ and $\dom\MM$ are finitely covering subsets of~$\fine P$.
\end{definition}

Observe that if $\E{D}'$ is a $\E{D}$-cover of~$X$, then any finite subset
of~$\E{D}$ that contains $\E{D}'$ is also a $\E{D}$-cover of~$X$.

Observe also that in the definition of a $\E{D}$-cover, we could have
replaced the inequality
 \[
 \bv{X=Y}\leq\JJ_{Z\in\E{D}'}\bv{X=Z}
 \]
by the inequality
 \[
 \bv{X=Y}\leq\JJ_{Z\in\E{D}'}\bv{Y=Z},
 \]
since $\bv{X=Y}\mm\bv{X=Z}=\bv{X=Y}\mm\bv{Y=Z}$ (by Lemma~\ref{L:XciYciZ}).

\begin{lemma}\label{L:CollDCov}
Let $P$ be a $D$-valued poset, let $\E{D}$ be a finitely covering subset of
$\fine P$. Then for all $X\in\fine P$, there exists a nonempty finite subset
$\E{D}'$ of~$\E{D}$ such that
 \[
 \bv{Y\ci X}\leq\JJ_{Z\in\E{D}'}\bv{Y=Z},\qq\text{for all }Y\in\E{D}.
 \]
\end{lemma}

\begin{proof}
Pick a common $\E{D}$-cover, $\E{D}'$, of all nonempty subsets of~$X$. We
compute:
 \begin{align*}
 \bv{Y\ci X}&=\JJ_{\es\sci T\ci X}\bv{Y=T}\\
\intertext{(by Lemma~\ref{L:XciYX=Z})}
 &=\JJ_{\es\sci T\ci X}\JJ_{Z\in\E{D}'}\bv{Y=T}\mm\bv{T=Z}\\
\intertext{(because $Y\in\E{D}$)}
 &\leq\JJ_{Z\in\E{D}'}\bv{Y=Z}
 \end{align*}
by Lemma~\ref{L:XciYciZ}.
\end{proof}

We observe that the condition that $P$ be finitely covering,
for a \CMPL{D} $P$, implies that both
$\dom\JJ$ and $\dom\MM$ are nonempty.

The following observation, although trivial, provides us with two important
classes of finitely covering \CMPL{D}s.

\begin{proposition}\label{P:LatFinDcov}
Let $P$ be a \CMPL{D}. Each of the following conditions
implies that $P$ is finitely covering:

\begin{enumerate}
\item $P$ is \emph{finite}, and both $\dom\JJ$ and $\dom\MM$ are nonempty.

\item $P$ is a \emph{lattice}, that is, $\dom\JJ=\dom\MM=\fine P$.
\end{enumerate}

\end{proposition}

\begin{proof}
(i) It is obvious that for every nonempty subset $\E{D}$ of~$\fine P$,
$\E{D}$ is a $\E{D}$-cover of every element of~$\fine P$. This holds, in
particular, for $\dom\JJ$ and $\dom\MM$.

(ii) For $X\in\fine P$, $\set{X}$ is simultaneously a $\dom\JJ$-cover and a
$\dom\MM$-cover of~$X$.
\end{proof}

Now we state the fundamental connection between \CMPL{D}s and $D$-valued
partial lattices:

\begin{proposition}\label{P:DmeasDval}
Let $P$ be a finitely covering \CMPL{D}. Then $P$ extends
to a $D$-valued partial lattice $\tP$ that satisfies
 \begin{align}
 \bv{a=\JJ X}&=\fJJe\famm{\bv{a=b}\mm\bv{X=Y}}
{\vv<b,Y>\in P\times\fine P,\ b=\JJ Y},\label{Eq:aJJXbJJY}\\
 \bv{a=\MM X}&=\fJJe\famm{\bv{a=b}\mm\bv{X=Y}}
{\vv<b,Y>\in P\times\fine P,\ b=\MM Y},\label{Eq:aMMXbMMY}
 \end{align}
for all $\vv<a,X>\in P\times\fine P$. Furthermore, $\tP$ is finitely sampled.
\end{proposition}

\begin{proof}
We first prove a claim.

\setcounter{claim}{0}
\begin{claim}
\hfill
\begin{itemize}
\item[(i)] $\bv{X\ci Y}\leq\bv{a\leq b}$, for all $a$, $b\in P$ and all $X$,
$Y\in\fine P$ such that $a=\JJ X$ and $b=\JJ Y$.

\item[(i*)] $\bv{X\ci Y}\leq\bv{b\leq a}$, for all $a$, $b\in P$ and all $X$,
$Y\in\fine P$ such that $a=\MM X$ and $b=\MM Y$.

\item[(ii)] $\bv{X=Y}\leq\bv{a=b}$, for all $a$, $b\in P$ and all $X$,
$Y\in\fine P$ such that $a=\JJ X$ and $b=\JJ Y$.

\item[(ii*)] $\bv{X=Y}\leq\bv{a=b}$, for all $a$, $b\in P$ and all $X$,
$Y\in\fine P$ such that $a=\MM X$ and $b=\MM Y$.
\end{itemize}

\end{claim}

\begin{cproof}
We first prove (i). Let $x\in X$. The inequality $y\leq b$ holds for all
$y\in Y$, thus $\bv{y\leq b}=1$, so
 \[
 \bv{x=y}=\bv{x=y}\mm\bv{y\leq b}\leq\bv{x\leq b}.
 \]
Hence,
 \[
 \bv{X\ci Y}\leq\bv{x\in Y}=\JJ_{y\in Y}\bv{x=y}\leq\bv{x\leq b},
 \]
so, since $a=\JJ X$,
 \[
 \bv{X\ci Y}\leq\MM_{x\in X}\bv{x\leq b}=\bv{a\leq b}.
 \]
(i*) is dual of (i), and (ii), (ii*) follow immediately.
\end{cproof}

We prove now that the equations \eqref{Eq:aJJXbJJY} and
\eqref{Eq:aMMXbMMY} are consistent definitions of $\bv{a=\JJ X}$ and
$\bv{a=\MM X}$. We do it, for example, for \eqref{Eq:aJJXbJJY}. So let
$\vv<a,X>\in P\times\fine P$. We put
 \[
 \ga=\JJ_{i<n}\bv{a=a_i}\mm\bv{X=X_i},
 \]
where $\setm{X_i}{i<n}$ is a $\dom\JJ$-cover of~$X$ and $a_i=\JJ X_i$, for
all $i<n$.
For $\vv<b,Y>\in P\times\fine P$ such that $b=\JJ Y$, we compute:
 \begin{align*}
 \bv{a=b}\mm\bv{X=Y}&=\JJ_{i<n}\bv{a=b}\mm\bv{X=Y}\mm\bv{X=X_i}\\
\intertext{(by the definition of a $\dom\JJ$-cover)}
 &=\JJ_{i<n}\bv{a=b}\mm\bv{X_i=Y}\mm\bv{X=X_i}\\
\intertext{(by an easy application of Lemma~\ref{L:XciYciZ})}
 &\leq\JJ_{i<n}\bv{a=b}\mm\bv{a_i=b}\mm\bv{X=X_i}\\
\intertext{(by Claim~1 applied to $\vv<X_i,Y>$)}
 &\leq\JJ_{i<n}\bv{a=a_i}\mm\bv{X=X_i}\\
 &=\ga,
 \end{align*}
hence, $\ga=\fJJe\famm{\bv{a=b}\mm\bv{X=Y}}
{\vv<b,Y>\in P\times\fine P,\ b=\JJ Y}$. This settles
\eqref{Eq:aJJXbJJY}. The proof for \eqref{Eq:aMMXbMMY} is dual.

We now verify that all items of Definition~\ref{D:DvalPL} are satisfied by
the Boolean values obtained above.
\smallskip

\noindent\textit{\textbf{Condition~1.}}
$\bv{a=\JJ X}\mm\bv{a\leq b}=\bv{a=\JJ X}\mm\MM_{x\in X}\bv{x\leq b}$.

Let $\setm{X_i}{i<n}$ be a $\dom\JJ$-cover of~$X$. Put $a_i=\JJ X_i$, for
all $i<n$. We compute:
 \begin{align*}
 \bv{a=\JJ X}\mm\bv{a\leq b}&=
 \JJ_{i<n}\bv{a=a_i}\mm\bv{X=X_i}\mm\bv{a\leq b}\\
 &=\JJ_{i<n}\bv{a=a_i}\mm\bv{X=X_i}\mm\bv{a_i\leq b}\\
\intertext{(because $\bv{a=a_i}\mm\bv{a\leq b}=\bv{a=a_i}\mm\bv{a_i\leq b}$)}
 &=\JJ_{i<n}\left(
 \bv{a=a_i}\mm\bv{X=X_i}\mm\MM_{x\in X_i}\bv{x\leq b}\right)\\
\intertext{(because $a_i=\JJ X_i$)}
 &=\JJ_{i<n}\left(
 \bv{a=a_i}\mm\bv{X=X_i}\mm\MM_{x\in X}\bv{x\leq b}\right)\\
\intertext{(by Lemma~\ref{L:XciYya})}
 &=\bv{a=\JJ X}\mm\MM_{x\in X}\bv{x\leq b}.
 \end{align*}

\noindent\textit{\textbf{Condition~2.}}
$\bv{a=\JJ X}\mm\bv{X=Y}\leq\bv{a=\JJ Y}$.

Again, let $\setm{X_i}{i<n}$ be a common $\dom\JJ$-cover of~$X$ and ~$Y$, and
put $a_i=\JJ X_i$, for all $i<n$. We compute:
 \begin{align*}
 \bv{a=\JJ X}\mm\bv{X=Y}
 &=\JJ_{i<n}\bv{a=a_i}\mm\bv{X=X_i}\mm\bv{X=Y}\\
 &\leq\JJ_{i<n}\bv{a=a_i}\mm\bv{Y=X_i}\\
\intertext{(by Lemma~\ref{L:XciYciZ})}
 &\leq\bv{a=\JJ Y}.
 \end{align*}

\noindent\textit{\textbf{Condition~3.}}
$\bv{a=\JJ X}\mm\bv{a=b}\leq\bv{b=\JJ X}$.

Again, let $\setm{X_i}{i<n}$ be a $\dom\JJ$-cover of~$X$, and put
$a_i=\JJ X_i$, for all $i<n$. We compute:
 \begin{align*}
 \bv{a=\JJ X}\mm\bv{a=b}
 &=\JJ_{i<n}\bv{a=a_i}\mm\bv{X=X_i}\mm\bv{a=b}\\
 &\leq\JJ_{i<n}\bv{b=a_i}\mm\bv{X=X_i}\\
 &=\bv{b=\JJ X}.
 \end{align*}

Hence we have verified items (1), (2), and (3) of Definition~\ref{D:DvalPL}.
The items (1*), (2*), and (3*) are dual.

At this point, we have verified that \eqref{Eq:aJJXbJJY} and
\eqref{Eq:aMMXbMMY} define a structure of $D$-valued partial lattice $\tP$
on~$P$.

It remains to prove that $\tP$ is finitely sampled. We verify, for
example, that $\tP$ is finitely join-sampled. So, let $X\in\fine P$. Let
$\setm{X_i}{i<n}$ be a $\dom\JJ$-cover of~$X$, and put
$a_i=\JJ X_i$, for all $i<n$. Put $U=\setm{a_i}{i<n}$.
For $a\in P$ and $i<n$,
 \[
 \bv{a=a_i}\mm\bv{X=X_i}\leq\bv{a_i=a_i}\mm\bv{X=X_i}\leq\bv{a_i=\JJ X},
 \]
so we obtain the following inequalities:
 \begin{align*}
 \bv{a=\JJ X}&=\JJ_{i<n}\bv{a=a_i}\mm\bv{X=X_i}\\
 &\leq\JJ_{i<n}\bv{a_i=\JJ X}\\
 &=\JJ_{u\in U}\bv{u=\JJ X}.
 \end{align*}
Therefore, $U$ is a join-sample of~$X$. Dually, $\tP$ is
finitely meet-sampled.
\end{proof}

\begin{remark}
In particular, in the context of Proposition~\ref{P:DmeasDval}, the
notation $\bv{a=\JJ X}$, for $\JJ X$ defined, is unambiguous, because
the singleton $\set{X}$ is a $\dom\JJ$-cover of~$X$, thus, if $b$ is defined
as $\JJ X$, then $\bv{a=\JJ X}$ (as defined in \eqref{Eq:aJJXbJJY}) equals
$\bv{a=b}$. Of course, the dual statement holds for the meet.
\end{remark}

We can now state a useful strengthening of Proposition~\ref{P:LatFinDcov}:

\begin{proposition}\label{P:LatFinDS}
Let $P$ be a \CMPL{D}. Each of the following conditions
implies that $P$ is finitely covering and that the associated $D$-valued
partial lattice satisfies \IDM, \IDJ, \FIM, and \FIJ:
\begin{enumerate}
\item $P$ is finite and both $\dom\JJ$ and $\dom\MM$ are nonempty.

\item $P$ is a lattice.
\end{enumerate}

\end{proposition}

\begin{proof}
In both cases, it follows from Proposition~\ref{P:LatFinDcov} that $P$ is
finitely covering. We denote by $\tP$ the associated finitely sampled
\CMPL{D}, see Proposition~\ref{P:DmeasDval}.

(i) We assume that $P$ is finite. For $X\in\fine P$, it is obvious that
$P$ is a finite \IDM-sample of~$X$. Hence, $\tP$ has \IDM.

We prove that $P$ is also an \IDJ-sample
of~$X$. Indeed, if $G$ is a prime filter of~$D$
and if $k<\go$, then the condition
 \[
 \Id_0(X/G,P/G)\sci\Id_1(X/G,P/G)\sci\cdots\sci\Id_k(X/G,P/G)
 \]
implies that $k<|P/G|$, thus, \emph{a fortiori}, $k<|P|$. In particular,
$\Id_{|P|-1}(X/G,P/G)=\Id_{|P|}(X/G,P/G)$. Hence, by
Lemma~\ref{L:EquivIdFilSG}, $P$ is an \IDJ-sample of~$X$, with index $|P|-1$.
So $\tP$ has \IDJ. The dual statements, about \FIM\ and \FIJ, are proved
similarly.

(ii) We assume that $P$ is a lattice. For $X\in\fine P$, if we put
$a=\MM X$, then $\bv{a=\MM X}=1$, thus $\Set{a}$ is an \IDM-sample
of~$X$. Put $b=\JJ X$, and let $U$ be a finite subset of $P$ containing $X$.
For every $x\in P$ and every $Z$ such that $\es\sci Z\ci U$,
$\bv{Z\ci\Id_0(X,U)}\leq\MM_{z\in Z}\bv{z\leq b}=\bv{\JJ Z\leq b}$, from
which it follows that
$\bv{x\leq\JJ Z}\mm\bv{Z\ci\Id_0(X,U)}\leq\bv{x\leq b}$. Since the value
$\bv{x\leq b}$ is reached for $Z=X$, we conclude that
$\bv{x\in\Id_1(X,U)}=\bv{x\leq b}$.
Now we observe that $x\mapsto\bv{x\leq b}$ is an ideal function on $\tP$.
It follows that $X$ is an \IDJ-sample of~$X$, with index $1$.
There are similar, dual statements, for \FIM\ and \FIJ.
\end{proof}

\section{Statement and proof of Theorem~A}\label{S:ThmA}

In order to relate \CMPL{D}s and congruence lattices of partial lattices,
we state the following simple result.

\begin{proposition}\label{P:DvalMPL}
Let $D$ be a distributive lattice with unit,
let $\vv<P,\leq,\JJ,\MM>$ be a partial lattice, endowed with a map
$P\times P\to D$, $\vv<x,y>\mapsto\bv{x\leq y}$. Then the following are
equivalent:

\begin{enumerate}
\item $\vv<P,\bv{\ghost},\leq,\JJ,\MM>$ is a \CMPL{D}.

\item There exists a homomorphism
$\gf\colon\vv<\Conc P,\jj,\zero_P>\to\vv<D,\mm,1>$ such that
 \[
 \bv{x\leq y}=\gf(\gQ_P^+(x,y)),\qq\text{for all }x,\,y\in P.
 \]
\end{enumerate}

\end{proposition}

\begin{proof}
(i)$\Rightarrow$(ii) It is sufficient to prove that if $n<\go$ and $a$, $b$,
$a_0$,\dots, $a_{n-1}$, $b_0$,\dots, $b_{n-1}\in P$, the condition
 \begin{equation}\label{Eq:gQ+abaibi}
 \gQ_P^+(a,b)\ci\JJ_{i<n}\gQ_P^+(a_i,b_i)
 \end{equation}
implies that
 \begin{equation}\label{Eq:bvabaibi}
 \bv{a\leq b}\geq\MM_{i<n}\bv{a_i\leq b_i}.
 \end{equation}
So put $\ga=\MM_{i<n}\bv{a_i\leq b_i}$, and define a binary relation
$\leq_\gq$ on~$P$ by the rule
 \[
 x\leq_\gq y\q\text{if{f}}\q\ga\leq\bv{x\leq y},\q\text{for all }x,\,y\in P.
 \]
We verify that $\gq$ is a congruence of the partial lattice $P$. It is
obvious that $\gq$ is a preordering of~$P$, see
Definition~\ref{D:DvalPoset}, and that $\gq$ contains the ordering $\leq$
of~$P$. Let $u\in P$, let $X\in\fine P$, we
verify that $u=\JJ X$ (resp., $u=\MM X$) implies that $u$ is the supremum
(resp., the infimum) of~$X$ with respect to $\gq$. We do it for example for
the join. {F}rom $X\leq u$, it follows that $X\leq_\gq u$. Now let $v\in P$
such that $X\leq_\gq v$, that is, $\ga\leq\bv{x\leq v}$, for all $x\in X$.
Since $u=\JJ X$, it follows from Definition~\ref{D:DvalMPL} that
$\ga\leq\bv{u\leq v}$, that is, $u\leq_\gq v$. This proves our assertion
about the supremum. The proof for the infimum is dual.

Now, by the definition of~$\ga$, the inequality $a_i\leq_\gq b_i$ holds for
all $i<n$, that is, $\gQ_P^+(a_i,b_i)\ci\gq$. Hence, by \eqref{Eq:gQ+abaibi},
$\gQ_P^+(a,b)\ci\gq$, that is, $\ga\leq\bv{a\leq b}$, in other words,
\eqref{Eq:bvabaibi} holds.

(ii)$\Rightarrow$(i) This direction of the proof follows immediately from the
identities
 \begin{align*}
 \gQ_P^+(a,b)&=\JJ_{x\in X}\gQ_P^+(x,b),&&\text{if }a=\JJ X,\\
 \gQ_P^+(b,a)&=\JJ_{x\in X}\gQ_P^+(b,x),&&\text{if }a=\MM X,
 \end{align*}
for all $a$, $b\in P$ and all $X\in\fine P$.
\end{proof}

\begin{definition}\label{D:D0PL}
Let $D$ be a distributive lattice with zero.
A \emph{\MPL{D}} is a pair $\vv<P,\gf>$, where $P$ is a partial lattice
and $\gf\colon\Conc P\to D$ is a \jzh.

If, in addition, $P$ is a lattice, we say that $\vv<P,\gf>$ is a
\emph{\ML{D}}.

We say that $\vv<P,\gf>$ is \emph{proper}, if $\gf$ isolates $0$
(that is, $\gf(\gq)=0$ if{f} $\gq=\zero_P$, for all $\gq$ in $\Conc P$).
\end{definition}

Hence, by Proposition~\ref{P:DvalMPL}, the notions of a \CMPL{D} and a
\MPL{D} are, up to dualization of~$D$, \emph{equivalent}.

\begin{definition}\label{D:SFSampled}
Let $D$ be a distributive lattice with unit, let $P$ be a
\emph{finitely covering} \CMPL{D}. We say that $P$ has \IDM\ (resp., \IDJ,
\FIM, \FIJ), if the $D$-valued partial lattice $\tP$ of
Proposition~\ref{P:DmeasDval} has \IDM\ (resp., \IDJ, \FIM, \FIJ).

We say that $P$ is \emph{balanced}, if it has \IDM, \IDJ, \FIM, and \FIJ.

If $D$ is a distributive lattice with zero, we say that a \MPL{D} $P$ is
balanced, if the associated \CMPL{D^\rd} is balanced. Similar definitions
hold for \IDM, \IDJ, \FIM, and \FIJ.
\end{definition}

By Proposition~\ref{P:LatFinDS}, every \ML{D}, or every finite
\MPL{D} with nonempty meet and join, is balanced.

Now we can provide a precise statement, and proof, for Theorem~A.
We recall that $j_P$ is the natural embedding from $P$ into $\FL(P)$, see
Section~\ref{S:DescFL(P)}.

\begin{all}{Theorem~A}
Let $D$ be a distributive lattice with zero, let $\vv<P,\gf>$ be a
balanced \MPL{D}.
Then there exists a \jzh\linebreak
$\gy\colon\Conc\FL(P)\to D$ such that
 \[
 \gy\circ\Conc j_P=\gf.
 \]
\end{all}

The remainder of this section will be devoted to the proof of Theorem~A.

We first endow $P$ with its natural structure of \CMPL{D^\rd}, see
Proposition~\ref{P:DvalMPL}. By assumption, this structure is balanced, that
is, it is finitely covering and it has \IDM, \IDJ, \FIM, and \FIJ, see
Definition~\ref{D:SFSampled}. We define elements
$\bv{\dx\leq\dy}$, for all elements $\dx$, $\dy$ of $\W(P)$, as in
Definition~\ref{D:bvdxleqdy}. We define binary relations $\peq^*$ and
$\equiv^*$ on $\W(P)$ by the rules
 \begin{align*}
 \dx\peq^*\dy&\qq\text{if{f}}\qq\bv{\dx\leq\dy}=1,\\
 \dx\equiv^*\dy&\qq\text{if{f}}\qq\bv{\dx=\dy}=1.
 \end{align*}
It follows from Corollary~\ref{C:BVleqTrans} that $\peq^*$ is a preordering
of $\W(P)$ and that $\equiv^*$ is the associated equivalence relation. Let
$L=\vv<\W(P),\peq^*>/{\equiv^*}$ be the quotient poset. For $\dx\in\W(P)$,
we denote by $[\dx]$ the $\equiv^*$-equivalence class of~$\dx$. For $x$,
$y\in L$ and $\dx\in x$, $\dy\in y$, the element $\bv{\dx\leq\dy}$ does not
depend of the choice of~$\vv<\dx,\dy>$, we denote it by $\bv{x\leq y}$.
Similarly, we define $\bv{x=y}=\bv{\dx=\dy}$.
Furthermore, it follows from Proposition~\ref{P:W(P)DvalLatt} that
$[\dx\jj\dy]$ (resp., $[\dx\mm\dy]$) is the supremum (resp., the infimum) of
$\set{x,y}$ in $L$.

Hence, \emph{$L$ is a lattice}. Furthermore, by
Proposition~\ref{P:W(P)DvalLatt}, the equality
 \[
 \bv{x_0\jj x_1\leq y}=\bv{x_0\leq y}\mm\bv{x_1\leq y}
 \]
holds for all $x_0$, $x_1$, $y\in L$. Symmetrically, the equality
 \[
 \bv{x\leq y_0\mm y_1}=\bv{x\leq y_0}\mm\bv{x\leq y_1}
 \]
holds for all $x$, $y_0$, $y_1\in L$.
Therefore, an easy induction shows that
$\bv{\JJ X\leq a}=\MM_{x\in X}\bv{x\leq a}$ and
$\bv{a\leq\MM X}=\MM_{x\in X}\bv{a\leq x}$, for all $a\in L$ and all
$X\in\fine L$. Hence, \emph{$L$ is a \CMPL{D^\rd}}, see
Definition~\ref{D:DvalMPL}. Therefore, by Proposition~\ref{P:DvalMPL}, there
exists a \jzh\ $\rho\colon\Conc L\to D$ such that
 \begin{equation}\label{Eq:rhogQ}
 \rho(\gQ_L^+(x,y))=\bv{x\leq y},\qq\text{for all }x,\,y\in L.
 \end{equation}
Furthermore, it is easy to verify that the rule $a\mapsto[a]$
defines a homomorphism of partial lattices from $P$ to $L$.
Thus, since $L$ is a lattice, there exists, by
Proposition~\ref{P:DescFL}, a unique lattice homomorphism
$f\colon\FL(P)\to L$ such that $f(a)=[a]$, for all $a\in P$.

We put $\gy=\rho\circ\Conc f$, a \jzh\ from $\Conc\FL(P)$ to $D$.
For $a$, $b\in P$,
 \[
 \gy\left(\gQ_{\FL(P)}^+(a,b)\right)=
 \rho\left(\gQ_L^+([a],[b])\right)=\bv{[a]\leq[b]}=\bv{a\leq b}=
 \gf\left(\gQ_P^+(a,b)\right),
 \]
so $\gy\circ\Conc j_P=\gf$. This concludes the proof of Theorem~A.

\section{Quotients of \CMPL{D}s by prime filters}

In this section, we fix a distributive lattice $D$ with unit.

If $P$ is a finitely covering \CMPL{D}, then, by
Proposition~\ref{P:DmeasDval}, $P$ extends canonically to a $D$-valued
partial lattice. So, by using Proposition~\ref{P:P/GPartLatt}, we can define
a partial lattice $P/G$, for every prime filter $G$ of~$D$. Our next result
describes the join and meet operations in $P/G$.

\begin{lemma}\label{L:PDmeasP/G}
Let $P$ be a finitely covering \CMPL{D}, let
$\xa\in P/G$, let $\xX\in\fine{P/G}$. Then the following assertions hold:
\begin{enumerate}
\item $\xa=\JJ\xX$ in $P/G$ if{f} there are $a\in P$ and $X\in\fine P$ such
that $\xa=\cls(a,G)$, $\xX=X/G$, and $a=\JJ X$.

\item $\xa=\MM\xX$ in $P/G$ if{f} there are $a\in P$ and $X\in\fine P$ such
that $\xa=\cls(a,G)$, $\xX=X/G$, and $a=\MM X$.
\end{enumerate}

\end{lemma}

\begin{proof}
We prove (i); (ii) is dual. Suppose first that $\xa=\JJ\xX$. Pick $a\in\xa$
and $X\in\fine P$ such that $\xX=X/G$. By the definition of the join
operation in $P/G$, $\bv{a=\JJ X}\in G$. Let $\setm{X_i}{i<n}$ (where $n>0$)
be a $\dom\JJ$-cover of~$X$. Put $a_i=\JJ X_i$, for all $i<n$.
So the equality
 \[
 \bv{a=\JJ X}=\JJ_{i<n}\bv{a=a_i}\mm\bv{X=X_i}
 \]
holds by definition, thus, since $G$ is prime, there exists $i<n$ such that
$\xa=\cls(a_i,G)$ and $\xX=X_i/G$. Since $a_i=\JJ X_i$, we have proved the
``if'' direction of the implication in (i).

To prove the converse, assume that $\xa=\cls(a,G)$, $\xX=X/G$, and
$a=\JJ X$. Since $\set{X}$ is a $\dom\JJ$-cover of~$X$,
$\bv{a=\JJ X}=1\in G$, hence $\cls(a,G)=\JJ(X/G)$ in $P/G$ by the definition
of~$\JJ$ in $P/G$, so $\xa=\JJ\xX$.
\end{proof}

\begin{definition}\label{D:IsomDmeas}
Let $P$ and $Q$ be \CMPL{D}s, let $f\colon P\to Q$ be a
homomorphism of partial lattices. We say that $f$ is
\begin{enumerate}
\item a \emph{uniform map}, if $\bv{x\leq y}\leq\bv{f(x)\leq f(y)}$, for all
$x$, $y\in P$.

\item an \emph{isometry}, if $f$ is an embedding of partial
lattices and $\bv{x\leq y}=\bv{f(x)\leq f(y)}$, for all $x$, $y\in P$.
\end{enumerate}
\end{definition}

\begin{lemma}\label{L:IsomGotoG}
Let $P$ and $Q$ be finitely covering \CMPL{D}s, let
$f\colon P\to Q$ be a uniform map. For any prime filter $G$ of~$D$, one
can define a homomorphism $f^G\colon P/G\to Q/G$
of partial lattices by the rule
 \begin{equation}\label{Eq:fG(x/G)f(x)/G}
 f^G(\cls(x,G))=\cls(f(x),G),\qq\text{for all }x\in P.
 \end{equation}
Furthermore, if $f$ is an isometry, then $f^G$ is an embedding of partial
lattices.
\end{lemma}

\begin{proof}
For $x$, $y\in P$, $\bv{x\leq y}\in G$ implies that $\bv{f(x)\leq f(y)}\in G$,
so, \eqref{Eq:fG(x/G)f(x)/G} defines a unique order-preserving map
$f^G\colon P/G\to Q/G$.

We prove that $f^G$ is a homomorphism of partial lattices. We do
it for example for the join. So let $\xa\in P/G$, $\xX\in\fine{P/G}$ such
that $\xa=\JJ\xX$. By Lemma~\ref{L:PDmeasP/G}, there are $a\in P$ and
$X\in\fine P$ such that $\xa=\cls(a,G)$, $\xX=X/G$, and $a=\JJ X$. Since $f$
is a homomorphism of partial lattices, $f(a)=\JJ f[X]$. Thus, again by
Lemma~\ref{L:PDmeasP/G}, $\cls(f(a),G)=\JJ f[X]/G$, that is,
$f^G(\xa)=\JJ f^G[\xX]$. The proof for the meet is dual.

Finally, if $f$ is an isometry, then
$\bv{x\leq y}\in G$ if{f} $\bv{f(x)\leq f(y)}\in G$, for all $x$, $y\in P$,
thus $f^G$ is an order-embedding.
\end{proof}

\section{Amalgamation of \CMPL{D}s above a finite lattice}

We extend in this section the results of Section~\ref{S:AmalgPart} to
\CMPL{D}s.

We fix a distributive lattice $D$ with unit.

\begin{definition}\label{D:DV-Form}
A \emph{V-formation} of \CMPL{D}s is a structure
$\vv<K,P,Q,f,g>$ subject to the following conditions:
\begin{itemize}
\item[(DV1)] $K$, $P$, $Q$ are \CMPL{D}s.

\item[(DV2)] $f\colon K\hookrightarrow P$ and $g\colon K\hookrightarrow Q$ are
\emph{isometries}.
\end{itemize}

A V-formation $\vv<K,P,Q,f,g>$ is \emph{standard}, if the following
conditions hold:
\begin{itemize}
\item[(SDV1)] $K$ is a \emph{finite lattice}.

\item[(SDV2)] $K=P\ii Q$ (set-theoretically), and $f$ and $g$ are,
respectively, the inclusion map from $K$ into $P$ and the inclusion map from
$K$ into $Q$.
\end{itemize}

\end{definition}

\begin{remark}
The definition of a standard V-formation of \CMPL{D}s is
\emph{not} a generalization of the definition of a standard V-formation of
partial lattices (Definition~\ref{D:V-Form}): indeed, observe the additional
requirement that $K$ be \emph{finite}.
\end{remark}

As in Section~\ref{S:AmalgPart}, we shall write $\vv<K,P,Q>$ instead of
$\vv<K,P,Q,f,g>$ for standard V-formations.

The following analogue of Lemma~\ref{L:V-FormSt} trivially holds:

\begin{lemma}\label{L:DV-FormSt}
Every V-formation $\vv<K,P,Q,f,g>$ of \CMPL{D}s, with $K$
a finite lattice, is isomorphic to a standard V-formation.
\end{lemma}

The definition of an amalgam or a pushout of a V-formation is, \emph{mutatis
mutandis}, exactly the same as in Definition~\ref{D:AmalgVform}. The
corresponding analogue of Proposition~\ref{P:PushStV} is then the following:

\begin{proposition}\label{P:DPushStV}
Let $\E{D}=\vv<K,P,Q,f,g>$ be a V-formation of \CMPL{D}s,
with $K$ a finite lattice. Then $\E{D}$ has a pushout in the category of
\CMPL{D}s and uniform maps.

Furthermore, assume that $\E{D}$ is a standard V-formation. Then the pushout
$\vv<R,f',g'>$ of~$\E{D}$ can be described by the following data:
\begin{itemize}
\item[(a)] $R=P\amalg_KQ$ as a partial lattice (see
Notation~\tup{\ref{Not:PamalgQ}}). In particular, $R=P\uu Q$
set-theoretically.

\item[(b)] For $x$, $y\in R$, the Boolean value $\bv{x\leq y}$ can be
computed as follows:

\begin{itemize}
\item[(b1)] $\bv{x\leq y}=\bv{x\leq y}_P$, if $x$, $y\in P$.

\item[(b2)] $\bv{x\leq y}=\bv{x\leq y}_Q$, if $x$, $y\in Q$.

\item[(b3)] $\bv{x\leq y}=\JJ_{z\in K}\bv{x\leq z}_P\mm\bv{z\leq y}_Q$, if
$x\in P$ and $y\in Q$.

\item[(b4)] $\bv{x\leq y}=\JJ_{z\in K}\bv{x\leq z}_Q\mm\bv{z\leq y}_P$, if
$x\in Q$ and $y\in P$.
\end{itemize}
\end{itemize}
Furthermore, both $f'$ and $g'$ are isometries.
\end{proposition}

\begin{proof}
We first prove the mutual compatibility of (b1)--(b4) above. Up to symmetry
between $P$ and $Q$, this amounts to considering the three following cases:
\smallskip

\noindent\textit{\textbf{Case~1.}} $x$, $y\in K$, prove that
$\bv{x\leq y}_P=\bv{x\leq y}_Q$. This follows immediately from the fact that
both maps $f\colon K\hookrightarrow P$ and $g\colon K\hookrightarrow Q$ are
isometries, thus $\bv{x\leq y}_P=\bv{x\leq y}_Q=\bv{x\leq y}_K$.

\smallskip
\noindent\textit{\textbf{Case~2.}} $x\in P$, $y\in K$, prove that
 \[
 \bv{x\leq y}_P=\JJ_{z\in K}\bv{x\leq z}_P\mm\bv{z\leq y}_Q.
 \]
For $z\in K$,
 \begin{align*}
 \bv{x\leq z}_P\mm\bv{z\leq y}_Q&=\bv{x\leq z}_P\mm\bv{z\leq y}_P\\
\intertext{(because $z$, $y\in K$)}
 &\leq\bv{x\leq y}_P,
 \end{align*}
and, for $z=y$, $\bv{x\leq z}_P\mm\bv{z\leq y}_Q=\bv{x\leq y}_P$, which
proves our assertion.

\smallskip
\noindent\textit{\textbf{Case~3.}} $x$, $y\in K$, prove that
 \[
 \JJ_{z\in K}\bv{x\leq z}_P\mm\bv{z\leq y}_Q
 =\JJ_{z\in K}\bv{x\leq z}_Q\mm\bv{z\leq y}_P.
 \]
By Case~2, the left hand side and the right hand side of the equality above
are both equal to $\bv{x\leq y}_P$ (and to $\bv{x\leq y}_Q$).

Now we verify that $\vv<x,y>\mapsto\bv{x\leq y}$ defines a structure of
$D$-valued poset on $R$. It is obvious that $\bv{a\leq a}=1$, for all
$a\in R$.

Now let $a$, $b$, $c\in R$, we prove the inequality
 \begin{equation}\label{Eq:bvaleblec}
 \bv{a\leq b}\mm\bv{b\leq c}\leq\bv{a\leq c}.
 \end{equation}
Up to symmetry between $P$ and $Q$, it is sufficient to consider the
three following cases:

\smallskip
\noindent\textit{\textbf{Case~1$'$.}} $a$, $b$, $c\in P$. Then
\eqref{Eq:bvaleblec} follows from the fact that $P$ is a $D$-valued poset.

\smallskip
\noindent\textit{\textbf{Case~2$'$.}} $a$, $b\in P$, and $c\in Q$.
We compute:
 \begin{align*}
 \bv{a\leq b}\mm\bv{b\leq c}&=
 \JJ_{x\in K}\bv{a\leq b}_P\mm\bv{b\leq x}_P\mm\bv{x\leq c}_Q\\
 &\leq\JJ_{x\in K}\bv{a\leq x}_P\mm\bv{x\leq c}_Q\\
 &=\bv{a\leq c}.
 \end{align*}

\smallskip
\noindent\textit{\textbf{Case~3$'$.}} $a\in P$, $b$, $c\in Q$. This case is
similar to Case~$2'$.

\smallskip
\noindent\textit{\textbf{Case~4$'$.}} $a\in P$, $b\in Q$, $c\in P$.
For $u$, $v\in K$,
 \begin{align*}
 \bv{a\leq u}_P\mm\bv{u\leq b}_Q\mm\bv{b\leq v}_Q\mm\bv{v\leq c}_P
 &\leq\bv{a\leq u}_P\mm\bv{u\leq v}_Q\mm\bv{v\leq c}_P\\
 &=\bv{a\leq u}_P\mm\bv{u\leq v}_P\mm\bv{v\leq c}_P\\
\intertext{(because $u$, $v\in K$)}
 &\leq\bv{a\leq c}_P\\
 &=\bv{a\leq c}.
 \end{align*}
It follows that
 \begin{align*}
 \bv{a\leq b}\mm\bv{b\leq c}&=
 \JJ_{u,\,v\in K}
 \bv{a\leq u}_P\mm\bv{u\leq b}_Q\mm\bv{b\leq v}_Q\mm\bv{v\leq c}_P\\
 &\leq\bv{a\leq c}.
 \end{align*}

This completes the proof that $R$ is a $D$-valued poset. Furthermore, it is
obvious that $x\leq y$ implies that $\bv{x\leq y}=1$, for all $x$, $y\in R$.

We now verify items (iv) and (v) of the definition
of a \CMPL{D}\ (see Definition~\ref{D:DvalMPL}). Let us
verify (iv). So, let $a\in R$ and $X\in\fine R$ such that $a=\JJ X$. We
verify that
 \begin{equation}\label{Eq:alebJJXleb}
 \bv{a\leq b}=\MM_{x\in X}\bv{x\leq b},\qq\text{for all }b\in R.
 \end{equation}
Without loss of generality, $\set{a}\uu X\ci P$ and $a=\JJ X$ in $P$.
If $b\in P$, then all the Boolean values involved in \eqref{Eq:alebJJXleb}
are computed in $P$, so \eqref{Eq:alebJJXleb} follows from the fact that $P$
is a \CMPL{D}.

So, suppose that $b\in Q$. For $x\in X$, $x\leq a$, thus $\bv{x\leq a}=1$,
so $\bv{a\leq b}=\bv{x\leq a}\mm\bv{a\leq b}\leq\bv{x\leq b}$, thus
 \[
 \bv{a\leq b}\leq\MM_{x\in X}\bv{x\leq b}.
 \]
To prove the converse inequality, we observe that, since $D$ is distributive
and both $X$ and $K$ are nonempty and finite, the following equalities
 \begin{align*}
 \MM_{x\in X}\bv{x\leq b}&=
 \MM_{x\in X}\JJ_{y\in K}\bv{x\leq y}_P\mm\bv{y\leq b}_Q\\
 &=\JJ_{\gn\colon X\to K}\MM_{x\in X}
 \bv{x\leq\gn(x)}_P\mm\bv{\gn(x)\leq b}_Q
 \end{align*}
hold. Thus, to complete the proof of \eqref{Eq:alebJJXleb}, it suffices to
prove that for any map $\gn\colon X\to K$, the inequality
 \begin{equation}\label{Eq:xgn(x)bleab}
 \MM_{x\in X}\bv{x\leq\gn(x)}_P\mm\bv{\gn(x)\leq b}_Q\leq\bv{a\leq b}
 \end{equation}
holds. Since $K$ is a lattice, $c=\JJ_{x\in X}\gn(x)$ is defined in $K$.
The inequality $\bv{x\leq\gn(x)}_P\leq\bv{x\leq c}_P$ holds for any $x\in X$
(because $\gn(x)\leq c$), and, since
$c=\JJ_{x\in X}\gn(x)$ in $Q$ and $Q$ is a \CMPL{D},
 \[
 \MM_{x\in X}\bv{\gn(x)\leq b}_Q=\bv{c\leq b}_Q.
 \]
It follows that
 \begin{align*}
 \MM_{x\in X}\bv{x\leq\gn(x)}_P\mm\bv{\gn(x)\leq b}_Q
 &=\MM_{x\in X}\bv{x\leq\gn(x)}_P\mm\MM_{x\in X}\bv{\gn(x)\leq b}_Q\\
 &\leq\MM_{x\in X}\bv{x\leq c}_P\mm\bv{c\leq b}_Q\\
 &=\bv{a\leq c}_P\mm\bv{c\leq b}_Q\\
\intertext{(because $c=\JJ_{x\in X}\gn(x)$ in $P$ and $P$
is a \CMPL{D})}
 &\leq\bv{a\leq b}.
 \end{align*}
This completes the proof of \eqref{Eq:alebJJXleb}. Therefore, $R$ satisfies
(iv) of Definition~\ref{D:DvalMPL}. The proof of (v) of
Definition~\ref{D:DvalMPL} is dual.

So, $R$ is a \CMPL{D}. The fact that both $f'$ and $g'$
are isometries is trivial.

If $S$ is a \CMPL{D}\ and $\oll{f}\colon P\to S$ and
$\ol{g}\colon Q\to S$ are uniform maps such that
$\oll{f}\circ f=\ol{g}\circ g$, then there exists a unique map
$h\colon R\to S$ such that $h\circ f'=\oll{f}$ and $h\circ g'=\ol{g}$,
namely, $h$ is defined by the rule
 \[
 h(x)=\oll{f}(x),\q\text{for any }x\in P,\qq\text{and}\q
 h(x)=\ol{g}(x),\q\text{for any }x\in Q.
 \]
Since $R$ is the pushout of~$P$ and $Q$ above $K$ in the category of partial
lattices (see Proposition~\ref{P:PushStV}), $h$ is a homomorphism of partial
lattices. It remains to prove that $\bv{x\leq y}_R\leq\bv{h(x)\leq h(y)}_S$,
for all $x$, $y\in R$. If $x$, $y\in P$, then
 \[
 \bv{h(x)\leq h(y)}_S=\bv{\oll{f}(x)\leq\oll{f}(y)}_S\geq\bv{x\leq y}_P
 =\bv{x\leq y}_R.
 \]
A similar proof applies to the case where $x$, $y\in Q$.
If $x\in P$ and $y\in Q$, we compute:
 \begin{align*}
 \bv{h(x)\leq h(y)}_S&=\bv{\oll{f}(x)\leq\ol{g}(y)}_S\\
 &\geq\JJ_{z\in K}
 \bv{\oll{f}(x)\leq\oll{f}(z)}_S\mm\bv{\ol{g}(z)\leq\ol{g}(y)}_S\\
\intertext{(because $\oll{f}\res_K=\ol{g}\res_K$)}
 &\geq\JJ_{z\in K}
 \bv{x\leq z}_P\mm\bv{z\leq y}_Q\\
\intertext{(because $\oll{f}$ and $\ol{g}$ are uniform)}
 &=\bv{x\leq y}_R.
 \end{align*}
The proof is similar in case $x\in Q$ and $y\in P$. Therefore, $h$ is
uniform. So, $R$ is the pushout of~$P$ and $Q$ above $K$ in the category of
\CMPL{D}s.
\end{proof}

Of course, in accordance with Notation~\ref{Not:PamalgQ}, we shall also
write $R=P\amalg_KQ$ in the context of Proposition~\ref{P:DPushStV}.

In Lemmas \ref{L:bv(x=y)} and \ref{L:ainPandQ}, let $\E{D}=\vv<K,P,Q,f,g>$
be a standard V-formation of \CMPL{D}s (so that $K$ is a finite lattice).
We put $R=P\amalg_KQ$, endowed with its structure of \CMPL{D}
described in Proposition~\ref{P:DPushStV}.

\begin{lemma}\label{L:bv(x=y)}
The equality
 \[
 \bv{x=y}=\JJ_{z\in K}\bv{x=z}\mm\bv{z=y}
 \]
holds, for all $x$, $y\in R$ such that either $x\in P$ and $y\in Q$, or
$x\in Q$ and $y\in P$.
\end{lemma}

\begin{proof}
We assume, for example, that $x\in P$ and $y\in Q$. We compute:
 \begin{align*}
 \bv{x=y}&=\bv{x\leq y}\mm\bv{y\leq x}\\
 &=\JJ_{u,\,v\in K}\bv{x\leq u}\mm\bv{u\leq y}\mm\bv{y\leq v}\mm\bv{v\leq x}.
 \end{align*}
For $u$, $v\in K$, $\bv{x\leq u}\mm\bv{v\leq x}\leq\bv{v\leq u}$, while
$\bv{u\leq y}\mm\bv{y\leq v}\leq\bv{u\leq v}$. Therefore,
 \begin{align*}
 \bv{x=y}&=\JJ_{u,\,v\in K}
 \bv{x\leq u}\mm\bv{u\leq y}\mm\bv{y\leq v}\mm\bv{v\leq x}\mm\bv{u=v}\\
 &=\JJ_{u,\,v\in K}
 \bv{x\leq u}\mm\bv{u\leq y}\mm\bv{y\leq u}\mm\bv{u\leq x}\mm\bv{u=v}\\
 &=\JJ_{u\in K}\bv{x=u}\mm\bv{u=y},
 \end{align*}
which concludes the proof.
\end{proof}

\begin{lemma}\label{L:ainPandQ}\hfill
\begin{enumerate}
\item For any $a\in P$ and any $Y\in\fine Q$, the inequality
$\bv{a\in Y}\leq\bv{a\in K}$ holds.

\item For any $X\in\fine P$ and any $Y\in\fine Q$, the inequality
$\bv{X\ci Y}\leq\bv{X\ci K}$ holds.
\end{enumerate}
\end{lemma}

\begin{proof}
(i) For all $y\in Y$, we compute, using Lemma~\ref{L:bv(x=y)}:
 \[
 \bv{a=y}=\JJ_{z\in K}\bv{a=z}\mm\bv{z=y}\leq\bv{a\in K},
 \]
since $\bv{a=z}\leq\bv{a\in K}$. This holds for all $y\in Y$, so the
conclusion follows.

(ii) is a trivial consequence of (i).
\end{proof}

We are now able to prove the following fundamental result:

\begin{proposition}\label{P:PQcovR}
Let $\E{D}=\vv<K,P,Q,f,g>$ be a standard V-formation of \CMPL{D}s. Put
$R=P\amalg_KQ$. If $P$ and $Q$ are finitely covering, then $R$ is finitely
covering.
\end{proposition}

\begin{proof}
We prove, for example, that the domain of the join in $R$ is finitely
covering. So let $Z\in\fine R$, we prove that $Z$ has a $\dom\JJ_R$-cover.

Write $Z=X\uu Y$, where $X\in\fin P$ and $Y\in\fin Q$. By
Lemma~\ref{L:CollDCov}, applied within~$P$,
there exists a finite subset $\setm{X_i}{i<m}$ (with
$m>0$) of~$\dom\JJ_P$ such that
 \[
 \bv{U\ci X\uu K}\leq\JJ_{i<m}\bv{U=X_i},\qq\text{for all }U\in\dom{\JJ}_P.
 \]
Therefore, for any $U\in\dom\JJ_P$,
 \begin{align}
 \bv{Z=U}&\leq\bv{U\ci X\uu Y}\mm\bv{Y\ci U}\notag\\
 &\leq\bv{U\ci X\uu Y}\mm\bv{Y\ci K}\notag\\
\intertext{(by Lemma~\ref{L:ainPandQ}(ii))}
 &=\bv{U\ci X\uu Y}\mm\bv{X\uu Y\ci X\uu K}\notag\\
 &\leq\bv{U\ci X\uu K}\notag\\
 &\leq\JJ_{i<m}\bv{U=X_i}.\label{Eq:Z=UindomP}
 \end{align}
Similarly, there exists a finite subset $\setm{X_i}{m\leq i<m+n}$ (with $n>0$)
of~$\dom\JJ_Q$ such that
 \begin{equation}\label{Eq:Z=UindomQ}
 \bv{Z=U}\leq\JJ_{j<n}\bv{U=X_{m+j}},\qq\text{for all }U\in\dom{\JJ}_Q.
 \end{equation}
Therefore, by \eqref{Eq:Z=UindomP} and \eqref{Eq:Z=UindomQ}
and since $\dom\JJ_R=\dom\JJ_P\uu\dom\JJ_Q$,
 \[
 \bv{Z=U}\leq\JJ_{i<m+n}\bv{U=X_i},\qq\text{for all }U\in\dom{\JJ}_R.
 \]
So, $\setm{X_i}{i<m+n}$ is a $\dom\JJ_R$-cover of~$Z$.
The proof for the meet is dual.
\end{proof}

In the context of Proposition~\ref{P:PQcovR}, since $R$ is a finitely
covering \CMPL{D}, it can, by Proposition~\ref{P:DmeasDval}, be
canonically extended into a $D$-valued partial lattice, namely, $\tR$.
Hence, for every prime filter $G$ of~$D$, $K/G$ is a finite lattice, and
$P/G$, $Q/G$, and $R/G$ are partial lattices, see
Proposition~\ref{P:P/GPartLatt}. Furthermore, by Lemma~\ref{L:IsomGotoG},
the canonical maps $K/G\hookrightarrow P/G$, $K/G\hookrightarrow Q/G$,
$P/G\hookrightarrow R/G$, and $Q/G\hookrightarrow R/G$ are \emph{embeddings}
of partial lattices. The question whether they form a pushout (in the
category of partial lattices) is answered naturally:

\begin{proposition}\label{P:Pushout/G}
Let $\vv<K,P,Q>$ be a standard V-formation of \CMPL{D}s,
with $P$ and $Q$ finitely covering. Then $R/G=P/G\amalg_{K/G}Q/G$, for any
prime filter $G$ of~$D$.
\end{proposition}

\begin{proof}
{F}rom $R=P\uu Q$ follows trivially that $R/G=P/G\uu Q/G$.

Let $\xa\in P/G$ and $\xb\in Q/G$ such that $\xa\leq\xb$. Pick $a\in\xa$ and
$b\in\xb$, then $\bv{a\leq b}\in G$, thus, by the definition of the Boolean
values in $R$ (see Proposition~\ref{P:DPushStV}), there exists $c\in K$ such
that $\bv{a\leq c}\mm\bv{c\leq b}\in G$. Hence, for $\xc=\cls(c,G)$, we
obtain that $\xa\leq\xc\leq\xb$. A similar statement, with $P$ and $Q$
exchanged, holds. This implies that $K/G=P/G\ii Q/G$ and that the ordering
of $R/G$ is the same as the ordering of $P/G\amalg_{K/G}Q/G$, see
Proposition~\ref{P:PushStV}. Finally, the fact that $R/G$ and
$P/G\amalg_{K/G}Q/G$ have the same join and meet operations follows easily
from Lemma~\ref{L:PDmeasP/G}.
\end{proof}

\section{\IDM- and \FIM-samples for $P\amalg_KQ$}\label{S:IdFiM}

In this section, we shall fix a distributive lattice $D$ with unit and a
standard V-formation $\vv<K,P,Q>$ of \CMPL{D}s, with $P$
and $Q$ finitely covering. We put $R=P\amalg_KQ$. By Propositions
\ref{P:DPushStV} and \ref{P:PQcovR}, $R$ is a finitely covering
\CMPL{D}.

This section will be devoted to the proof of the following result:

\begin{proposition}\label{P:AmPresIDM}
Suppose that $P$ and $Q$ have \IDM\ (resp., \FIM). Then $R$ has \IDM\
(resp., \FIM).
\end{proposition}

\begin{proof}
We provide a proof for \IDM; the proof for \FIM\ is dual.

Let $a$, $b\in R$, we shall find an \IDM-sample of $\set{a,b}$. Up to
symmetry between $P$ and $Q$, there are two cases to consider.
\smallskip

\noindent\textit{\textbf{Case~1.}} $a$, $b\in P$.

Let $U$ be an \IDM-sample of $\set{a,b}$ in $P$, we prove that $U$ is also
an \IDM-sample of $\set{a,b}$ in $R$. This amounts to proving the inequality
 \begin{equation}\label{Eq:MMsama,b}
 \bv{x\leq a}\mm\bv{x\leq b}\leq
 \JJ_{u\in U}\bv{x\leq u}\mm\bv{u\leq a}\mm\bv{u\leq b},
 \qq\text{for all }x\in R
 \end{equation}
(the converse inequality of \eqref{Eq:MMsama,b} is trivial). First, for
$x\in P$, all the Boolean values involved in \eqref{Eq:MMsama,b} are
Boolean values in $P$, so, since $U$ is an \IDM-sample of $\set{a,b}$ in $P$,
\eqref{Eq:MMsama,b} holds.

Now suppose that $x\in Q$. We compute:
 \begin{align*}
 \bv{x\leq a}\mm\bv{x\leq b}&=\JJ_{u,\,v\in K}
 \bv{x\leq u}\mm\bv{u\leq a}\mm\bv{x\leq v}\mm\bv{v\leq b}\\
 &=\JJ_{w\in K}\bv{x\leq w}\mm\bv{w\leq a}\mm\bv{w\leq b}\\
\intertext{(by putting $w=u\mm v$ in $K$)}
 &=\JJ_{w\in K,\ u\in U}
 \bv{x\leq w}\mm\bv{w\leq u}\mm\bv{u\leq a}\mm\bv{u\leq b}\\
\intertext{(because $U$ is an \IDM-sample of $\set{a,b}$ in $P$)}
 &\leq\JJ_{u\in U}\bv{x\leq u}\mm\bv{u\leq a}\mm\bv{u\leq b},
 \end{align*}
so \eqref{Eq:MMsama,b} is established in this case.
\smallskip

\noindent\textit{\textbf{Case~2.}} $a\in P$, $b\in Q$.

Let $U$ (resp., $V$) be a common \IDM-sample in $P$ (resp., in $Q$) of all
pairs of the form $\set{a,z}$ (resp., $\set{b,z}$), for $z\in K$. We put
$W=U\uu V$, and we prove that $W$ is an \IDM-sample of $\set{a,b}$ in $R$.

For any $x\in P$, we compute:
 \begin{align*}
 \bv{x\leq a}\mm\bv{x\leq b}&=
 \JJ_{z\in K}\bv{x\leq a}\mm\bv{x\leq z}\mm\bv{z\leq b}\\
 &=\JJ_{z\in K,\ u\in U}
 \bv{x\leq u}\mm\bv{u\leq a}\mm\bv{u\leq z}\mm\bv{z\leq b}\\
\intertext{(because $U$ is an \IDM-sample of all pairs $\set{a,z}$ for
$z\in K$)}
 &\leq\JJ_{u\in U}\bv{x\leq u}\mm\bv{u\leq a}\mm\bv{u\leq b}.
 \end{align*}
Similarly, we can obtain that
 \[
 \bv{x\leq a}\mm\bv{x\leq b}\leq
 \JJ_{v\in V}\bv{x\leq v}\mm\bv{v\leq a}\mm\bv{v\leq b},
 \]
for any $x\in Q$. Hence, 
 \[
 \bv{x\leq a}\mm\bv{x\leq b}\leq
 \JJ_{w\in W}\bv{x\leq w}\mm\bv{w\leq a}\mm\bv{w\leq b},
 \]
for any $x\in R$.
\end{proof}

\section{\IDJ- and \FIJ-samples in $P\amalg_KQ$}\label{S:IdFiJ}

In this section, we shall fix, as in Section~\ref{S:IdFiM}, a distributive
lattice $D$ with unit and a standard V-formation $\vv<K,P,Q>$ of
\CMPL{D}s, with
$P$ and $Q$ finitely covering. We put $R=P\amalg_KQ$. By Propositions
\ref{P:DPushStV} and \ref{P:PQcovR}, $R$ is a finitely covering
\CMPL{D}.

This section will be devoted to the proof of the following result:

\begin{proposition}\label{P:AmPresIDJ}
Suppose that $P$ and $Q$ have \IDJ\ (resp., \FIJ). Then $R$ has \IDJ\
(resp., \FIJ).
\end{proposition}

\begin{proof}
We provide a proof for \IDJ; the proof for \FIJ\ is dual.

Let $Z\in\fine R$. We put $X=Z\ii P$ and $Y=Z\ii Q$. Observe that $Z=X\uu Y$.

Let $X^*$ be a common \IDJ-sample of all subsets of $X\uu K$ in $P$.
Symmetrically, let $Y^*$ be a common \IDJ-sample of all subsets of $Y\uu K$
in $Q$. Let $m$ be a common index for both samples, see
Definition~\ref{D:IDFIJsample}.

We denote by $h$ the \emph{height} of~$K$, and we put $k=(h+2)m+h+1$. We
shall prove the following assertion:
 \begin{equation}\label{Eq:SampleStat}
 Z^*=X^*\uu Y^*\text{ is an \IDJ-sample of }Z\text{ in }R,
 \text{ with index }k.
 \end{equation}

Let $G$ be a prime filter of~$D$. We put $\ol{T}=T/G$, for every subset $T$
of~$R$. Recall that $\ol{R}=\ol{P}\amalg_{\ol{K}}\ol{Q}$, see
Proposition~\ref{P:Pushout/G}.

We define $I_0\ci\ol{P}$ and $J_0\ci\ol{Q}$ as follows:
 \begin{align*}
 I_0&=\Id_m^{\ol{P}}(\ol{X},\ol{X^*});\\
 J_0&=\Id_m^{\ol{Q}}(\ol{Y},\ol{Y^*}).
 \end{align*}
Of course, the superscript $\ol{P}$ (or $\ol{Q}$) on the math operator $\Id$
indicates in which partial lattice the $\Id_n(U,V)$ function (see
Definition~\ref{D:IdFiln01}) is computed.

Since $m$ is an index for the \IDJ-sample $X^*$ of~$X$, it follows from
Lemma~\ref{L:TruthIdFiln} that
 \[
 \Id_m^{\ol{P}}(\ol{X},\ol{X^*})=\Id_{m+1}^{\ol{P}}(\ol{X},\ol{T}),
 \]
for every $T\in\fine P$ such that $X^*\ci T$. In particular, $I_0$ is an
ideal of~$P$ (it is empty if $X=\es$). Similarly, $J_0$ is an ideal of~$Q$.

\setcounter{claim}{0}

\begin{claim}
Assume that $(I_0\uu J_0)\ii\ol{K}=\es$. Then $I_0\uu J_0$ is an ideal of
$\ol{R}$. Furthermore, $I_0\uu J_0=\Id_n^{\ol{R}}(\ol{Z},\ol{T})$ holds for
all $n\geq m$ and all $T\ce Z^*$ in $\fine R$.
\end{claim}

\begin{cproof}
Since $\JJ_{\ol{R}}=\JJ_{\ol{P}}\uu\JJ_{\ol{Q}}$ and since $I_0$ (resp.,
$J_0$) is an ideal of~$\ol{P}$ (resp., $\ol{Q}$), $I_0\uu J_0$ is closed under
$\JJ_{\ol{R}}$. By the assumption that $(I_0\uu J_0)\ii\ol{K}=\es$, no pair
of $I_0\times\ol{Q}$ and $\ol{P}\times J_0$ is comparable, so, since $I_0$
(resp., $J_0$) is a lower subset of~$\ol{P}$ (resp., $\ol{Q}$), $I_0\uu J_0$
is a lower subset of $\ol{R}$. Every ideal of~$\ol{R}$ that contains $\ol{Z}$
contains $I_0\uu J_0$, thus, since $\ol{Z}\ci I_0\uu J_0$, $I_0\uu J_0$ is
the ideal of~$\ol{R}$ generated by $\ol{Z}$. The second part of the statement
of Claim~1 follows immediately.
\end{cproof}

Now we assume that $(I_0\uu J_0)\ii\ol{K}$ is nonempty. Since it is a
nonempty subset of the finite lattice $\ol{K}$, it admits a supremum, that
we denote by $c_0$. We observe that both $I_0$ and $J_0$ are contained in
$\Id_m^{\ol{R}}(\ol{Z},\ol{Z^*})$, thus
$c_0\in\Id_{m+1}^{\ol{R}}(\ol{Z},\ol{Z^*})$. We extend this construction by
defining inductively $I_n$, $J_n$, and $c_n$, for any $n<\go$, by
 \begin{align*}
 I_{n+1}&=\Id_m^{\ol{P}}(\ol{X}\uu\set{c_n},\ol{X^*}),\\
 J_{n+1}&=\Id_m^{\ol{Q}}(\ol{Y}\uu\set{c_n},\ol{Y^*}),\\
 c_n&=\JJ((I_n\uu J_n)\ii K),
 \end{align*}
for all $n<\go$.

Since $X^*$ is a common \IDJ-sample of all subsets of $X\uu K$ in $P$, with
index $m$, $I_{n+1}$ is an ideal of~$P$, for all $n<\go$. Similarly,
$J_{n+1}$ is an ideal of~$Q$. Furthermore, $c_n\in I_{n+1}\ii J_{n+1}$, thus
$c_n\leq c_{n+1}$. Therefore, $I_n\ci I_{n+1}$ and $J_n\ci J_{n+1}$ for all
$n$.

An easy inductive generalization of the argument above showing that
$c_0\in\Id_{m+1}^{\ol{R}}(\ol{Z},\ol{Z^*})$ leads to the following:

\begin{claim}
$I_n\uu J_n\ci\Id_{(n+1)m+n}^{\ol{R}}(\ol{Z},\ol{Z^*})$, for all $n<\go$.
\end{claim}

For $n<\go$, if $c_0<c_1<\cdots<c_n$, then $n<\hgt(\ol{K})\leq\hgt K=h$. If
$c_n=c_{n+1}$, then $c_n=c_l$ for all $l\geq n$. It follows from this that
$c_h=c_{h+1}$, thus $I_{h+1}=I_{h+2}$ and $J_{h+1}=J_{h+2}$. We put $c=c_h$,
$I=I_{h+1}$, and $J=J_{h+1}$.

\begin{claim}
$I\uu J$ is an ideal of~$\ol{R}$.
\end{claim}

\begin{cproof}
Since $I$ is an ideal of~$\ol{P}$ and $J$ is an ideal of~$\ol{Q}$, $I\uu J$
is closed under $\JJ_{\ol{R}}$. Now we prove that $I\uu J$ is a lower subset
of~$\ol{R}$. Since $c_0\in I_1\ii J_1\ci I\ii J$, both $I\ii\ol{K}$ and
$J\ii\ol{K}$ are nonempty, hence there are elements $a$ and $b$ of~$\ol{K}$
defined by
 \[
 a=\JJ(I\ii\ol{K})\qq\text{and}\qq b=\JJ(J\ii\ol{K}).
 \]
Since $c=c_h\in I\ii J\ii\ol{K}$, $c\leq a$ and $c\leq b$. On the other hand,
 \[
 c=\JJ((I\uu J)\ii\ol{K})=a\jj b,
 \]
whence $c=a=b$. So we have established that
 \begin{equation}\label{Eq:cIKJK}
 c=\JJ(I\ii\ol{K})=\JJ(J\ii\ol{K}).
 \end{equation}
Now let $x\in\ol{R}$, $y\in I\uu J$ such that $x\leq y$, we prove that
$x\in I\uu J$. By symmetry, we may assume that $y\in I$. If $x\in\ol{P}$,
then, since $I$ is a lower subset of~$\ol{P}$, $x\in I$ and we are done. If
$x\in\ol{Q}$, then there exists $z\in\ol{K}$ such that $x\leq z\leq y$.
Since $z\in\ol{K}$ and $z\leq y\in I$, $z\leq c$ by \eqref{Eq:cIKJK}, so
$x\leq c$. But $c\in J$ and $J$ is a lower subset of~$\ol{Q}$, thus
$x\in J$, and we are done again.

So $I\uu J$ is a lower subset of~$\ol{R}$, hence an ideal of~$\ol{R}$.
\end{cproof}

By Claims 2 and 3, $I\uu J=\Id_{(h+2)m+h+1}^{\ol{R}}(\ol{Z},\ol{T})$, for
all finite $T$ containing $Z^*$, is the ideal of~$\ol{R}$ generated by
$\ol{Z}$. In particular,
 \begin{equation}\label{Eq:Id(h+2)m+h+1}
 \Id_{(h+2)m+h+1}^{\ol{R}}(\ol{Z},\ol{Z^*})=
 \Id_{(h+2)m+h+2}^{\ol{R}}(\ol{Z},\ol{T}),
 \end{equation}
for all finite $T$ containing $Z^*$. This also holds in the context of
Claim~1, since one can, in that case, replace $(h+2)m+h+1$ by $m$.
Therefore,
\eqref{Eq:Id(h+2)m+h+1} holds for \emph{every} prime filter $G$ of~$D$.
By Lemma~\ref{L:EquivIdFilSG}, this proves \eqref{Eq:SampleStat}.
\end{proof}

\part{Congruence amalgamation with distributive target}
\label{Pt:CongDistr}

\section{Proof of Theorem B}

We first observe the following obvious restatement of Theorem~B in terms of
\MPL{D}s:

\begin{all}{Theorem B}
Let $D$ be a distributive lattice with zero.
Let $\vv<K,\gl>$, $\vv<P,\gm>$, and $\vv<Q,\gn>$ be \MPL{D}s, with $K$ a
finite lattice and each of $P$ and $Q$ either a finite partial lattice or a
lattice. Let $f\colon\vv<K,\gl>\to\vv<P,\gm>$ and
$g\colon\vv<K,\gl>\to\vv<Q,\gn>$ be homomorphisms.

Then there exist a \ML{D} $\vv<L,\gf>$ and
homomorphisms $\oll{f}\colon\vv<P,\gm>\to\vv<L,\gf>$ and
$\ol{g}\colon\vv<Q,\gn>\to\vv<L,\gf>$ such that
$\oll{f}\circ f=\ol{g}\circ g$. Furthermore, the construction can be done in
such a way that the following additional properties hold:

\begin{enumerate}
\item $L$ is generated, as a lattice, by~$\oll{f}[P]\uu\ol{g}[Q]$.

\item The map $\gf$ isolates $0$.
\end{enumerate}
\end{all}

By Proposition~\ref{P:DvalMPL}, if $D$ is a distributive lattice with
zero, then the notions of \MPL{D} (see Definition~\ref{D:D0PL}) and of
\CMPL{D^\rd} (see Definition~\ref{D:DvalMPL}) are, essentially, equivalent.
It is, in fact, easy to see that this is a category equivalence. The
corresponding notion of \emph{morphism of \MPL{D}} is given by the following
very easy result:

\begin{lemma}\label{L:DualDPL01}
Let $D$ be a distributive lattice with zero,
let $\vv<P,\gm>$ and $\vv<Q,\gn>$ be \MPL{D}s, let
$f\colon P\hookrightarrow Q$ be an embedding of partial lattices. Then the
following are equivalent:
\begin{enumerate}
\item The equality $\gn\circ\Conc f=\gm$ holds.

\item If $\vv<P,\gm>$ and $\vv<Q,\gn>$ are viewed as \CMPL{D^\rd}s, then
$f$ is an isometry (see Definition~\tup{\ref{D:IsomDmeas}}).
\end{enumerate}
\end{lemma}

\begin{proof}
Endow each of the structures $P$ and $Q$ with its map $\bv{\ghost}$, with
target $D^\rd$. An explicit definition of $\bv{\ghost}_P$ and $\bv{\ghost}_Q$
is the following:
 \begin{align*}
 \bv{x\leq y}_P&=\gm(\gQ_P^+(x,y)),&&\text{for all }x,\,y\in P;\\
 \bv{x\leq y}_Q&=\gn(\gQ_Q^+(x,y)),&&\text{for all }x,\,y\in Q.
 \end{align*}
Hence, for $x$, $y\in P$,
 \begin{align}
 \bv{f(x)\leq f(y)}_Q&=\gn(\gQ_Q^+(f(x),f(y))\notag\\
 &=\gn\circ(\Conc f)(\gQ_P^+(x,y)).\label{Eq:gngQ+P(x,y)}
 \end{align}
But the principal congruences $\gQ_P^+(x,y)$, for $x$, $y\in P$, generate
the \jz-semilattice $\Conc P$. Hence, $\gn\circ\Conc f=\gm$ if{f} both maps
$\gn\circ\Conc f$ and $\gm$ agree on all principal congruences of~$P$, that
is, by \eqref{Eq:gngQ+P(x,y)}, $\bv{f(x)\leq f(y)}_Q=\bv{x\leq y}_P$, for all
$x$, $y\nobreak\in\nobreak P$.
\end{proof}

\begin{definition}\label{D:EmbDPL}
Let $D$ be a distributive lattice with zero, let $\vv<P,\gm>$ and
$\vv<Q,\gn>$ be \MPL{D}s. A \emph{homomorphism} from $\vv<P,\gm>$ to
$\vv<Q,\gn>$ is a homomorphism $f\colon P\to Q$ of partial
lattices such that $\gn\circ\Conc f=\gm$. If, in addition, $f$ is an
embedding of partial lattices, we say that $f$ is an \emph{embedding} of
\MPL{D}s.
\end{definition}

\begin{definition}\label{D:Kernel}
Let $D$ be a distributive lattice with zero, let $\vv<P,\gm>$ be a \MPL{D}.
The \emph{kernel} of $\vv<P,\gm>$ is the congruence $\gq$ of~$P$ defined by
the rule
 \[
 x\leq_\gq y\text{ if{f} }\gm\gQ_P^+(x,y)=0,\qquad\text{for all }x,\,y\in P.
 \]
The \emph{kernel projection} of $\vv<P,\gm>$ is the canonical projection from
$P$ onto $P/{\gq}$.
\end{definition}

In other words, the kernel of $\vv<P,\gm>$ is the largest congruence $\gq$
of~$P$ such that $\gm(\gx)=0$ for all $\gx\leq\gq$ in $\Conc P$. In
particular, $\vv<P,\gm>$ is proper (see Definition~\ref{D:D0PL}) if{f} its
kernel projection is trivial. 

The following lemma states that the kernel projection of $\vv<P,\gm>$ is the
universal projection of $\vv<P,\gm>$ onto a proper \MPL{D}:

\begin{lemma}\label{L:Kernel}
Let $D$ be a distributive lattice with zero, let $\vv<P,\gm>$ be a \MPL{D}.
We denote by $p\colon P\twoheadrightarrow P'$ the kernel projection of~$P$.
Then there exists a unique \jzh\ $\gm'\colon\Conc P'\to D$ such that
$\gm'\circ\Conc p=\gm$. Furthermore, the following assertions hold:
\begin{enumerate}
\item $\vv<P',\gm'>$ is a proper \MPL{D}.

\item For every proper \MPL{D} $\vv<Q,\gn>$ and every homomorphism
$f\colon\vv<P,\gm>\to\vv<Q,\gn>$, there exists a unique homomorphism
$f'\colon\vv<P',\gm'>\to\vv<Q,\gn>$ such that $f'\circ p=f$, and $f'$ is an
embedding of \MPL{D}s (see Definition~\textup{\ref{D:EmbDPL}}).
\end{enumerate}
\end{lemma}

\begin{proof}
(i) The compact congruences of $P'=P/\gq$ are exactly the congruences of the
form $\ga\jj\gq/\gq$, where $\ga$ is a compact congruence of $P$. By the
definition of $\gq$, $\ga\jj\gq\leq\gb\jj\gq$ implies that
$\gm(\ga)\leq\gm(\gb)$, for all $\ga$, $\gb\in\Conc P$. Hence we can define
a \jzh\ $\gm'\colon\Conc P'\to D$ by the rule
 \[
 \gm'(\ga\jj\gq/\gq)=\gm(\ga),\quad\text{for all }\ga\in\Conc P.
 \]
Observe that $\gm'\circ\Conc p=\gm$. The uniqueness assertion about
$\gm'$ follows from the surjectivity of the map $\Conc p$.

(ii) For all $x$, $y\in P$, $\cls(x,\gq)\leq\cls(y,\gq)$ if{f}
$\gm\gQ_P^+(x,y)=0$, that is, $\gn\gQ_Q^+(f(x),f(y))=0$, or, since
$\vv<Q,\gn>$ is proper, $f(x)\leq f(y)$. This makes it possible to define an
embedding $f'\colon P'\hookrightarrow Q$ of partial lattices by the rule
$f'(\cls(x,\gq))=f(x)$, for all $x\in P$, and $f'\circ p=f$. Observe that
$\gn\circ\Conc f'=\gm'$. The uniqueness assertion about $f'$ follows from the
surjectivity of the map $p$.
\end{proof}

\begin{proposition}\label{P:SquareAmalg}
Let $D$ be a distributive lattice with zero.
Let $\vv<K,\gl>$, $\vv<P,\gm>$, and $\vv<Q,\gn>$ be \MPL{D}s, with $K$ a
finite lattice and $P$, $Q$ balanced.
Let $f\colon\vv<K,\gl>\to\vv<P,\gm>$ and
$g\colon\vv<K,\gl>\to\vv<Q,\gn>$ be homomorphisms.

Then there exists a proper \ML{D} $\vv<L,\gf>$, together with
homomorphisms $\oll{f}\colon\vv<P,\gm>\to\vv<L,\gf>$ and
$\ol{g}\colon\vv<Q,\gn>\to\vv<L,\gf>$, such that
$\oll{f}\circ f=\ol{g}\circ g$ and $L$ is generated, as a lattice,
by~$\oll{f}[P]\uu\ol{g}[Q]$.
\end{proposition}

\begin{proof}
We first consider the case where $f$ and $g$ are embeddings and
$\gl$, $\gm$, and $\gn$ isolate $0$.
We view $\vv<K,\gl>$, $\vv<P,\gm>$, and $\vv<Q,\gn>$ as \CMPL{D^\rd}s.
By Lemma~\ref{L:DV-FormSt}, we can assume without loss of generality that
$\vv<K,P,Q,f,g>$ is a standard V-formation. Now we put $R=P\amalg_KQ$, as
defined in Proposition~\ref{P:DPushStV}, with the corresponding embeddings
$f'$ and $g'$. Let $\vv<R,\gr>$ be the corresponding \MPL{D}.
{}From the fact that both $\gm$ and $\gn$ isolate $0$ and the description of
$R$ (Proposition~\ref{P:DPushStV}) follows that $\gr$ isolates $0$.
By Propositions \ref{P:PQcovR}, \ref{P:AmPresIDM}, and
\ref{P:AmPresIDJ}, $R$ is balanced.
By Theorem~A, there exists a \jzh\ $\gf'\colon\Conc\FL(R)\to D$
such that $\gf'\circ\Conc j_R=\gr$. We denote by
$p\colon\FL(R)\twoheadrightarrow L$ the kernel projection of
$\vv<\FL(R),\gf'>$, see Definition~\ref{D:Kernel}, and we put $j=p\circ j_R$.
By Lemma~\ref{L:Kernel}, there exists a unique \jzh\ $\gf\colon\Conc L\to D$
such that $\gf\circ\Conc p=\gf'$, and $\vv<L,\gf>$ is proper. Furthermore,
 \[
 \gf\circ\Conc j=\gf\circ\Conc p\circ\Conc j_R=\gf'\circ\Conc j_R=\gr.
 \]
We put $\oll{f}=j\circ f'$ and $\ol{g}=j\circ g'$. {F}rom
$f'\circ f=g'\circ g$ follows that $\oll{f}\circ f=\ol{g}\circ g$. Since
$j_R[R]$ generates $\FL(R)$, $\oll{f}[P]\uu\ol{g}[Q]$ generates $L$.
Since $j$ is a homomorphism from $\vv<R,\gr>$ to $\vv<L,\gf>$ and
since both $\gr$ and $\gf$ isolate $0$, $j$ is an embedding, thus $\oll{f}$
and $\ol{g}$ are embeddings.

Now we consider the general case. Let
$h\colon\vv<K,\gl>\twoheadrightarrow\vv<K',\gl'>$,
$p\colon\vv<P,\gm>\twoheadrightarrow\vv<P',\gm'>$,
$q\colon\vv<Q,\gn>\twoheadrightarrow\vv<Q',\gn'>$ be the kernel projections.
By Lemma~\ref{L:Kernel}, there are embeddings
$f'\colon\vv<K',\gl'>\hookrightarrow\vv<P',\gm'>$ and
$g'\colon\vv<K',\gl'>\hookrightarrow\vv<Q',\gn'>$ such that
$f'\circ h=p\circ f$ and $g'\circ h=q\circ g$. By the result of the previous
paragraph, there exist a proper \ML{D} $\vv<L,\gf>$ and embeddings
$f''\colon\vv<P',\gm'>\hookrightarrow\vv<L,\gf>$ and
$g''\colon\vv<Q',\gn'>\hookrightarrow\vv<L,\gf>$ such that
$f''\circ f'=g''\circ g'$ and $L$ is generated by $f''[P']\uu g''[Q']$. We
put $\oll{f}=f''\circ p$ and $\ol{g}=g''\circ q$.
\end{proof}

\begin{remark}\label{Rk:EmbD}
In the context of Proposition~\ref{P:SquareAmalg}, we shall later make use
of the following simple fact: \emph{If $\vv<Q,\gn>$ is proper, then $\ol{g}$
is an embedding}.

Indeed, for all $x$, $y\in Q$,
 \begin{align*}
 \ol{g}(x)\leq\ol{g}(y)&\Rightarrow\gf\gQ_L^+(\ol{g}(x),\ol{g}(y))=0&&
 (\text{because }\gf\text{ isolates }0)\\
 &\Leftrightarrow(\gf\circ\Conc\ol{g})(\gQ_Q^+(x,y))=0&&\\
 &\Leftrightarrow\gn\gQ_Q^+(x,y)=0&&\\
 &\Leftrightarrow x\leq y&&(\text{because }\vv<Q,\gn>\text{ is proper}),
 \end{align*}
which proves our assertion.
\end{remark}

\begin{proof}[Proof of Theorem~B]
By Proposition~\ref{P:LatFinDS}, every \MPL{D} $\vv<R,\gr>$ such that
either $R$ is finite with nonempty join and meet operations or $R$ is a
lattice is balanced. Theorem~B follows immediately as a
particular case of Proposition~\ref{P:SquareAmalg}.
\end{proof}

\section{Saturation properties of \MPL{D}s}

We start with a definition.

\begin{definition}\label{D:ExtMPL}
Let $D$ be a distributive lattice with zero,
let $\vv<P,\gm>$ and $\vv<L,\gf>$ be \MPL{D}s, with $L$ a lattice.
We say that an embedding $f\colon\vv<P,\gm>\hookrightarrow\vv<L,\gf>$ is a
\emph{lower embedding} (resp., \emph{upper embedding}, \emph{internal
embedding}), if the filter (resp., ideal, convex sublattice) of $L$ generated
by $P$ equals~$L$.
\end{definition}

We refer to Definition~\ref{D:EmbDPL} for the definition of an embedding of
\MPL{D}s.

\begin{definition}\label{D:EqnSat}
Let $D$ be a distributive lattice with zero.
A proper \ML{D} $\vv<L,\gf>$ is \emph{saturated}
(resp., \emph{lower saturated}, \emph{upper saturated}, \emph{internally
saturated}), if for every embedding (resp., lower embedding, upper embedding,
internal embedding)
$e\colon\vv<K,\gl>\hookrightarrow\vv<P,\gm>$ of finite proper \MPL{D}s, with
$K$ a \emph{lattice}, and every homomorphism
$f\colon\vv<K,\gl>\to\vv<L,\gf>$, there exists a homomorphism
$g\colon\vv<P,\gm>\to\vv<L,\gf>$ such that $g\circ e=f$.
\end{definition}

\begin{proposition}\label{P:EmbSat}
Let $D$ be a distributive lattice with zero.
Every proper balanced \MPL{D} $\vv<P,\gf>$ admits an embedding (resp.,
a lower embedding, an upper embedding, an
internal embedding) into a saturated (resp., 
lower saturated, upper saturated, internally saturated) \ML{D} $\vv<L,\gy>$
such that $|L|=|P|+|D|+\aleph_0$.
\end{proposition}

\begin{proof}
A standard increasing chain argument. We present the proof for saturated, the
proofs for lower, upper, or internally saturated are similar. We put
$\gk=|P|+|D|+\aleph_0$. By Theorem~A, there exist a \ML{D} $\vv<K,\gz>$ and an
embedding $f\colon\vv<P,\gf>\hookrightarrow\vv<K,\gz>$. Furthermore, by
replacing $K$ by the image of its kernel projection (use
Lemma~\ref{L:Kernel}), we can suppose that $\vv<K,\gz>$ is proper. Hence,
without loss of generality, \emph{$P$ is a lattice}.

Let $e_\gx\colon\vv<K_\gx,\gl_\gx>\hookrightarrow\vv<Q_\gx,\gn_\gx>$,
$f_\gx\colon\vv<K_\gx,\gl_\gx>\to\vv<P,\gf>$, for $\gx<\gk$, enumerate, up
to isomorphism, all embeddings
$e\colon\vv<K,\gl>\hookrightarrow\vv<Q,\gn>$ and homomorphisms
$f\colon\vv<K,\gl>\to\vv<P,\gf>$ with $\vv<K,\gl>$ and $\vv<Q,\gn>$ finite,
proper \MPL{D}s, and with $K$ a lattice. It is easy to construct, by using
Proposition~\ref{P:SquareAmalg} and Remark~\ref{Rk:EmbD}, a transfinite chain
$\famm{\vv<L_\gx,\gf_\gx>}{\gx<\gk}$ of proper \ML{D}s, together with
embeddings
$f_{\gx,\gh}\colon\vv<L_\gx,\gf_\gx>\hookrightarrow\vv<L_\gh,\gf_\gh>$, for
$\gx<\gh<\gk$, satisfying the following properties:

\begin{enumerate}
\item $\vv<L_0,\gf_0>=\vv<P,\gf>$;

\item $f_{\gx,\gz}=f_{\gh,\gz}\circ f_{\gx,\gh}$, for $\gx<\gh<\gz<\gk$;

\item for any $\gx<\gk$, there exists a homomorphism
$g_\gx\colon\vv<Q_\gx,\gn_\gx>\to\vv<L_{\gx+1},\gf_{\gx+1}>$ such that
$f_{0,\gx+1}\circ f_\gx=g_\gx\circ e_\gx$, that is, the following
diagram is commutative:
 \[
 \begin{CD}
 \vv<Q_\gx,\gn_\gx> @>g_\gx>> \vv<L_{\gx+1},\gf_{\gx+1}>\\
 @Ae_{\gx}AA @AAf_{\gx,\gx+1}A\\
 \vv<K_\gx,\gl_\gx> @>>f_{0,\gx}\circ f_\gx> \vv<L_\gx,\gf_\gx>
 \end{CD}
 \]
\end{enumerate}

We denote by $\vv<P,\gf>'$ the direct limit of all $\vv<L_\gx,\gf_\gx>$,
with transition maps $f_{\gx,\gh}$, for $\gx<\gh<\gk$.
Let $f_{\vv<P,\gf>}\colon\vv<P,\gf>\hookrightarrow\vv<P,\gf>'$ be the
limiting map associated with the direct system above. Observe that
$f_{\vv<P,\gf>}$ is an embedding.

The \ML{D}
$\vv<P,\gf>'$ has the property that for every embedding\linebreak
$e\colon\vv<K,\gl>\hookrightarrow\vv<Q,\gn>$ of finite proper \MPL{D}s, with
$K$ a \emph{lattice}, and every homomorphism
$f\colon\vv<K,\gl>\to\vv<P,\gf>$, there exists a homomorphism
$g\colon\vv<Q,\gn>\to\vv<P,\gf>'$ such that
$g\circ e=f_{\vv<P,\gf>}\circ f$.

To conclude the proof, it suffices to iterate the process $\go$ times: put
$\vv<P^{(0)},\gf^{(0)}>=\vv<P,\gf>$, and, for $n<\go$, put
$\vv<P^{(n+1)},\gf^{(n+1)}>=\vv<P^{(n)},\gf^{(n)}>'$, with the embedding
$f_{\vv<P^{(n)},\gf^{(n)}>}\colon\vv<P^{(n)},\gf^{(n)}>
\hookrightarrow\vv<P^{(n+1)},\gf^{(n+1)}>$. The direct limit $\vv<L,\gf>$ of
all the $\vv<P^{(n)},\gf^{(n)}>$, with respect to the transition maps
$f_{\vv<P^{(n)},\gf^{(n)}>}$, satisfies the required conditions.
\end{proof}

\section{Proofs of Theorems C and D}

We first recall the statement of Theorem~C:

\begin{all}{Theorem C}
Let $K$ be a lattice, let $D$ be a distributive lattice with zero, and let
$\gf\colon\Conc K\to D$ be a \jzh. There are a relatively complemented
lattice $L$ of cardinality $|K|+|D|+\aleph_0$, a lattice homomorphism
$f\colon K\to L$, and an isomorphism $\ga\colon\Conc L\to D$ such that the
following assertions hold:
\begin{enumerate}
\item $\gf=\ga\circ\Conc f$.

\item The range of~$f$ is coinitial (resp., cofinal) in $L$.

\item If the range of~$\gf$ is cofinal in $D$, then the range of~$f$ is
internal in $L$.

\end{enumerate}
\end{all}

In this section, we shall fix a distributive lattice $D$ with zero and an
internally saturated \ML{D} $\vv<L,\gf>$.

\begin{lemma}\label{L:RelCpl}
The lattice $L$ is relatively complemented.
\end{lemma}

\begin{proof}
Let $a<b<c$ in $L$, we prove that there exists $x\in L$ such that $a=b\mm x$
and $c=b\jj x$.

Put $K=\set{a,b,c}$, the three-element chain, let
$f\colon K\hookrightarrow L$ be the inclusion map. If we put
$\gl=\gf\circ\Conc f$, then $\vv<K,\gl>$ is a finite, proper
(see Definition~\ref{D:EqnSat}) \ML{D} and
$f$ is an embedding from $\vv<K,\gl>$ into $\vv<L,\gf>$.

Next, we put $P=\set{a,b,c,t}$, the two-atom Boolean lattice, with zero
element $a$, unit element $c$, and atoms $b$ and $t$, endowed with the
homomorphism $\gm\colon\Conc P\to D$ defined by
 \begin{align*}
 \gm\gQ_P(a,b)=\gm\gQ_P(t,c)&=\gf\gQ_L(a,b),\\
 \gm\gQ_P(a,t)=\gm\gQ_P(b,c)&=\gf\gQ_L(b,c).
 \end{align*}
Then $\vv<P,\gm>$ is a proper \ML{D}, and
the inclusion map $j\colon K\hookrightarrow P$ is an embedding from
$\vv<K,\gl>$ into $\vv<P,\gm>$. The lattices $K$ and $P$ can be visualized
on Figure~\ref{Fi:RelCpl}.

\begin{figure}[htb]
\begin{picture}(100,65)(40,0)
\thicklines
\put(0,0){\circle{6}}
\put(0,30){\circle{6}}
\put(0,60){\circle{6}}
\put(100,0){\circle{6}}
\put(100,30){\circle{6}}
\put(100,60){\circle{6}}
\put(130,30){\circle{6}}

\put(-5,0){\makebox(0,0)[r]{$a$}}
\put(-5,30){\makebox(0,0)[r]{$b$}}
\put(-5,60){\makebox(0,0)[r]{$c$}}
\put(105,0){\makebox(0,0)[l]{$a$}}
\put(105,30){\makebox(0,0)[l]{$b$}}
\put(105,60){\makebox(0,0)[l]{$c$}}
\put(135,30){\makebox(0,0)[l]{$t$}}

\put(-35,20){\makebox(0,0)[r]{\LARGE $K$}}

\put(0,3){\line(0,1){24}}
\put(0,33){\line(0,1){24}}
\put(100,3){\line(0,1){24}}
\put(100,33){\line(0,1){24}}
\put(102.12,2.12){\line(1,1){25.76}}
\put(127.88,32.12){\line(-1,1){25.76}}

\put(165,20){\makebox(0,0)[r]{\LARGE $P$}}

\end{picture}
\caption{Adding a relative complement of~$b$ in $[a,c]$}
\label{Fi:RelCpl}
\end{figure}
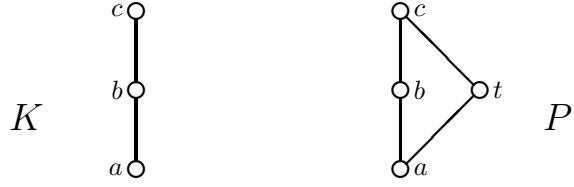

By assumption on $\vv<L,\gf>$, there exists a homomorphism
$g\colon\vv<P,\gm>\to\vv<L,\gf>$ such that $g\circ j=f$.
Put $x=g(t)$. Then $a=b\mm x$ and $c=b\jj x$.
\end{proof}

\begin{definition}\label{D:Persp}
Let $o\leq i$ be elements of a lattice $K$.
We say that the elements $a$, $b$ of the interval
$[o,i]$ are \emph{perspective in $[o,i]$}, if there exists $x\in[o,i]$ such
that $x\mm a=x\mm b$ and $x\jj a=x\jj b$.
\end{definition}

\begin{lemma}\label{L:Persp}
Let $o$, $i$, $a$, $b\in L$ such that $o\leq\set{a,b}\leq i$. Then the
following conditions are equivalent:
\begin{enumerate}
\item $a$ and $b$ are perspective in $[o,i]$.

\item $\gf\gQ_L(o,a)=\gf\gQ_L(o,b)$ and $\gf\gQ_L(a,i)=\gf\gQ_L(b,i)$.

\end{enumerate}
\end{lemma}

\begin{proof}
(i)$\Rightarrow$(ii) If $a$ and $b$ are perspective in $[o,i]$, then the
intervals $[o,a]$ and $[o,b]$ are projective, hence $\gQ_L(o,a)=\gQ_L(o,b)$.
Similarly, $\gQ_L(a,i)=\gQ_L(b,i)$.

(ii)$\Rightarrow$(i)
Let $K=\set{0,u\mm v,u,v,u\jj v,1}$ be the lattice diagrammed on
Figure~\ref{Fi:PictK}, and let $f\colon K\to L$ be the unique lattice
homomorphism sending $0$ to $o$, $1$ to $i$, $u$ to $a$, and $v$ to $b$. We
put $\gl=\gf\circ\Conc f$.

\begin{figure}[htb]
\begin{picture}(100,110)(40,0)
\thicklines
\put(20,0){\circle{6}}
\put(20,30){\circle{6}}
\put(0,50){\circle{6}}
\put(40,50){\circle{6}}
\put(20,70){\circle{6}}
\put(20,100){\circle{6}}

\put(15,0){\makebox(0,0)[r]{$0$}}
\put(15,30){\makebox(0,0)[r]{$u\mm v$}}
\put(-5,50){\makebox(0,0)[r]{$u$}}
\put(45,50){\makebox(0,0)[l]{$v$}}
\put(15,70){\makebox(0,0)[r]{$u\jj v$}}
\put(15,100){\makebox(0,0)[r]{$1$}}

\put(-35,40){\makebox(0,0)[r]{\LARGE $K$}}

\put(20,3){\line(0,1){24}}
\put(17.88,32.12){\line(-1,1){15.76}}
\put(22.12,32.12){\line(1,1){15.76}}
\put(2.12,52.12){\line(1,1){15.76}}
\put(37.88,52.12){\line(-1,1){15.76}}
\put(20,73){\line(0,1){24}}

\put(120,0){\circle{6}}
\put(120,30){\circle{6}}
\put(100,50){\circle{6}}
\put(140,50){\circle{6}}
\put(170,50){\circle{6}}
\put(120,70){\circle{6}}
\put(120,100){\circle{6}}

\put(115,0){\makebox(0,0)[r]{$0$}}
\put(115,30){\makebox(0,0)[r]{$u\mm v$}}
\put(95,50){\makebox(0,0)[r]{$u$}}
\put(145,50){\makebox(0,0)[l]{$v$}}
\put(180,50){\makebox(0,0)[r]{$x$}}
\put(115,70){\makebox(0,0)[r]{$u\jj v$}}
\put(115,100){\makebox(0,0)[r]{$1$}}
\put(170,100){\makebox(0,0)[l]{$x\jj u=x\jj v=1$}}
\put(170,80){\makebox(0,0)[l]{$x\mm u=x\mm v=0$}}

\put(205,40){\makebox(0,0)[r]{\LARGE $P$}}

\put(120,3){\line(0,1){24}}
\put(117.88,32.12){\line(-1,1){15.76}}
\put(122.12,32.12){\line(1,1){15.76}}
\put(102.12,52.12){\line(1,1){15.76}}
\put(137.88,52.12){\line(-1,1){15.76}}
\put(120,73){\line(0,1){24}}
\put(167.88,52.12){\line(-1,1){45.76}}
\put(167.88,47.88){\line(-1,-1){45.76}}

\end{picture}
\caption{Making $u$ and $v$ perspective in $[0,1]$}
\label{Fi:PictK}
\end{figure}
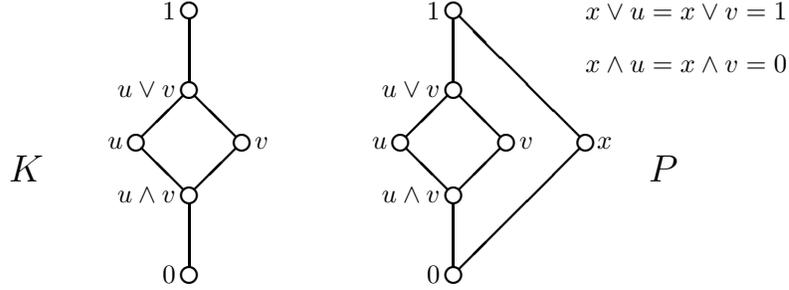

If we could find a finite proper \MPL{D} $\vv<P,\gm>$ and an internal
embedding $j\colon\vv<K,\gl>\hookrightarrow\vv<P,\gm>$ such that $u$ and
$v$ are perspective in $P$, then an argument similar to the
one used in the proof of Lemma~\ref{L:RelCpl} would conclude the proof.

To this end, we simply put $P=K\cup\set{x}$, for an element $x$ not in $K$,
with the ordering of~$K$ extended by the relations
$0<x<1$, together with the following additional joins and meets:
 \begin{equation}\label{Eq:NewCongP}
 x\jj u=x\jj v=1;\qq x\mm u=x\mm v=0,
 \end{equation}
see Figure~\ref{Fi:PictK}.
We denote by $j\colon K\hookrightarrow P$ the canonical embedding. Observe
that $j$ is internal.

We claim that the map $\Conc j$ is surjective.
Indeed, it is easy to verify that the following equalities hold
 \begin{gather*}
 \gQ^+_P(x,u\jj v)=\gQ_P(u\jj v,1);\\
 \gQ^+_P(x,u)=\gQ^+_P(x,v)=\gQ^+_P(x,u\mm v)=\gQ^+_P(x,0)
 =\gQ_P(u,1)=\gQ_P(v,1),
 \end{gather*}
thus all the congruences $\gQ^+_P(x,w)$, for $w\in K$, belong to the range of
$\Conc j$. A similar statement applies to the congruences $\gQ^+_P(w,x)$,
for $w\in K$, which establishes our claim.

We now define congruences $\gx$, $\gh$, $\ga$, and $\gb$ of~$P$ by
 \[
 \gx=\gQ_P(0,u\mm v);\q\gh=\gQ_P(u\jj v,1);\q
 \ga=\gQ^+_P(u,v);\q\gb=\gQ^+_P(v,u).
 \]
It follows from \eqref{Eq:NewCongP} that $\gx\jj\ga=\gx\jj\gb$---denote it
by $\ol{\gx}$, and that $\gh\jj\ga=\gh\jj\gb$---denote it by $\ol{\gh}$.
Therefore, by using the surjectivity of $\Conc j$, we obtain that
 \begin{equation}\label{Eq:ThatsConP}
 \Conc P=\set{\zero_P,\gx,\gh,\ga,\gb,\ga\jj\gb,
 \gx\jj\gh,\ol{\gx},\ol{\gh},\one_P},
 \end{equation}
with all the elements of the right hand side of \eqref{Eq:ThatsConP}
pairwise distinct. The lattice $\Conc P$ is diagrammed on
Figure~\ref{Fi:ConP}.

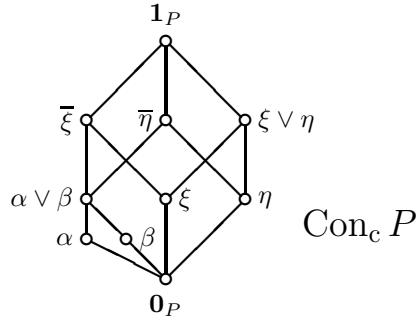
\begin{figure}[htb]
\begin{picture}(100,100)(0,0)
\thicklines
\put(0,0){\circle{4}}
\put(0,30){\circle{4}}
\put(-30,30){\circle{4}}
\put(30,30){\circle{4}}
\put(0,60){\circle{4}}
\put(-30,60){\circle{4}}
\put(30,60){\circle{4}}
\put(0,90){\circle{4}}
\put(-15,15){\circle{4}}
\put(-30,15){\circle{4}}

\put(0,-10){\makebox(0,0){$\zero_P$}}
\put(-35,15){\makebox(0,0)[r]{$\ga$}}
\put(-10,15){\makebox(0,0)[l]{$\gb$}}
\put(-35,30){\makebox(0,0)[r]{$\ga\jj\gb$}}
\put(5,30){\makebox(0,0)[l]{$\gx$}}
\put(35,30){\makebox(0,0)[l]{$\gh$}}
\put(-35,60){\makebox(0,0)[r]{$\ol{\gx}$}}
\put(-5,60){\makebox(0,0)[r]{$\ol{\gh}$}}
\put(35,60){\makebox(0,0)[l]{$\gx\jj\gh$}}
\put(0,100){\makebox(0,0){$\one_P$}}

\put(0,2){\line(0,1){26}}
\put(0,62){\line(0,1){26}}
\put(-30,32){\line(0,1){26}}
\put(30,32){\line(0,1){26}}
\put(-30,17){\line(0,1){11}}
\put(1.41,1.41){\line(1,1){27.17}}
\put(1.41,31.41){\line(1,1){27.17}}
\put(-28.59,31.41){\line(1,1){27.17}}
\put(-28.59,61.41){\line(1,1){27.17}}
\put(-1.41,31.41){\line(-1,1){27.17}}
\put(28.59,31.41){\line(-1,1){27.17}}
\put(28.59,61.41){\line(-1,1){27.17}}
\put(-1.41,1.41){\line(-1,1){12.17}}
\put(-16.41,16.41){\line(-1,1){12.17}}
\put(-1.79,0.89){\line(-2,1){26.42}}

\put(95,20){\makebox(0,0)[r]{\LARGE $\Conc P$}}

\end{picture}
\caption{The congruence lattice of $P$}
\label{Fi:ConP}
\end{figure}

Hence $\Conc P$ is the \jz-semilattice freely generated by
$\check{\gx}$, $\check{\gh}$, $\check{\ga}$, $\check{\gb}$, subject to the
relations
 \begin{equation}\label{Eq:Relgxyab}
 \check{\gx}\jj\check{\ga}=\check{\gx}\jj\check{\gb};\q
 \check{\gh}\jj\check{\ga}=\check{\gh}\jj\check{\gb}.
 \end{equation}
To prove that there exists a \jzh\ $\gm\colon\Conc P\to D$ that satisfies the
equalities
 \begin{multline}\label{Eq:Imgm}
 \gm(\gx)=\gf\gQ_L(o,a\mm b);\q\gm(\gh)=\gf\gQ_L(a\jj b,i);\q\\
 \gm(\ga)=\gf\gQ_L^+(a,b);\q\gm(\gb)=\gf\gQ_L^+(b,a),
 \end{multline}
it suffices to prove that the elements of~$D$ that lie on the right hand
sides of the four equalities in
\eqref{Eq:Imgm} satisfy the relations \eqref{Eq:Relgxyab}, which is
an easy verification. Hence the map $j$ is a homomorphism from $\vv<K,\gl>$
to $\vv<P,\gm>$.
\end{proof}

\begin{notation}
For $o$, $i$, $a$, $b$, $c\in L$ such that $o\leq\set{a,b}\leq i$, we define
$c=a\oplus b$ to hold in $[o,i]$, if $a\mm b=o$ and $a\jj b=c$.
\end{notation}

\begin{lemma}\label{L:DecompEq}
Let $o$, $a$, $b$, $i\in L$ such that $o\leq\set{a,b}\leq i$ .
Then there exist $a_0$, $a_1$, $b_0$, $b_1\in L$ such that the following
conditions hold:
\begin{enumerate}
\item $a=a_0\oplus a_1$ and $b=b_0\oplus b_1$ in $[o,i]$;

\item $\gQ_L(o,a_0)=\gQ_L(o,a_1)=\gQ_L(o,a)$ and
$\gQ_L(o,b_0)=\gQ_L(o,b_1)=\gQ_L(o,b)$;

\item $\gQ_L(a_l\jj b_l,a\jj b)=\gQ_L(o,a\jj b)$, for all $l<2$.
\end{enumerate}
\end{lemma}

\begin{proof}
We put $K=\set{o,a\mm b,a,b,a\jj b,i}$, and we let
$f\colon K\hookrightarrow L$ be the inclusion map. Put $\gl=\gf\circ\Conc f$.
As in the proofs of Lemmas \ref{L:RelCpl} and \ref{L:Persp}, it suffices to
find a finite partial lattice $P$, endowed with a \jzh\
$\gm\colon\Conc P\to D$, an internal embedding
$j\colon\vv<K,\gl>\hookrightarrow\vv<P,\gm>$, and elements $a_0$, $a_1$,
$b_0$, $b_1$ of~$P$ satisfying (i)--(iii) above in $P$.

We use Schmidt's well-known $M_3[K]$ construction, see \cite{Schm68}: namely,
we put
 \[
 P=M_3[K]=\setm{\vv<x,y,z>\in K^3}{x\mm y=x\mm z=y\mm z},
 \]
endowed with the componentwise ordering.
Since $K$ is finite, $P$ is a lattice. Furthermore, the canonical
embedding $j\colon K\hookrightarrow P$, $x\mapsto\vv<x,x,x>$ is
internal and congruence-preserving, see \cite{Schm68} or \cite{GrWe}. Put
$\gm=\gl\circ(\Conc j)^{-1}$. So, $j$ is an internal embedding from
$\vv<K,\gl>$ into $\vv<P,\gm>$.

Now we put $a_0=\vv<a,o,o>$, $a_1=\vv<o,a,o>$, $b_0=\vv<b,o,o>$,
$b_1=\vv<o,b,o>$. Hence $a_0\mm a_1=\vv<o,o,o>$ and $a_0\jj a_1$ is the
least element of~$P$ above $\vv<a,a,o>$, namely, $\vv<a,a,a>$, that is,
$j(a)$. Hence $j(a)=a_0\oplus a_1$. Similarly, $j(b)=b_0\oplus b_1$.
So (i) follows.

For $e=\vv<o,o,i>$, $a_0\oplus e=a_1\oplus e$ and $b_0\oplus e=b_1\oplus e$,
thus $\gQ_P(j(o),a_0)=\gQ_P(j(o),a_1)=\gQ_P(j(o),j(a))$. Similarly,
$\gQ_P(j(o),b_0)=\gQ_P(j(o),b_1)=\gQ_P(j(o),j(b))$.
So (ii) follows.

Finally, $a_0\jj b_0=\vv<a\jj b,o,o>$ and
$a_1\jj b_1=\vv<o,a\jj b,o>$, whence
 \[
 j(a\jj b)=(a_0\jj b_0)\oplus(a_1\jj b_1)\q\text{in }P.
 \]
It follows that
 \[
 \gQ_P(a_l\jj b_l,j(a\jj b))=\gQ_P(j(o),a_{1-l}\jj b_{1-l})
 =\gQ_P(j(o),j(a\jj b)),
 \]
for any $l<2$, so (iii) follows.
\end{proof}

\begin{lemma}\label{L:gfEmb}
Let $o$, $a$, $b$, $i\in L$ such that $o\leq\set{a,b}\leq i$. If
$\gf\gQ_L(o,a)=\gf\gQ_L(o,b)$, then $\gQ_L(o,a)=\gQ_L(o,b)$. More precisely,
there are $a_0$, $a_1$, $b_0$, $b_1\in L$ such that
\begin{enumerate}
\item $a=a_0\oplus a_1\text{ and }b=b_0\oplus b_1$ in $[o,i]$;

\item $a_0$ and $b_0$ (resp., $a_1$ and $b_1$) are perspective in $[o,i]$.
\end{enumerate}
\end{lemma}

\begin{proof}
Let $a_0$, $a_1$, $b_0$, and $b_1$ be as in Lemma~\ref{L:DecompEq}. By (ii)
of Lemma~\ref{L:DecompEq}, $\gQ_L(o,a_l)=\gQ_L(o,a)$ and
$\gQ_L(o,b_l)=\gQ_L(o,b)$, for all $l<2$. It follows from our assumptions
that $\gf\gQ_L(o,a_l)=\gf\gQ_L(o,b_l)$. Furthermore,
 \[
 \gQ_L(a_l,a\jj b)=\gQ_L(a_l,a_l\jj b_l)\jj\gQ_L(a_l\jj b_l,a\jj b)
 =\gQ_L(o,a\jj b),
 \]
for all $l<2$, and, similarly, $\gQ_L(b_l,a\jj b)=\gQ_L(o,a\jj b)$.

It follows then from Lemma~\ref{L:Persp} that $a_0$ and $b_0$ (resp.,
$a_1$ and $b_1$) are perspective in $[o,i]$.
\end{proof}

\begin{lemma}\label{L:gfVemb}
The map $\gf$ is an isomorphism from $\Conc L$ onto an ideal of~$D$.
If, in addition, $\vv<L,\gf>$ is either lower saturated
or upper saturated, then $\gf$ is an isomorphism from $\Conc L$
onto~$D$.
\end{lemma}

\begin{proof}
We first prove that $\gf$ is one-to-one. Let $\ga$, $\gb\in\Conc L$
such that $\gf(\ga)=\gf(\gb)$.
By Lemma~\ref{L:RelCpl}, $L$ is relatively complemented, thus there are
$o$, $a$, $b\in L$ such that $o\leq a$, $o\leq b$, $\ga=\gQ_L(o,a)$, and
$\gb=\gQ_L(o,b)$. In particular, $\gf\gQ_L(o,a)=\gf\gQ_L(o,b)$.
By Lemma~\ref{L:gfEmb}, $\ga=\gb$.

We prove next that the range of~$\gf$ is an ideal of~$D$. Since it is a
\jz-subsemilattice of~$D$, it suffices to prove that the range of~$\gf$ is a
lower subset of~$D$. So let $\ga$ be an element of the lower subset of~$D$
generated by the range of~$\gf$, we prove that $\ga$ belongs to the range of
$\gf$. There are elements $o\leq i$ of~$L$ such that $\ga\leq\gf\gQ_L(o,i)$.
If $\ga=0$ or $\ga=\gf\gQ_L(o,i)$, then $\ga$ belongs to the range of~$\gf$.
Now suppose that $0<\ga<\gf\gQ_L(o,i)$. Put $K=\set{o,i}$, let
$f\colon K\hookrightarrow L$ be the inclusion map, and let
$\gl=\gf\circ\Conc f$. Let $P=\set{o,x,i}$ be the three-element chain, with
$o<x<i$, and let $j\colon K\hookrightarrow P$ be the inclusion map. Endow
$P$ with the \jzh\ $\gm\colon\Conc P\to D$ defined by $\gm\gQ_P(o,x)=\ga$
and $\gm\gQ_P(x,i)=\gf\gQ_L(o,i)$. Observe that $\vv<P,\gm>$ is a proper
\MPL{D} and that $j$ is an internal embedding from $\vv<K,\gl>$ into
$\vv<P,\gm>$. Since $\vv<L,\gf>$ is internally saturated, there exists a
homomorphism $g\colon\vv<P,\gm>\to\vv<L,\gf>$ such that $g\circ j=f$.
Hence the element $\ga=\gm\gQ_P(o,x)=(\gf\circ\Conc g)(\gQ_P(o,x))$
belongs to the range of $\gf$.

Assume, finally, that $\vv<L,\gf>$ is either lower saturated
or upper saturated. Let $\ga\in D$, we prove that $\ga$ belongs to
the range of~$\gf$. We do it, for example, for lower saturated
$\vv<L,\gf>$. The conclusion is obvious if
$\ga=0$, so suppose that $\ga>0$. Pick any element $o$ of~$L$, and put
$K=\set{o}$, endowed with the zero homomorphism from $\Conc K$ to $D$. Let
$P=\set{o,x}$ be the two-element chain, with $o<x$, endowed with the \jzh\
$\gm\colon\Conc P\to D$ defined by $\gm\gQ_P(o,x)=\ga$. Then $j$ is a lower
embedding from $\vv<K,\gl>$ into the proper \MPL{D} $\vv<P,\gm>$, with $\ga$
in the range of~$\gm$. We conclude as in the previous paragraph that $\ga$
belongs to the range of $\gf$.
\end{proof}

We record in Proposition~\ref{P:SummSat} the information that we gathered
in this section about internally saturated \MPL{D}s:

\begin{proposition}\label{P:SummSat}
Let $D$ be a distributive lattice with zero,
let $\vv<L,\gf>$ be an internally saturated \MPL{D}. Then the
following assertions hold:
\begin{enumerate}
\item $L$ is relatively complemented.

\item The map $\gf$ is an isomorphism from $\Conc L$ onto an ideal of~$D$.

\item For $o$, $a$, $b$, $i\in L$ such that $o\leq\set{a,b}\leq i$,
$\gQ_L(o,a)=\gQ_L(o,b)$ if{f} there are $a_0$, $a_1$, $b_0$, $b_1\in[o,i]$
such that the following conditions hold:
\begin{enumerate}
\item $a=a_0\oplus a_1$ and $b=b_0\oplus b_1$ in $[o,i]$.

\item $a_0$ and $b_0$ (resp., $a_1$ and $b_1$) are perspective in $[o,i]$.
\end{enumerate}

\item If, in addition, $\vv<L,\gf>$ is either lower saturated or
upper saturated, then $\gf$ is an isomorphism from $\Conc L$ onto
$D$.
\end{enumerate}
\end{proposition}

Now let $K$, $D$, and $\gf$ be
as in the statement of Theorem~C. We replace $K$ by the image of its kernel
projection (see Lemma~\ref{L:Kernel}), so that without loss of generality,
$\gf$ isolates $0$. By Proposition~\ref{P:EmbSat}, there exist a lower
saturated \MPL{D} $\vv<L,\ga>$ and a lower embedding
$f\colon\vv<K,\gf>\hookrightarrow\vv<L,\ga>$ such that
$|L|=|K|+|D|+\aleph_0$. Since $\vv<L,\ga>$ is lower saturated, it follows from
Proposition~\ref{P:SummSat}(iv) that $\gf$ is an isomorphism from
$\Conc L$ onto $D$. The proof is similar if ``lower'' is replaced by
``upper''.

If the range of~$\gf$ is cofinal in $D$, a similar argument to the one above
works if we replace ``lower embedding'' by ``internal embedding'' and ``lower
saturated'' by ``internally saturated''. This completes the proof of
Theorem~C.

Now we can prove Theorem~D. Indeed, let $K$, $P$, $Q$, $f$, $g$, $\gm$, and
$\gn$ satisfy the assumption of Theorem~D (which is the same as the
assumption of Theorem~B). We first use Theorem B to find a lattice $L'$,
homomorphisms of partial lattices $f'\colon P\to L'$ and $g'\colon Q\to L'$,
and  a \jzh\ $\gf'\colon\Conc L'\to D$ isolating zero such that
$f'\circ f=g'\circ g$, $\gm=\gf'\circ\Conc f'$, $\gn=\gf'\circ\Conc g'$, and
$L'$ is generated, as a lattice, by~$\oll{f}[P]\uu\ol{g}[Q]$.
Then we apply Theorem~C to $\gf'\colon\Conc L'\to D$, to find a relatively
complemented lattice $L$ of cardinality $|L'|+|D|+\aleph_0$, a lattice
homomorphism $h\colon L'\to L$, and an isomorphism $\gf\colon\Conc L\to D$
such that the following assertions hold:
\begin{enumerate}
\item $\gf'=\gf\circ\Conc h$.

\item The range of~$h$ is coinitial (resp., cofinal) in $L$.

\item If the range of~$\gf'$ is cofinal in $D$, then the range of~$h$ is
internal in $L$.
\end{enumerate}
Then $\oll{f}=h\circ f'$ and $\ol{g}=h\circ g'$ satisfy the required
conditions. This completes the proof of Theorem~D.

\section{A few consequences of Theorem~C}\label{S:ConsThmC}

As a special case of Theorem~C (for the case where $\gf$ is an isomorphism),
we obtain the following result:

\begin{corollary}\label{C:EmbRelCpl}
Every lattice $K$ such that $\Conc K$ is a lattice has an internal,
congruence-preserving embedding into a relatively complemented lattice.
\end{corollary}

The other extreme application case of Theorem~C is for $K$ being the trivial
lattice and $\gf$ the zero map:

\begin{corollary}\label{C:ReprDistr}
Let $D$ be a distributive lattice with zero. Then there exists a relatively
complemented lattice $L$ with zero such that $\Conc L\iso D$. Furthermore,
if $D$ is bounded, then one can take $L$ bounded.
\end{corollary}

Actually, by using more of Theorem~C, we can obtain a better representation
result than Corollary~\ref{C:ReprDistr}:

\begin{corollary}\label{C:omegaDistr}
Let $S$ be a distributive \jz-semilattice that can be expressed as the
direct limit of a \emph{countable} sequence of distributive lattices with
zero and \jz-homomorphisms. Then there exists a relatively complemented
lattice $L$ with zero such that $\Conc L\iso S$. If, in addition, $S$ is
bounded, then one can take $L$ bounded.
\end{corollary}

\begin{proof}
We assume that $S$ is the direct limit of~$\famm{D_n}{n<\go}$, with
transition \jzh s $\gf_n\colon D_n\to D_{n+1}$, for $n<\go$. If, in
addition, $S$ is bounded, then we can suppose that the $D_n$-s are bounded
and that the $\gf_n$-s are $\set{\jj,0,1}$-homomorphisms. We construct
by induction a relatively complemented lattice $L_n$, a lattice
homomorphism $f_n\colon L_n\to L_{n+1}$, and an isomorphism
$\ga_n\colon\Conc L_n\to D_n$.

By Corollary~\ref{C:ReprDistr}, there exists a relatively complemented
lattice $L_0$ with zero such that $\Conc L_0\iso D_0$; let
$\ga_0\colon\Conc L_0\to D_0$ be any isomorphism. If $D_0$ has a unit, then
we can suppose that $L_0$ is bounded.

Suppose having constructed a lattice $L_n$ and an isomorphism
$\ga_n\colon\Conc L_n\to D_n$. We apply Theorem~C to the \jzh\
$\gf_n\circ\ga_n\colon\Conc L_n\to D_{n+1}$. We obtain a relatively
complemented lattice $L_{n+1}$, a zero-preserving lattice homomorphism
$f_n\colon L_n\to L_{n+1}$, and an isomorphism
$\ga_{n+1}\colon\Conc L_{n+1}\to D_{n+1}$ such that the following diagram is
commutative.
 \[
 \begin{CD}
 \Conc L_n @>\Conc f_n>> \Conc L_{n+1}\\
 @V\ga_nVV @VV\ga_{n+1}V\\
 D_n @>>\gf_n> D_{n+1}
 \end{CD}
 \]
Furthermore, in case $S$ is bounded, the map $\gf_n\circ\ga_n$ is
cofinal, so we can take $f_n$ with internal range.

Hence the sequence $\famm{L_n}{n<\go}$ of lattices, endowed with the sequence
of transition maps $f_n\colon L_n\to L_{n+1}$, determines a direct limit
system, whose image under the $\Conc$ functor is isomorphic, \emph{via} the
$\ga_n$-s, to the direct system $\famm{D_n}{n<\go}$ with the
\allowbreak$\gf_n$-s. Since the $\Conc$ functor preserves direct limits,
it follows from this that
$\Conc L$ is isomorphic to $S$. In case $S$ is bounded, all the $L_n$-s are
bounded and all the $f_n$-s are $\set{0,1}$-embeddings, thus $L$ is bounded.
\end{proof}

\section{Open problems}

Let $p$ be either a prime number or zero. We denote by $\mathbf{V}_p$ the
quasivariety of all lattices that embed into the subspace lattice of a
vector space over the prime field~$\mathbb{F}_p$ of characteristic $p$.

\begin{problem}\label{Pb:FpExt}
Does every lattice in $\mathbf{V}_p$ have a congruence-preserving
relatively complemented extension in $\mathbf{V}_p$?
\end{problem}

It may be the case that a more natural context for
Problem~\ref{Pb:FpExt} is not provided by the congruence lattice, but
the \emph{dimension monoid}, see \cite{WDim}. The corresponding
reformulation of Problem~\ref{Pb:FpExt} is then the following:

\begin{problem}\label{Pb:FpExtDim}
Does every lattice in $\mathbf{V}_p$ have a dimension-preserving
relatively complemented extension in $\mathbf{V}_p$?
\end{problem}

As in \cite{WDim}, we say that a lattice homomorphism $f\colon K\to L$ is
\emph{dimension preserving}, if the map $\Dim f\colon\Dim K\to\Dim L$ is
an isomorphism.

\begin{problem}
Let $S$ be the \jz-direct limit of a countable sequence of distributive
lattices with zero. Does there exist a relatively complemented lattice $L$
in $\mathbf{V}_p$ such that $\Conc L\iso S$?
\end{problem}

If $K$ is a sublattice of a lattice $L$, we say that $L$ is an
\emph{automorphism-preserving extension} of~$K$, if every automorphism of
$K$ extends to a unique automorphism of~$L$ and $K$ is closed under all
automorphisms of~$L$.

\begin{problem}
Let $K$ be a lattice such that $\Conc K$ is a lattice. Does $K$ have a
relatively complemented, congruence-preserving, automorphism-preserving
extension?
\end{problem}

\section*{Acknowledgments}

It is my great pleasure to thank Kira Adaricheva, Mikhail Sheremet, and
Ji\v{r}\'\i\ T\r{u}ma for having contributed to make this paper alive by
their many suggestions and moral support. Ji\v{r}\'\i\ T\r{u}ma read the
paper quite carefully and suggested a large number of remarks that improved
its readability, and, at one point, corrected an important oversight.

\end{document}